\documentclass{article}

\usepackage[french,english]{babel}
\usepackage{arxiv}
\usepackage[utf8]{inputenc} % allow utf-8 input
\usepackage[T1]{fontenc}    % use 8-bit T1 fonts
\usepackage{url}            % simple URL typesetting
\usepackage{booktabs}       % professional-quality tables
\usepackage{amsfonts}       % blackboard math symbols
\usepackage{nicefrac}       % compact symbols for 1/2, etc.
\usepackage{microtype}      % microtypography
\usepackage{lipsum}
\usepackage{float}
\usepackage{multicol}
\usepackage{algorithm}%
\usepackage{algorithmicx}%
\usepackage{algpseudocode}%

\usepackage{amsmath}
\usepackage{amsthm}
\usepackage{upgreek}

\usepackage{graphics,graphicx}

\usepackage{xcolor}

\usepackage{hyperref}
\newcommand{\N}{\mathbb N}  
 
\newcommand{\K}{\mathbb K}  
\newcommand{\R}{\mathbb R}

\newcommand{\Q}{\mathbb Q}

\newcommand{\T}{\mathbf T}
\newcommand{\J}{\mathbf J}
\newcommand{\cO}{\mathcal O}
\newcommand{\cC}{\mathcal C}
\newcommand{\rd}{\mathrm d}

\newcommand{\Id}{\mbox{\rm Id}}

\newcommand{\ord}{\mbox{\rm ord}}
\newcommand{\rank}{\mbox{\rm rank}}
\newcommand{\diag}{\mbox{\rm diag}}

\newcommand{\im}{\mbox{im}}

\newcommand{\Lie}{\mbox{Lie}}

\newcommand{\Deg}{\hbox{$^{\circ}$}}

\newcommand{\m}{\mbox{m}}
\newcommand{\Newton}{\mbox{N}}
\newcommand{\s}{\mbox{s}}
\newcommand{\rad}{\mbox{rad}}

\newcommand{\dd}{\mathrm{d}}
\newcommand{\rj}{\mathrm{j}}

\newcommand{\kmh}{\mathrm{km}/\mathrm{h}}

\newcommand{\cotan}{\mathrm{cotan}}

\def\lbb{[\![}
\def\rbb{]\!]}

\long\def\PATCH#1{{%\color{blue}
#1}}

\def\meth{\text{\dh}}

\newtheorem{theorem}{Theorem}  
\newtheorem{corollary}[theorem]{Corollary}
\newtheorem{definition}[theorem]{Definition}
\newtheorem{example}[theorem]{Example}
\newtheorem{remark}[theorem]{Remark}
\newtheorem{lemma}[theorem]{Lemma}
\newtheorem{proposition}[theorem]{Proposition}

\def\includegraphicsh[#1]#2#3{{\renewcommand\tabcolsep{0pt}
    \begin{tabular}[b]{c}\href{http://www.lix.polytechnique.fr/~ollivier/AERODYNAMICS/#2}{\includegraphics[#1]{#2}}\\#3\end{tabular}}}

\def\graph[#1]#2#3#4{{\renewcommand\tabcolsep{0pt}
  \begin{tabular}[b]{c}\href{http://www.lix.polytechnique.fr/~ollivier/AERODYNAMICS/#2}{\begin{tabular}[b]{l} {\lower2mm\hbox{\footnotesize $#3$}}\\
\includegraphics[width=#1]{#2}\end{tabular}}\\ #4\lower -2mm\hbox to
    0pt{\hbox to 7.5mm{\hfill}\footnotesize time (s)\hss}\end{tabular}}}

\def\Ast#1{#1\hbox to 0pt{\hbox to 0pt{\hfill}\hskip-0.5pt\small$^{\ast}$\hss}}
\def\AAst#1{#1\hbox to 0pt{\hbox to 0pt{\hfill}\hskip-1.5pt\small$^{\ast}$\hss}}
\def\AAAst#1{#1\hbox to 0pt{\hbox to 0pt{\hfill}\hskip-.75pt\small$^{\ast}$\hss}}

\begin{document}
\title{Flat singularities of chained systems, illustrated with an
  aircraft model}
\author{Yirmeyahu J. Kaminski\\
Holon Institute of Technology\\
Holon, Israel\\
{\tt\small kaminsj@hit.ac.il}
\And 
François Ollivier\\
LIX, CNRS--École Polytechnique \\
91128 Palaiseau Cedex, France\\
{\tt \small ollivier@lix.polytechnique.fr}
}

\maketitle

\begin{abstract} {\small We consider flat differential control systems
    for which there exist flat outputs that 
are part of the state variables and study them using Jacobi bound. We
introduce a notion of saddle Jacobi bound for an ordinary differential
system of $n$ equations in $n+m$ variables. Systems with saddle
Jacobi number equal to $0$ generalize various notions of chained and diagonal
systems and form the widest class of systems admitting subsets of
state variables as flat output, for which flat parametrization may be
computed without differentiating the initial equations. We investigate
apparent and intrinsic flat singularities of such systems.  As an
illustration, we consider the case of a simplified aircraft model,
providing new flat outputs and showing that it is flat at all points
except possibly in stalling conditions. Finally, we present numerical
simulations showing that a feedback using those flat outputs is robust
to perturbations and can also compensate model errors, when using a
more realistic aerodynamic model.

}
\end{abstract}

{\selectlanguage{french}
  
\renewenvironment{abstract}
{
  \centerline
  {\large \bfseries \scshape Résumé}
  \begin{quote}
}
{
  \end{quote}
}
\begin{abstract} {\small Nous considérons des systèmes 
    différentiellement plats pour lesquels il existe des sorties
    plates qui font partie des variables d'état et nous les étudions
    en utilisant la borne de Jacobi. Nous introduisons une notion
    de nombre-selle de Jacobi pour un système pour $n$ équations en
    $n+m$ variables. Les systèmes avec un nombre-selle de Jacobi égal
    à $0$ généralisent diverses notions de systèmes chaînés et
    diagonaux et forment la classe la plus large de systèmes admettant
    des sous-ensembles de variables d'état en tant que sortie plate,
    pour lesquels une paramétrisation plate peut être calculée sans
    différencier les équations initiales. Nous étudions les
    singularités plates apparentes et intrinsèques de ces systèmes.  À
    titre d'illustration, nous considérons le cas d'un modèle d'avion
    simplifié, en fournissant de nouvelles sorties plates et en
    montrant qu'il est plat en tout point, sauf éventuellement en
    situation de décrochage. Enfin, nous présentons des simulations
    numériques montrant qu'un bouclage utilisant ces sorties plates
    est robuste aux perturbations et peut également compenser les
    erreurs de modèle, lors de l'utilisation d'un modèle aérodynamique
    plus réaliste.

  }
  
\end{abstract}

}

\date{}

\noindent AMS classifiaction: 93-10, 93B27, 93D15, 68W30, 12H05, 90C27


\section{Introduction}

\subsection{Mathematical context}

\emph{Differentially flat systems}, introduced by Fliess, L{\'e}vine, Martin and
Rouchon~\cite{FLMR95,FLMR99} are ordinary differential systems
$P_{i}(x_{1}, \cdots, x_{n+m})$, $1\le i\le n$, the solutions of which can be
parametrized in a simple way. Indeed, they admit differential
functions $\zeta_{j}(x)$, $1\le i\le m$, such that for all $1\le i\le
n+m$, $x_{i}=X_{i}(\zeta)$, where $X_{i}$ is a differential
function, i.e. also depending on the derivatives of $\zeta$ up to some finite order. Examples of such systems where considered by
Monge~\cite{Monge1787}, and Monge problem, studied by
Hilbert~\cite{Hilbert1912} and Cartan~\cite{Cartan1914,Cartan1915} is
precisely to test if a differential system satisfies this
property\footnote{Allowing change of independent variable,
  \textit{i.e.} in control change of time, so that Monge problem is
  more precisely to test \textit{orbital flatness}~\cite{FLMR99}}.
Flat systems, for which motion planning and feed-back stabilization
are very easy have proven their importance in nonlinear control.  We
use the theoretical framework of \emph{diffiety
  theory}~\cite{Vinogradov1986,Zharinov1992}, and extend to it the
notion of \emph{defect}, introduced by Fliess \textit{et
  al.}~\cite{FLMR95} in the setting of Ritt's \emph{differential
  algebra}~\cite{Ritt50}. The defect is a nonnegative integer that express
the distance of a system to flatness. It is $0$ iff the system is flat.

\emph{Jacobi's bound}~\cite{Jacobi1,Jacobi2} is an \textit{a priori}
bound on the order of a differential system $P_{i}(x_{1}, \ldots,
x_{n})$, $1\le i\le n$ that is expressed by the \emph{tropical
  determinant} of the order matrix $(a_{i,j})$, with
$a_{i,j}:=\ord_{x_{j}}P_{i}$, which is given by the formula
$\cO_{\Sigma}:=\max_{\sigma\in
  S_{n}}\sum_{i=1}^{n}a_{i,\sigma(i)}$. This bound is conjectural in
the general case, but was proved by Kondrateva \textit{et
  al.}~\cite{Kondratieva2009} \PATCH{for quasi-regular systems, a
  genericity hypothesis} that stands when Jacobi's truncated
determinant does not vanish. Then, the bound is precisely the
order.

\PATCH{Quasi-regular systems include linear systems
  (Ritt~\cite{Ritt1935}). Linear systems with constant coefficients
    were considered by Chrystal~\cite{Chrystal1897} (see also
    Duffin~\cite{Duffin1963}). Such a system is defined by $M(\dd/\dd
    t)X=0$, where $X=(x_{1}, \ldots, x_{n})^{\rm t}$ and $M$ is a
    square matrix of linear operators $m_{i,j}(\dd/\dd t)$, that are
    polynomials of degree $a_{i,j}$ in the derivation $\dd/\dd
    t$. Chrystal shows that the order of the system is the degree of
    the characteristic polynomial of the system, i.e.\ the determinant
    $|M(\lambda)|$. This degree is at most the tropical determinant
    $\cO_{A}$ of the order matrix $A=(a_{i,j})$, the basic idea of
    tropicalization being to replace products by sums and sums by
    max. The truncated determinant is then the coefficient of
    $\lambda^{\cO_{A}}$ in $|M(\lambda)|$.}

Harold Kuhn's \emph{Hungarian method}~\cite{Kuhn2012},
discovered independently, is very close to Jacobi's polynomial time
algorithm for computing the bound and is an important step in the
history of combinatorial optimization (see Burkard \textit{et
  al.}~\cite{Burkard2012}). One may notice that this result was much
probably inspired to Jacobi by \emph{isoperimetrical systems},
satisfied by functions $x_{i}(t)$ such that $\int_{t_{1}}^{t_{2}} L(x)
\dd t$ is extremal.

\subsection{Aims of this paper}
We continue the investigation of intrinsic and apparent singularities
of flat control system~\cite{FLMR_95,FLMR_99,Levine-09,Levine-11},
initiated in our previous
papers~\cite{Kaminski-et-al-2018,Kaminski-et-al-2020}, with a study of
block triangular systems that generalizes \emph{extended chained
  form}~\cite{Gstottner2022} and an application to aircraft
control. We recall that flat systems are systems for which the
trajectory can be parametrized using a finite set of state functions,
called \emph{flat outputs}, and a finite number of their
derivatives. This important notion of nonlinear control simplifies
motion planning, feed-back design and also
optimization~\cite{Oldenburg2002,Ross2002,Flores2006,Faulwasser2011,Beaver2024}.

The class of block triangular system is important in practice as it
includes various notions of ``chained systems'' or triangular
systems~\cite{Li2010,Li2012,Li2016,Kolar2015,Silveira2015}, containing many
classical examples such as robot arms~\cite{Franch2003,Yueksel2016},
cars with many trailers~\cite{Rouchon1993} or discretizations of PDE
flat systems~\cite{Ollivier2001,Strohle2022}. For them, testing
flatness reduces to computing the rank of Jacobian matrices and
finding the flat outputs to a combinatorial problem.

Our goal is to investigate if it is possible to choose flat outputs
among the state functions, and to describe the associated regularity
conditions. This is a main difference with preceding papers on chained
systems that investigated the existence of a change of variables and
a static feedback allowing to reduce to some more restricted class of
chained or triangular form.

We try to enlarge as much as possible the class of systems for which
flatness can be tested by polynomial time combinatorics and
computations of rank of Jacobian matrices. For such systems, that we
call \emph{regular oudephippical\footnote{From the Greek
    \lower
    2.pt\hbox{\includegraphics[scale=1]{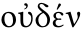}},
    ``nothing'', or ``zero'' for Iamblichus, and
    \lower 4.pt\hbox{\includegraphics[scale=1]{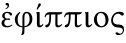}}, ``saddle''.}
  systems}, some lazy flat parametrization may be computed without
differentiating the initial equations, meaning that we have a block
decomposition of variables $\Xi=\bigcup_{i=1}^{r}\Xi_{i}$, such that
$\Xi_{1}=F_{1}(Z)$, where $F_{1}$ is a differential function of the
flat outputs $Z$, then $\Xi_{2}=F_{2}(Z,\Xi_{1})$, where $F_{2}$ is a
differential function of $Z$ and the first block $\Xi_{1}$, \dots,
$\Xi_{i+1}=F_{i+1}(Z,\Xi_{1}, \ldots, \Xi_{i})$, \dots\ On the other
hand, we do not require those systems to be in normal form: they can
be implicit systems and, if we impose that some kind of flat
parametrization can be computed without further differentiation, we
can nevertheless consider from a theoretical standpoint a much larger
class of systems than the usual state space representation, that may
require to be computed a great number of derivation, again bounded by
some Jacobi number.

\subsection{Main theoretical results}

Considering underdetermined systems $P_{i}(x_{1}, \ldots, x_{n+m})$,
we define the saddle Jacobi bound $\hat\cO_{\Sigma}$ as being the minimal
Jacobi bound $\cO(Y,\Sigma)$, for all $Y\subset \Xi$ with $\sharp Y=n$. If
$\hat\cO_{\Sigma}=\cO(\hat Y, \Sigma)$ and the corresponding truncated
determinant $\nabla_{\hat Y, \Sigma}$ does not vanish, then the
\emph{defect} of the system, as defined in~\cite{FLMR95} is at most
$\hat\cO_{\Sigma}$. This implies that if the saddle Jacobi bound is equal
to $0$ and the associated truncated determinant does not identically
vanish, the system is flat, which defines \emph{regular oudephippical systems}.

We give a sufficient condition of flat singularity for some classes of
chained systems, that is enough to prove that the aircraft simplified
model admits an intrinsic flat singularity in some stalling
condition and some sufficient condition of regularity for block
diagonal systems that are enough to show that the simplified aircraft
is flat when not in stalling condition.

We show that previously known classes of chained and diagonal systems
enter the wider class of oudephippical systems. We further prove that a
system $\Sigma$, such that a subset $Z\subset \Xi$ of the state variable is a
flat output and a lazy flat parametrization can be computed without
using any strict derivative of $\Sigma$ is oudephippical.

\subsection{Flat outputs for the aircraft and regularity conditions}

These theoretical results are illustrated with a study of a simplified
aircraft model. With $12$ states, $4$ controls and about $50$
parameters, this model is already more complicated than most flat
models in the literature, although it is among the first to have been
considered. 

Martin~\cite{Martin-PHD,Martin-96} has shown that a simplified
aircraft model where the thrusts related to the actuators and angular
velocities are neglected is flat and given the flat outputs $x$, $y$,
$z$, $\beta$, where $(x,y,z)$ are the coordinates of the center of
gravity and $\beta$ the sideslip angle. We show that the bank angle
$\mu$, the angle of attack $\alpha$ and the engine thrust $F$ can also
be used instead of $\beta$.

We explicit regularity conditions for those choices of flat outputs
and show that the regularity condition for $\mu$ is related to some
kind of stalling condition. The discovery of $3$ new
sets of flat outputs just by a systematic application of our
theoretical results on chained systems illustrate their usefulness.

\subsection{Numerical simulations, models and implementations}
In our simulations, we used the aircraft model and sets of parameters
provided by Grauer and Morelli~\cite{Grauer-Morelli} for various types
of aircrafts: fighter F16C, STOL utility aircraft DHC-6 Twin Otter and
NASA Generic Transport Model (GTM), a subscale airliner model.  Such
aerodynamics models are not known to be flat, unless one neglects some
terms, as Martin did, such as the thrusts created by the control
surfaces (ailerons, elevators, rudder) or related to angular speeds.

We illustrate in two ways the importance of the block
decomposition by providing two implementations, using two different
kinds of feedbacks: a first one in Python is able to reject
perturbations and relies on the difference of
dynamics speeds between the blocks, the second in Maple uses fast
computations allowed by the lazy parametrization to work out a
feed-back able to reject model errors, keeping the values of the flat
outputs close to the planed trajectory, with an acceptable
computational complexity.

We investigated first the robustness of the flat control with respect
to some failures and some perturbations, for the simplified model,
using simulations performed in Python. In a second stage, a Maple
implementation was used to test the ability of a suitable feed-back to
keep the trajectories of the full model close to the theoretical
trajectories computed with the simplified flat one.

We investigate flight situations that are intrinsic singularities for
$\beta$, such as gravity-free flight, for which we use alternative
flat outputs, including bank angle $\mu$. A set of flat outputs
including the thrust $F$ may also be used when $\beta\neq0$ and is
suitable to control a slip-forward maneuver for dead-stick emergency
landing~\cite{NTSB1973,NTSB2010}. See simulations
in~\cite{Ollivier2022}.

\PATCH{The spirit of this study is not at this stage to provide realistic
simulations, but to show that our mathematical methodology is able to
consider models of some complexity, like the GNA, and so to be adapted
to more realistic settings.}

\subsection{Plan of the paper}

We present flat systems in sec.~\ref{sec:flatness}, giving first
definitions and main properties sec.~\ref{sec:flatness}, considering examples
sec.~\ref{subsec:examples_flat_sys} and providing characterization of
flat singularities sec.~\ref{subsec:charact_intrinsic} and
generalizing the notion of defect sec.~\ref{subsec:defect}.

We then present Jacobi's bound sec.~\ref{sec:Jacobi}, starting with
combinatorial definition sec.~\ref{subsec:Jacobi_comb_def} with some
emphasis on Kőnig's theorem sec.~\ref{subsec:Jacobi_Koenig} before
coming to evaluation of the order and computations of normal forms
sec.~\ref{subsec:Jacobi_order_normal}. We then introduce the saddle
Jacobi number sec.~\ref{subsec:saddle}.

We can then define \=o-system sec.~\ref{sec:bock-diag} and provide an
algorithmic criterion sec.~\ref{subsubsec:char_o-sys} with a
sufficient condition
of regularity sec.~\ref{subsec:suf_reg_cond} followed by an
algorithmic criterion for an \=o-system to be regular at a given point
sec.~\ref{subsec:char_reg_o-sys}. We review examples of chained or
triangular systems sec.~\ref{subsec:examples_chained} that enter our
category of \=o-systems and conclude with special sufficient conditions
of regularity or singularity for block triangular systems
sec.~\ref{subsubsec:suf_sing_cond}.

We then consider applications to an aircraft model
sec.~\ref{sec:aircraft}, first describing the equations
sec.~\ref{subsec:model} and the GNA aerodynamic model
sec.~\ref{sub:GNA}. We show that the model is block triangular under
some simplification sec.~\ref{sub:parametrization}. We then
investigate the four main choices of flat outputs
sec.~\ref{subsec:choices-FO} and consider stalling conditions and
their relation to flat singularity sec.~\ref{subsec:stalling}.

The last sections are devoted to simulations using first the
simplified model sec.~\ref{sec:simplified}, then using the full model
sec.~\ref{sec:full_model}. 

\subsection{Notations}

This paper mixes theoretical results from various fields, with
different habits for notations, and an aircraft model coming with
classical notations from aircraft engineering. We tried to use uniform
notations, as long as it did not become an obstacle to readability or
made the access to references too difficult. Regarding diffieties in
some abstract setting, it is convenient to denote the derivation
operators by different symbols, such as $\delta$, $\dd$ or even
$\dd_{t}$ in the case of jet space.

When we start considering control systems, we prefer to use the more
comfortable notations $x'$, $x''$, \dots, $x^{(k)}$. We will also use
$\partial_{x}$ for the derivation $\partial/\partial x$.

Considering control systems, it is natural to denote by $\bar n$ the
number of state variables, which is also the number of state equations
and by $\bar m$ the number of control. When considering abstract systems,
it is easier to denote by $n(=\bar n+\bar m)$ the total number of variables
and by $s(=\bar n)$ the number of equations.

Considering the aircraft equations, we need use then notations that
are common to most textbooks and technical papers : so $x$, $y$, $z$
are space coordinates, $X$, $Y$, $Z$ are the coordinates of the thrust
in the wind referential and not sets of state variables,
$\delta_{\ell}$ is not a derivation but the rudder control etc. To
avoid conflicts, we managed to restrict the notations of previous
sections used in section~\ref{sec:aircraft} to the sets of variables
$\Xi_{h}$.

\section{Flatness}\label{sec:flatness}

For more details on flat systems, we refer to Fliess \textit{et
  al.}~\cite{FLMR_95,FLMR_99} or L{\'e}vine~\cite{Levine-09,Levine-11}.
Roughly speaking, the solutions of flat systems are parametrized by
$m$ differentially independent functions, called flat outputs, and a
finite number of their derivatives. This property, which characterize
them, is specially important for motion planning. We present here flat
systems in the framework of diffiety
theory~\cite{Vinogradov1986,Zharinov1992}.

\subsection{Definitions and properties}\label{subsec:flat-def}

We will be concerned here with systems of the following shape:
\begin{equation}\label{eq:sys}
x_{i}'=f_{i}(x,u,t),\>\hbox{for}\> 1\le i\le \bar n,
\end{equation}
where $x_{1}$, \dots, $x_{\bar n}$ are the state variables and $u_{1}$,
\ldots, $u_{\bar m}$ the controls.

In the sequel, we may sometimes denote $\partial/\partial_{x}$ by
$\partial_{x}$, for short.

\begin{definition}\label{def:diff}
  A diffiety is a $\cC^{\infty}$ manifold $V$ of
  denumerable dimension equipped with a global derivation $\delta$ (that is a vector field), the \emph{Cartan derivation} of the diffiety. The ring
  of functions $\cO(V)$ is the ring of $\cC^{\infty}$ function on $V$
  depending on a \emph{finite} number of coordinates. The topology on
  the diffiety is the coarsest topology that makes coordinate
  functions continuous, i.e. the topology defined by open sets on
  submanifolds of finite dimensions.
\end{definition}

We can give a few example as an illustration.
\begin{example}
  The point $0$ with derivation $\delta:=0$ is considered as a diffiety.
\end{example}

\begin{example}
  The \emph{trivial diffiety} $\T^{m}$ is $\left(\R^{\N}\right)^{m}$
  equipped with the derivation
$$
  \delta:=\sum_{i=1}^{m}\sum_{k\in\N}
  u_{i}^{(k+1)}\partial/\partial_{u_{i}^{(k+1)}}.
$$
\end{example}

\begin{example}
  The \emph{time diffiety} $\R_{\partial_{t}}$ is $\R$ equipped with the
  derivation $\delta_{t}:=\partial/\partial t$.
\end{example}

\begin{definition}\label{def:diff_morph}
  A morphism of diffiety\footnote{Or Lie-B{\"a}cklund transform.}
  $\phi:V_{1}\mapsto V_{2}$ is a smooth map between manifolds such
  that $\phi^{\ast}\circ\delta_{2}=\delta_{1}\circ\phi^{\ast}$, where
  $\phi^{\ast}:\cO(V_{2})\mapsto\cO(V_{1})$ is the dual application,
  defined by $\phi^{\ast}(f) = f \circ \phi$ for $f \in \cO(V_2)$.
\end{definition}

We may illustrate this definition with the next example.

\begin{example}   The product diffiety $\R_{\partial_{t}} \times \T^{m}$ is
  isomorphic to the jet space $\J(\R,\R^{m})$. Indeed, points of this
  jet space can be seen as couples $(t,S)\in\R\times\R\lbb\tau\rbb$, with
\begin{equation}\label{eq:Si}
  S_{i}:=\sum_{k\in\N}\frac{u_{i}^{(k)}(t)}{k!}\tau^{k}.
\end{equation}
  There is a natural action of the derivation $\rd_{t}$ on the
  ring of function on the jet space $\cO(\J(\R,\R^{m}))$ that defines a
  diffiety structure on it. Using $(t,
  u_{1}, u_{1}', \ldots, u_{m}, u_{m}', \ldots)$, as coordinates,
  there is a natural bijection $\phi$ between $\R_{\partial_{t}} \times \T^{m}$ and
  $\J(\R,\R^{m})$, given by:
$\phi(t,u) = (t,S)
$, with $S_{i}$ defined by \eqref{eq:Si}.
The derivation $\rd_{t}$ on the jet space
  is defined by
$$
  \dd_{t}:=\partial_{t}+\sum_{i=1}^{m}\sum_{k\in\N}u_{i}^{(k+1)}\frac{\partial}{\partial
    u_{i}^{(k)}},
$$ so that $\phi$ is
  compatible with the derivations on both diffieties and is a diffiety
  morphism. Moreover, we have
\begin{equation}\label{eq:ddS}
  \dd_{t} S=(\partial_{t}+\partial_{\tau}) S.
\end{equation}
\end{example}

We can now define flat diffieties.

\begin{definition}\label{def:diff_flat}
   A point $\eta$ of a diffiety $V$ is called flat, if
  it admits a neighborhood $O$ that is diffeomorphic by $\phi$ to an open set of
  $\J(\R,\R^{m})$. Let the generators of $\T^{m}$ be $u_{i}$,
  their images by the dual automorphism $z_{i}:=\phi^{\ast}(u_{i})$ are called
  \emph{linearizing outputs} or \emph{flat outputs}. 
  A diffiety $V$ is \emph{flat} if there exists a dense open set
  $W \subset V$ of \emph{flat points}.
  
  A set of such flat outputs defines a Lie-Backl{\"u}nd atlas, as
  defined in~\cite{Kaminski-et-al-2018}.
  
  For a given set of flat outputs $Z$, a point $\eta$ is called
  \emph{singular} related to this set if it is outside the domain of
  definition of the local diffiety diffeomorphism that defines
  it. Such a point in called an \emph{intrinsic flat singularity} if
  none of its neighborhood is isomorphic to an open subspace of
  $\J(\R,\R^{m})$. Otherwise it is called an \emph{apparent
    singularity}.
\end{definition}

\subsection{Examples}\label{subsec:examples_flat_sys}
We illustrate this definition by associating a diffiety to the
system~\eqref{eq:sys} considered above.

\begin{example} Any system~\eqref{eq:sys} defines a
  diffiety $U\times \left(\R^{\N}\right)^{\bar m}$, where
  $U\subset\R^{\bar n+\bar m+1}$ is the domain of definition of the functions
  $f_{i}$, equipped with the Cartan derivation
  \begin{equation}\label{eq:Cartan}
    \dd_{t}:=
    \partial_{t}+\sum_{i=1}^{\bar n}f_{i}(x,u,t)\partial_{x_{i}}
    + \sum_{j=1}^{\bar m}\sum_{k\in\N}u_{j}^{(k+1)}\partial_{u_{j}^{(k)}}
  \end{equation}
  Such a system is a \emph{normal form} defining the diffiety.
\end{example}

\begin{remark} Making
    the abstract terminology of def.~\ref{def:diff_flat} more
    concrete, flatness means that both the state and input variables
    $x_i$, $u_i$ are functions of the flat outputs $z_i$ and a finite number of
    their derivatives on one hand. On the other hand, this also means
    that the $z_i$ are functions of the state and input variables and
    a finite number of their derivatives, and that the differential
    $\dd z_i$ and all their derivatives are linearly independent and
    generate the vector space of differentials
    $\langle \dd x_{i},\> 1\le i\le n; \dd u_{j}^{(k)},\> 1\le j\le m,\> k\in\N\rangle$.
\end{remark}

Flatness may be illustrated by the classical car example.

\begin{example}\label{exam:car}

A very simplified car model is the following, where $\theta$ is
defined modulo $2\pi$:
    
\begin{equation}\label{carsys:eq}
\theta'  =  \frac{\cos(\theta)y'-\sin(\theta)x'}{\ell}
\end{equation}

The state vector is made of the coordinates $(x,y)$ of a point at
distance $\ell$ of the rear axle's
center and of the angle $\theta$ between the car's axis and the
$x$-axis. The controls may be taken to be $u=x'$ and $v=y'$. 
    
One can define different sets of flat outputs depending on the actual
open set, where they are defined, as follows.

\begin{enumerate}
\item Over $U_1 = \{\zeta_{1}' \neq 0\}$, we take $Z_{(1)} =
  \{\zeta_{1}:=x-\ell\cos(\theta),\zeta_{2}:=y-\ell\sin(\theta)\}$ and
  the inverse Lie-B\"acklund transforms given by:
  $\theta=\tan^{-1}(\zeta_{2}'/\zeta_{1}')$
  or $\theta=\pi+\tan^{-1}(\zeta_{2}'/\zeta_{1}')$,
  $x=\zeta_{1}+\ell\cos(\theta)$ and $y=\zeta_{2}+\ell\sin(\theta)$.

\item Over $U_2 = \{\zeta_{2}' \neq 0\}$, we take again $Z_{(1)} =
  \{\zeta_{1}:=x-\ell\cos(\theta),\zeta_{2}:=y-\ell\sin(\theta)\}$ and
  the inverse Lie-B\"acklund transforms given by:
  $\theta=\cotan^{-1}(\zeta_{1}'/\zeta_{2}')$ or
  $\theta=\pi+\cotan^{-1}(\zeta_{1}'/\zeta_{2}')$,
  $x=\zeta_{1}+\ell\cos(\theta)$ and $y=\zeta_{2}+\ell\sin(\theta)$.

\item Over $U_3 = \{\theta' \neq 0\}$, we take
  $Z_{(3)}
  =\{\zeta_{1}=\theta,\zeta_{2}=cos(\theta)y-\sin(\theta)x\}$.
  Here, the inverse Lie-B{\"a}cklund
  transform is given by:
$x=-\sin(\zeta_{1})\zeta_{2}-\cos(\zeta_{1})(\zeta_{2}'-\ell\zeta_{1}')/\zeta_{1}'$,
  $y=\cos(\zeta_{1})\zeta_{2}-\sin(\zeta_{1})(\zeta_{2}'-\ell\zeta_{1}')/\zeta_{1}'$
  and $\theta=\zeta_{1}$.
\end{enumerate}

See~\cite{Kaminski-et-al-2018} for details, using a more realistic model.
\end{example}

\subsection{Characterization of intrinsic
  singularities}\label{subsec:charact_intrinsic}

The linearized tangent system of~\eqref{eq:sys} at a point $(x_0,u_0)$
is classically defined as:
\begin{equation}\label{eq:lin}
\dot{\delta x} = \frac{\partial f}{\partial x}(x_0,u_0,t)\delta x + \frac{\partial f}{\partial u}(x_0,u_0,t) \delta u. 
\end{equation}

We use here a slightly different definition, that allows a more
precise local study.

\begin{definition}\label{def:power-series}
Let $\eta$ be a point of a diffiety $V$. For any function
$g\in\cO(V)$, we denote by $\rj_{\eta}(g)$ the power series
$\sum_{k\in\N}g^{(k)}(\eta)\tau^{k}\in\R\lbb \tau \rbb$.

  For any system $\Sigma$ defined by~\eqref{eq:sys}, we define the
  \emph{linearized system at the point $\eta$}, denoted by
  $\dd_{\eta}\Sigma$, to be
\begin{equation}\label{eq:linearized}
\dd \dot x_{i}=\sum_{i=1}^{n} \rj_{\eta} \left ( \frac{\partial f_{i}(x,u\PATCH{,t})}{\partial
  x_{i}} \right ) \dd x_{i} +
\sum_{j=1}^{m} \rj_{\eta} \left ( \frac{\partial f_{i}(x,u\PATCH{,t})}{\partial
  u_{j}} \right ) \dd u_{j}.
\end{equation}
\end{definition}

There exists a whole algebraic approach to flat systems and their
linear tangent systems. For details, we refer
to~\cite{fliess-all-93}. We will limit ourselves here to mention that
for a flat system, the module generated by the differentials of the
flat outputs is free, as stated by the following theorem, which provides
a necessary condition for local flatness. One may notice that for a
linear system, controllable means flat: from the algebraic standpoint,
the associated \emph{module} (that is the module generated by the
differentials of the flat outputs) is \emph{free}, which just means
that it is generated by a \emph{basis}.

\begin{theorem}\label{th:linearized} At any flat regular point $\eta$,
  the linearized system defines a free module.
\end{theorem}
\begin{proof} If $z$ is a
  flat output, then at any flat point $\eta$, $\dd_{\eta}z$ is a basis
  of the module defined by the linearized system. Indeed, for any
  function $H(z)$ depending on $z$ and its derivatives up to order
  $r$, we have
$$
\dd H(z)=\sum_{i=1}^{m}\sum_{k=0}^{r} \frac{\partial H}{\partial
  z_{i}^{(k)}}\dd z_{i}^{(k)}, 
$$
so that it is a linear combination of derivatives of the $\dd z_{i}$.
\end{proof}

This criterion may be illustrated by the car example.

\begin{example}
The system defined at ex.~\ref{exam:car} is nonflat at all
trajectories such that $x'=y'=\theta' = 0$. In~\cite{Li2012}, the authors
have shown that no flat output depending only on $x$, $y$, $z$ and not
on their derivatives can be regular on such points.

The linearized system is 
\begin{equation}\label{carsys:lin_eq}
  \dd\theta'=\frac{\cos(\theta)\dd y'-\sin(\theta)\dd x
  +(x'\cos(\theta)-y'\sin(\theta))\dd\theta}{\ell}.
\end{equation}

This implies that $(\ell\dd\theta-\sin(\theta)\dd x +\cos(\theta)\dd
y)'= 0$, so the associated module contains a torsion element and is not
flat, according to th.~\ref{th:linearized}. (See
also~\cite{Kaminski-et-al-2018}.)
\end{example}

A general criterion of freeness for modules other power-series would
allow a wider use of the theorem.

\subsection{Defect}\label{subsec:defect}

We propose the following definition to extend the notion of
defect~\cite{FLMR95}\PATCH{, first introduced in differential
  algebra,} to the framework of diffiety theory.

\begin{definition} Any diffiety $V$ defined by a finite set of
  equations, as a subdiffiety of the jet space $\J(\R,\R^{r})$, may
  locally be described, in the neighborhood of some point $\eta$, as
  an open subset of $\R^{n}\times\T^{m}$, for suitable integers $n$
  and $m$, with coordinate functions $x_{i}$, $1\le i\le n$, for
  $\R^{n}$ and $u_{i,k}$, $(i,k)\in\N\times[1,m]$, for
  $\left(\R^{\N}\right)^{m}$, with a derivation defined by
  $x_{i}'=f_{i}(x,u)$ with $f_{i}$ of order $0$ in the $x_{j}$ and of
  arbitrary order in the $u_{h}$.  This is called a
  \emph{representation of $V$ at $\eta$}, and $n$ is the \emph{order
    of the representation}. (See~\cite{Ollivier2007} for more
  details.)

The \emph{defect of $V$ at $\eta$} is the smallest integer $n$ such
that $V$ admits a representation of order $n$ in a neighborhood of $\eta$. 
\end{definition}

\begin{remark} If the diffiety is defined by an explicit normal form
  \begin{equation}\label{eq:norm_form}
    x_{i}^{(r_{i})}=f_{i}(x),\quad 1\le i\le s,
  \end{equation}
  where the $f_{i}$ do no depend of derivatives of the $x_{i}$, $1\le
  i\le s$ with order greater or equal to $r_{i}$, then, reducing to an
  order $1$ system by adding new variables $x_{i,k}$ standing to
  $x_{i}^{(k)}$, for $0\le k<r_{i}$, and new equations
  $x_{i,k}'=x_{i,k+1}$, for $0\le k< r_{i}-1$, we see that the order
  of the representation is
  $\sum_{i=1}^{s}r_{i}$. See~\cite[4.2]{Ollivier2022b} for more
  details.
\end{remark}

It is obvious that if the defect of $V$ at $\eta$ is $0$, then $\eta$
is a flat point of $V$.
\bigskip\bigskip

\section{Jacobi's bound}\label{sec:Jacobi}

Jacobi's bound was introduced by Jacobi in posthumous
manuscripts~\cite{Ollivier2009a,Ollivier2009b}. It is a bound on the
order of a differential system, that is still conjectural in the
general case, but was proved by Kondratieva \textit{et
  al.}~\cite{Kondratieva2009} under regularity hypotheses in the
framework of differential algebra. A proof in the framework of diffiety
theory is available in~\cite{Ollivier2007} and one may find
complete proofs of all main results of Jacobi in~\cite{Ollivier2022},
in the setting of differential algebra. We will refer to this paper
for combinatorial aspects which are the same for differential algebra
and diffiety theory.

If $P_{i}$ is a function of the $x_{j}$ and their derivatives, we
denote by $\ord_{x_{j_{0}}}P_{i}$ the order of $P_{i}$ considered as a
function of $x_{j_{0}}$ and its derivatives, which is the maximal order of derivation at which $x_{j_0}$ appears in $P_i$.
\bigskip

\textbf{Warning.} From now on, all functions will be assumed to be
analytic on their definition domain, so that the order does not depend
on the considered point.

\PATCH{
\subsection{Jacobi's bound and Smith normal forms}\label{subsec:Smith}
It may be usefull to illustrate the tropical nature of Jacobi's bound
in the simple case of linear differential systems, with constant
coefficients. The basic idea of tropical geometry is to reduce the
study of the algebraic equations that define an algebraic variety to
the study of the set of degrees or multidegrees of polynomials in the
associated ideal. We may consider a system $M(\dd/\dd t)X=0$, with
$X=(x_{1}, \ldots, x_{n})^{\rm t}$ and $M$ a square matrix of linear
operators $m_{i,j}(\dd/\dd t)$, i.e.  polynomials of degree $a_{i,j}$
in the derivation $\dd/\dd t$.

Assume that $BMC=\diag(e_{1}, \ldots, e_{n})$, with $B$ and $C$
inversible, is a Smith normal form of $M$, with all $e_{i}$
nonzero. The order of the linear system $MX=0$ is the sum of orders of
the operators $e_{i}(\dd/\dd t)$, that is the sum of the degrees of
the polynomials $e_{i}$ or the degree of the characteristic polynomial
$|\diag(e_{1}(\lambda), \ldots, e_{n}(\lambda))|$. This is also the
degree of the characteristic polynomial
\PATCH{$|M(\lambda)|$}\footnote{\PATCH{In the case of a system
  $X'=AX$, the characteristic polynomial of $A$ is the determinant of
  $M-\lambda\Id$.}}. \PATCH{This idea appears explicitely in the proof
  schetched by Jacobi~\cite{Jacobi1}.}

We may write 
$$
\PATCH{|M(\lambda)|}=\sum_{\sigma\in S_{n}}
\epsilon(\sigma)\prod_{i=1}^{n}m_{i,\sigma(i)}(\lambda). 
$$
The degree of $\prod_{i=1}^{n}m_{i,\sigma(i)}(\lambda)$ is equal to
$\sum_{i=1}^{n}a_{i,j}:=\deg m_{i,j}$. So, the order of the system is
a most
$$
\max_{\sigma\in S_{n}}\sum_{i=1}^{n}a_{i,\sigma(i)}.
$$ This expression is a \emph{tropical determinant}, obtained from the
order matrix $A=(a_{i,j})$ by replacing in the determinant formula
sums by $\max$ and products by sums.\footnote{The general case is more
  complicated, already for time-varying linear systems (see
  Ritt~\cite{Ritt1935}). Then, there exists an analog of Smith normal
  form, due to Jacobson~\cite{Jacobson1937}, but no suitable notion of
  divisors, as factorization in $\R(t)[\dd/\dd t]$ is not
  unique. Indeed, $(\dd/\dd t)^{2}$ is equal to $(\dd/\dd
  t+1/(x+\alpha))(\dd/\dd t-1/(x+\alpha))$, for any
  $\alpha\in\R$~\cite{Chyzak2022}. One must also notice that
  $\diag(\dd/\dd t, (\dd/\dd t)(\dd/\dd t+1))(x_{1},x_{2})^{\rm t}$ is
  a Smith normal form with $\R[\dd/\dd t]$ as the base ring, but not a
  Jacobson normal form with base ring $\R(t)[\dd/\dd t]$, as then the
  quotient module may be generated by a single element
  $x_{1}+tx_{2}$.}  }

\subsection{Combinatorial definitions}\label{subsec:Jacobi_comb_def}
We recall briefly a few basic definitions and properties.

\begin{definition}\label{def:jac_num}
We denote by $S_{s,n}$ the set of injections from $\{1, \ldots, s\}$ to
$\{1, \ldots, n\}$.

  Let $P_{i}$, $1\le i\le s$ be a differential system in $n\ge s$
  variables $x_{j}$. By convention, if $P_{i}$ is free from $x_{j}$
  and its derivatives, \textit{i.e.} if $P_{i}$ does not depend on
  $x_{j}$ and its derivatives, we define
  $\ord_{x_{j}}P_{i}=-\infty$\footnote{This convention, introduced by
    Ritt (See~\cite[\S~4]{Ollivier2022b} for details), is
    known as the \emph{strong bound}. The convention
    $\ord_{x_{j}}P_{i}=0$ is the \emph{weak bound}.}.  
    
    With this convention we define the \emph{order matrix} of $\Sigma$, denoted
  $A_{\Sigma}:=(a_{i,j})$, where $a_{i,j}:=\ord_{x_{j}}P_{i}$. The
  \emph{Jacobi number} of the system $\Sigma$ is:
$$
\cO_{\Sigma}:=\cO_{A_{\Sigma}}:=\max_{\sigma\in S_{s,n}}\sum_{i=1}^{s} a_{i,\sigma(i)},
$$ which, when $s=n$, is the \emph{tropical determinant} of
$A_{\Sigma}$.

Let $Y\subset X:=\{x_{1}, \ldots, x_{n}\}$ be a subset of $s$
variables, then $\cO_{Y,\Sigma}$ denotes the Jacobi number of $\Sigma$,
considered as a system in the variables of $Y$ alone. 
\end{definition}

The tropical determinant may be computed in polynomial time using
Jacobi's algorithm~\cite[2.2]{Ollivier2022b} that relies on the notion
of \emph{canon}, and is equivalent to
Kuhn's Hungarian method~\cite{Kuhn2012} that relies on the notion of
\emph{minimal cover}. We shall now detail these notions. 

\begin{definition}\label{def:canon}
For a $s\times n$ matrix of integers $A$, denoting by $S_{s,n}$ the
set of injections of integer sets $[1,s]\mapsto[1,n]$, a \emph{canon}
is a vector of integers $(\ell_{1}, \ldots, \ell_{s})$, such that
there exists $\sigma_{0}\in S_{s,n}$ that satisfies, for all $j \in \im \sigma_0$,
$a_{\sigma_{0}^{-1}(j),j}+\ell_{\sigma_{0}^{-1}(j)}=\max_{i=1}^{n}(a_{i,j}+\ell_{i})$. The
$a_{i,\sigma(i)}$, for $1 \le i \le s$, are a \emph{maximal family of
  transversal maxima}.
\end{definition}

The following proposition is easy and yet important. 

\begin{proposition}
When $s=n$, if $\ell$ is a canon with a maximal family of transversal
maxima described by the permutation $\sigma_{0}$, then the tropical
determinant of $A$ is $\sum_{i=1}^{n}a_{i,\sigma_{0}(i)}$.
\end{proposition}

\PATCH{
\begin{example}\label{ex:canon}
In order to compute a maximal transversal sum in the left matrix
below, one may add the integers $(1,0,4,2,3)$ to its rows, so that we
get the right matrix, which is a canon: one may find a transversal
family of maximal elements (in their column), the sum of which is maximal.
$$
  \left(
  \begin{array}{rrrrr}
     1&2&7&3&4\\
     10&4&9&3&5\\
     2&3&2&3&0\\
     8&7&5&4&1\\
     1&6&2&4&2
  \end{array}
  \right)
  \begin{array}{c}
    1\\0\\4\\2\\3
  \end{array}
\quad\quad\quad
  \left(
  \begin{array}{rrrrr}
     2&3&8&4&\Ast{5}\\
     10&4&\Ast{9}&3&5\\
     6&7&6&\Ast{7}&4\\
     \Ast{10}&9&7&6&3\\
     4&\Ast{9}&5&7&5
  \end{array}
  \right)
$$
\end{example}
}

\begin{remark}
This proposition is the first of the two main reasons to introduce the
concept of canon. Indeed computing such a canon can be performed in
polynomial time, while computing directly the Jacobi number has
exponential complexity. The other main justification for using canons
will appear below in the context of application to flatness in
proposition~\ref{prop::2ndreason}.
\end{remark} 

\begin{definition}\label{def:cover}
Assuming $s=n$, a \emph{cover} is a couple of integer vectors $\mu,\nu$, such that
$a_{i,j}\le\mu_{i}+\nu_{j}$. A \emph{minimal cover} is a cover such
that the tropical determinant of $A$ satisfies: $\cO_{A}=\sum_{i=1}^{n}(\mu_{i}+\nu_{i})$.
\end{definition}

The next proposition describes the equivalence between minimal covers
and canons.

\begin{proposition}\label{prop:canon_cover}
Assuming $s=n$, to any canon $\ell$, one may associate a minimal cover
$\mu_{i}:=\max_{\kappa=1}^{n}\ell_{\kappa}-\ell_{i}$ and
$\nu_{j}:=\max_{i=1}^{n}(a_{i,j}-\mu_{i})$.

Reciprocally, to any minimal cover $\mu$, $\nu$, one may associate a
canon $\ell_{i}=\max_{h=1}^{n}\mu_{h}-\mu_{i}$.
\end{proposition}
\begin{proof}
  See~\cite[prop.~20]{Ollivier2022b}.
\end{proof}

\PATCH{
\begin{example}\label{ex:cover} Considering the canon of the matrix defined in
  ex.~\ref{ex:canon}, 
 the minimal cover associated to it is $\alpha=(3,4,0,2,1)$,
  $\beta=(6,5,5,3,1)$. The starred terms in the matrix below are
  entries $a_{i,j}$ such that $a_{i,j}=\alpha_{i}+\beta_{j}$ and
  provide a maximal transversal sum. For all entries, one has
  $a_{i,j}\le\alpha_{i}+\beta_{j}$:
  $$\begin{array}{l}
    \begin{array}{rrrrrr}
      \phantom{1}&\phantom{1}6&5&5&3&1
    \end{array}\\
  \left(
  \begin{array}{rrrrr}
     1&2&7&3&\Ast{4}\\
     10&4&\Ast{9}&3&5\\
     2&3&2&\Ast{3}&0\\
     \Ast{8}&7&5&4&1\\
     1&\Ast{6}&2&4&2
  \end{array}
  \right)
  \begin{array}{c}
    3\\4\\0\\2\\1\hbox to 0pt {.\hss}
  \end{array}
  \end{array}
$$
\end{example}
}

The following theorem shows the existence of a unique minimal canon,
that is computed in polynomial time by Jacobi's
algorithm~\cite[alg.~9]{Ollivier2022b}.

\begin{theorem}\label{th:min_canon}
Using the partial order defined by $\ell\le\ell'$ if
$\ell_{i}\le\ell_{i}'$ for all $1\le i\le s$, there exists a unique
minimal canon $\lambda$, that satisfies $\lambda\le\ell$ for any canon
$\ell$. 
\end{theorem}
\begin{proof}
  See~\cite[th.~13]{Ollivier2022b} for more details.
\end{proof}

The minimal cover associated to the minimal canon will be used for
prop.~\ref{prop:maximum_R}. 

\begin{definition}\label{def:Jacobi_cover}
Assuming $s=n$, to this minimal canon, we associate the minimal cover $\alpha$,
$\beta$ that we call \emph{Jacobi's cover}.
\end{definition}

\PATCH{The canon of ex.~\ref{ex:canon} is the minimal canon and so the
cover of ex.~\ref{ex:cover} is Jacobi's cover.}

\subsection{Kőnig's theorem}\label{subsec:Jacobi_Koenig}

Matrices of $0$ and $1$ are a case of special interest that has been
considered by Frobenius~\cite{Frobenius1917} and
Kőnig~\cite{Konig1931,Szarnyas2020}, to whom is due the following theorem.

\begin{theorem}\label{th:Konig} Let $A=(a_{i,j})$ be a $s\times n$
  matrix such that $a_{i,j}\in\{0,1\}$, $(i,j)\in[1,s]\times[1,n]$, then
  $\cO_{A}$ is the smallest integer $r$ such that all nonzero elements in
  $A$ belong to the union of $p$ rows and $r-p$ columns.
\end{theorem}
\begin{proof} See~\cite[th.~17]{Ollivier2022}.
\end{proof}

\begin{example}
For the matrix $A = \left (  \begin{array}{ccc} 1 & 1 & 1 \\ 1 & 0 & 0 \\ 1 & 0 & 0  \end{array} \right )$, we have $p=1$ and $r=2$, since all nonzero entries appear in the union of the first row and the first column. It is also cleat that the tropical determinant of $A$ is $\cO_A = 2 = a_{1,2} + a_{2,1} + a_{3,3}$. 
\end{example}

In sec.~\ref{subsubsec:char_o-sys}, we will be concerned with order
matrices $A$ containing $0$ and $-\infty$ entries. According to Kőnig
theorem~\ref{th:Konig} and changing $0$ to $1$ and $-\infty$ to $0$,
if one may find at most $r$ entries equal $0$ located in all different rows
and columns of $A$, then one may find $p$ rows $R$ and $r-p$ columns
$C$ such that all entries $0$ belong to a row in $R$ or a column in
$C$.

\begin{definition}\label{def:transversal_0} We call such $0$ or any
  family of elements placed in all mutually different rows and
  columns \emph{transversal elements} and $r$ the
  \emph{maximum number of transversal $0$}.
\end{definition}

We can even be more precise with the following result.

\begin{proposition}\label{prop:maximum_R}
  Considering a matrix $A$ the entries of which are either $0$ or
  $-\infty$, there exists a unique set $R_{0}$, maximal for inclusion
  and a unique set $C_{0}$, minimal for inclusion, such that $\sharp
  R_{0}+\sharp C_{0}=r$, where $r$ is the maximal number of
  transversal $0$, and all entries $0$ belong to rows in $R_{0}$
  or columns in $C_{0}$\footnote{Here $\sharp A$ denotes the
    \PATCH{number of elements of
    a set $A$.}}. 
\end{proposition}
\begin{proof} The basic
  idea is to transform the matrix of $0$ and $-\infty$ to a matrix of
  $1$ and $0$ entries. So $r$ is now the number of nonzero entries. 
  
We refer to~\cite[prop.~58]{Ollivier2022b} for details. First, if a
matrix $A$ is a matrix of $0$ and $1$, we may restricts to covers
$\mu$, $\nu$ that are vectors of $0$ and $1$, as well as the
associated canon $\ell$ (see~\cite[prop.~26]{Ollivier2022b}.

  Then, we can make the matrix $A$ a $n\times n$ square matrix by
  adding rows or columns of $0$, which does not change the value of
  $r$. Then, let $\lambda$ be the minimal canon of $A$. Provided that
  some $\lambda_{i_{0}}$ is equal to $1$, the rows in $R_{0}$ are the
  rows of index $i$, with $\lambda_{i}=0$, which is equivalent to
  $\alpha_{i}=1$ for the Jacobi cover
  (def.~\ref{def:Jacobi_cover}). Without loss of generality, as there
  exist $r$ transversal $1$, we may
  assume that $a_{i,i}=1$, for $1\le i\le r$. The column $j$ belongs
  to $C$ iff $j \notin R$, so the minimality of $\lambda$ imply
  that $R_{0}$ is maximal for inclusion and implies the minimality of
  $C_{0}$.

  The result is straightforward when $r=s$. When $r<s$ and all
  $\lambda_{i}$ are $0$, we only have to add a row of $n$ ones
  and a column of $n+1$ zeros to reduce to the previous case. 
\end{proof}

\PATCH{
\begin{example}\label{ex:Rmax} In the following matrix of zeros and
  ones, the maximal sum is $4$ (starred entries), so that all ones
  belong to the union of $p$ rows and $4-p$ colums. The maximal set of
  such rows contains the first, second and third rows, that must be
  completed with the first column. They correspond to the rows with
  $\lambda_{i}=0$ in the minimal canon. All ones are also included in
  the union of the first two rows and the first two columns.
  $$
  \left(
  \begin{array}{rrrrr}
     1&0&1&\AAst{1}&1\\
     1&0&\AAst{1}&1&0\\
     1&\AAst{1}&0&0&0\\
     \AAst{1}&0&0&0&0\\
     1&0&0&0&\Ast{0}
  \end{array}
  \right)
  \begin{array}{c}
    0\\0\\0\\1\\1
  \end{array}\quad\quad\quad
  \left(
  \begin{array}{rrrrr}
     1&0&1&\AAst{1}&1\\
     1&0&\AAst{1}&1&0\\
     1&\AAst{1}&0&0&0\\
     \Ast{2}&1&1&1&1\\
     2&1&1&1&\AAst{1}
  \end{array}
  \right)
$$
\end{example}
}

We will also need to use the \emph{path relation} associated to the
minimal canon, which is the key ingredient of Jacobi's algorithm to
compute a minimal canon and $\cO_{A}$ in polynomial time~\cite{Jacobi1}.

\begin{definition}\label{def:path_relation}
  Let $\ell$ be a canon associated to a square $s\times s$ matrix $A$
  with $\cO_{A}\neq-\infty$. Let $\sigma:[1,s]\mapsto[1,s]$ be
  permutation such that $\cO_{A}=\sum_{i=1}^{s}a_{i,\sigma(i)}$.  We
  say that \emph{there is an elementary path from row $i_{1}$ to row
    $i_{2}$} if $a_{i_{1},\sigma(i_{1})}=a_{i_{2},\sigma(i_{1})}$ (which just expresses the fact the maximal element $a_{i_{1},\sigma(i_{1})}$ in column $\sigma(i_1)$ also appears in row $i_2$). We
    define the \emph{path relation defined by $\ell$} as being the
    reflexive and transitive closure of the elementary path relation.
\end{definition}

The path relation does not depend on the choice of a permutation
$\sigma$, provided that
$\cO_{A}=\sum_{i=1}^{s}a_{i,\sigma(i)}$~\cite[prop.~54]{Ollivier2022b}.
We have the following characterization of the minimal
canon~\cite[lem.~51~i)]{Ollivier2022b}.

\begin{lemma}\label{lem:charact_min}
  A canon $\ell$ is the minimal canon iff for any row $i_{1}$, there is a path
  from it to a row $i_{2}$ with $\ell_{i_{2}}=0$.
\end{lemma}

\PATCH{
\begin{example}\label{ex:path}
We illustrate the path relation with example~\ref{ex:canon}. Entries
in the maximal sums are starred, and entries in other rows equal to
an starred term in the same column are italicized. The starred
term in row~$3$ is equal to the italicized term in row~$5$, so that
there is a path from row~$3$ to row~$5$. In the same way, there is a
path from row~$5$ to row~$4$, row~$4$ to row~$2$ and row~$1$ to
row~$2$, with $\lambda_{2}=0$, so that the canon is minimal.
$$
  \left(
  \begin{array}{rrrrr}
     1&2&7&3&4\\
     10&4&9&3&5\\
     2&3&2&3&0\\
     8&7&5&4&1\\
     1&6&2&4&2
  \end{array}
  \right)
  \begin{array}{c}
    1\\0\\4\\2\\3
  \end{array}
\quad\quad\quad
  \left(
  \begin{array}{rrrrr}
     2&3&8&4&\Ast{5}\\
     \hbox{\it 10}&4&\Ast{9}&3&\hbox{\it 5}\\
     6&7&6&\Ast{7}&4\\
     \Ast{10}&\hbox{\it 9}&7&6&3\\
     4&\Ast{9}&5&\hbox{\it 7}&\hbox{\it 5}
  \end{array}
  \right)
$$
\end{example}
}

We can now conclude these combinatorial preliminaries by stating the
following algorithmic result.

\begin{theorem}\label{th:HK}
Let $A$ be a $s\times n$ matrix of $0$ and $1$, with $\cO_{A}=r$,
there exists an algorithm to construct the sets of rows $R_{0}$ and
columns $C_{0}$ of prop.~\ref{prop:maximum_R} in $O(r^{1/2}sn)$
elementary operations.
\end{theorem}
\begin{proof} See~ \cite[algo.~60]{Ollivier2022b}.

  We sketch here the idea of the proof. Using Hopcroft and
Karp~\cite{HK1973} algorithm, we may build a maximal set of
transversal $1$ in $O(r^{1/2}sn)$ operations (see
also~\cite[\S~3]{Ollivier2022b}). This algorithm does not compute the
minimal canon. Using Jacobi's algorithm~\cite[\S~2.2]{Ollivier2022b},
we need to compute third class rows~\cite[\S~2.2]{Ollivier2022b}, and
increase them by $1$, which is to be done only once, as the minimal
canon only contains $0$ and $1$ entries, as already stated in the
proof of prop.~\ref{prop:maximum_R}. The computation of third class
row can be done in $O(sn)$
operations. See~\cite[algo.~9~e)]{Ollivier2022b} for more details.
\end{proof}

\begin{remark}\label{rem:HK}
We assume that this process is implemented in the procedure
\textsc{HK}.
\end{remark}

\PATCH{
\begin{example} We go back to ex.~\ref{ex:Rmax} In the following
  matrix, Hopcroft and Karp algorithm provides a maximal set of
  transversal ones (starred). 
  Row $5$ belongs to the third class, as it contains no starred
  elements, and row $4$ as there is a path from it to row $5$. This is
  just the third class definition: rows containing no starred elements and all
  the rows from which there is a path to a row in the third class. 
  We then obtain the canon by increasing third class rows by $1$.
  $$
  \left(
  \begin{array}{rrrrr}
     \hbox{\it 1}&0&\hbox{\it 1}&\AAst{1}&1\\
     \hbox{\it 1}&0&\AAst{1}&\hbox{\it 1}&0\\
     \hbox{\it 1}&\AAst{1}&0&0&0\\
     \AAst{1}&0&0&0&0\\
     \hbox{\it 1}&0&0&0&\Ast{0}
  \end{array}
  \right)
  \begin{array}{c}
    0\\0\\0\\1\\1
  \end{array}
$$
\end{example}
}

\subsection{Order and normal forms}\label{subsec:Jacobi_order_normal}

\begin{definition}\label{def:trunc-det}
Let $\Sigma = \{ P_{i}, 1\le i\le n\}$ be a square system in $n$ differential
indeterminates $x_{1}$, \dots, $x_{n}$. The \emph{system determinant}
or \emph{truncated determinant}\footnote{Jacobi named it
  \textit{determinans mancum sive determinans mutilatum} because only
the terms $\partial P_{i}/\partial x_{j}^{(a_{i,j})}$ such that
$a_{i,j}=\alpha_{i}+\beta_{j}$ appear in it.} is
$$
\nabla_{\Sigma}:=\left|\frac{\partial P_{i}}{\partial x_{j}^{\alpha_i  + \beta_j}
}\right|,
$$
\end{definition}
where $(\alpha_i)_i$ and $(\beta_j)_j$ define the Jacobi cover. 

With this definition, we may state the following result, due to Jacobi.

\begin{theorem}\label{th:BJ}
Let $P_{i}$, $1\le i\le n$ be a system in $n$ differential
indeterminates $x_{1}$, \dots, $x_{n}$ that defines a diffiety $V$ in
a neighborhood of a point $\eta\in\J(\R,\R^{m})$.

If $\nabla_{P}$ does not vanish at $\eta$, there exists
$\sigma\in S_{n}$ and an open set $W\ni\eta$ such that the
diffiety admits in $W$ a normal form
$$
x_{j}^{(\alpha_{\sigma^{-1}(j)}+\beta_{j})}=f_{j}(x),
$$
so that the order of the diffiety is $\cO_{\Sigma}$.

This normal form may be computed using derivatives of $P_{i}$ of order
at most $\lambda_{i}$. 
\end{theorem}
\begin{proof}
The results relies on Jacobi's bound and Jacobi's shortest reduction
method~\cite[\S~1 and \S~3]{Jacobi2}. See also~\cite[7.3 or
  9.2]{Ollivier2022b} in the algebraic case. We give a sketch of the
proof in the framework of
diffieties. See~\cite[th.~0.3~(ii)]{Ollivier2007} for more details.

First, let $\{P_{1}, \ldots, P_{n}\}=\bigcup_{k=1}^{q}\Sigma_{k}$ be a
partition of the set $P$ of equations such that
$\lambda_{i_{1}}=\lambda_{i_{2}}$, where $\lambda$ is the minimal
canon, iff $P_{i_{1}}$ and $P_{i_{2}}$ belong to the same subset
$\Sigma_{k}$ of the partition. One may find a permutation $\sigma\in
S_{n}$ such that, for all $1\le h\le q$,
$$
D_{h}:=\left|\frac{\partial P}{\partial x};
(P,x)\in\bigcup_{k=1}^{h}\Sigma_{k}\times\sigma(\bigcup_{k=1}^{h}\Sigma_{k})\right|
$$
does not vanish at point $\eta$.

We recall that $\alpha_{i}=\Lambda(:=\max_{k=1}^{n}\lambda_{i})-\lambda_{i}$
and that $\lambda_{i}=\max_{k=1}^{n}\alpha_{i}(=\Lambda)-\alpha_{i}$ by
def.~\ref{def:Jacobi_cover}. We may
then consider the set of equations 
$E:=\{P_{i}^{(k)}|1\le i\le n,\>0\le k\le\lambda_{i}\}$ and
the set of derivatives
$U:=\{x_{j}^{(k+\beta_{j})}|\alpha_{\sigma^{-1}(j)}\le
k\le\Lambda\}$. Easy computations show that
$$
\Delta:=\left|\frac{\partial Q}{\partial \upsilon};
(Q,\upsilon)\in E\times U\right|=
\prod_{h=1}^{q}D_{h}^{\Lambda-\hat\lambda_{h}},
$$ where $\hat\lambda_{h}=\lambda_{i}$ for $P_{i}\in\Sigma_{h}$. So
$\Delta$ does not vanish in a neighborhood $W$ of $\eta$, where we may
define a local parametrization of the variety defined by the system
$E$, using the implicit function theorem:
$x_{j}^{(\alpha_{\sigma^{-1}(j)}+k+\beta_{j})}=f_{k,j}(x)$, for $1\le
j\le n$ and $0\le k\le \lambda_{\sigma^{-1}(j)}$,
where $f_{k,j}$ only depends on derivatives of $x_{j'}$ of order
smaller than
$\alpha_{\sigma^{-1}(j')}+\beta_{j'}$, for $1\le j'\le n$.

It is then easily seen that the equations
$x_{j}^{(\alpha_{\sigma^{-1}(j)}+\beta_{j})}-f_{k,j}(x)=0$, for $1\le j\le n$, and  and their derivatives locally
define $W\cap V$, so that
$x_{j}^{(\alpha_{\sigma^{-1}(j)}+\beta_{j})}=f_{k,j}(x)$ is a normal form
of the diffiety $W\cap V$.

We may check then that the order of the diffiety is the sum
$\sum_{j=1}^{n}(\alpha_{\sigma^{-1}(j)}+\beta_{j})
=\sum_{i=1}^{n}\alpha_{i}+\beta_{i}=\cO_{P}$.
\end{proof}

When the system determinant vanishes, Jacobi's number provides a
majoration of the order, under genericity hypotheses.

Considering algebraic systems, one may also refer to \cite[6 and
  7.1]{Ollivier2022b}. The basic idea is to use a new ordering on
derivatives, compatible with Jacobi's cover: $\ord^{\rm
  J}_{x_{j}}P:=\ord_{x_{j}}P-\beta_{j}$. The nonvanishing of the
system determinant is then precisely the condition required to get
a normal formal by applying the implicit function theorem.
We illustrate the result with a linear example to make computations
easier.

\begin{example}\label{ex:iso}
  Assume that one wants to minimize or maximize the
  integral
  \begin{equation}\label{eq:iso}
\int_{a}^{b} U(x_{1}(t), \ldots, x_{n}(t)) \dd t.
  \end{equation}

Here we have used a shortened notation, the function $U$ actually also depends on the derivatives of the functions $x_1, \cdots x_n$ up to certain orders. 

The functions $x_{j}$ such that this integral is extremal are
solutions of the isoperimetrical system $\Sigma$ defined by the equations
\begin{equation}\label{eq:iso2}
P_{i}(x):=\sum_{k=0}^{e_{i}}(-1)^{k}\frac{\dd^{k}\frac{\partial U}{\partial x_{i}^{(k)}}}{\dd t^{k}}=0,
\end{equation}
with $e_{i}:=\ord_{x_{i}}U$. We have $a_{i,j}:=\ord_{x_{j}}
P_{i}=e_{i}+e_{j}$, so that the minimal canon is
$\lambda_{i}=\max_{k=1}^{n}e_{k}-e_{i}$ and the Jacobi cover is
$\alpha_{i}=e_{i}-\min_{k=1}^{n}e_{k}$,
$\beta_{j}=e_{j}+\min_{k=1}^{n}e_{k}$. The order is equal to Jacobi's
bound $\cO_{\Sigma} = 2\sum_{i=1}^{n}e_{i}$ when the system determinant
$\nabla_{\Sigma}$, which is equal to the Hessian determinant
$|\partial^{2}U/\partial x_{i}\partial x_{j}|$ does not vanish.  See
\cite[\S~1.2]{Ollivier2022b} for details.
\end{example}

\begin{example}\label{ex:Jac}
Consider the system $P_{1}:=x_{1}+x_{2}'$,
$P_{2}:=x_{1}'-x_{2}''+x_{3}$, $P_{3}:=x_{2}'''+x_{3}'$. We have
$\cO_{P}=3$, $\alpha=(0,1,2)$ and $\beta=(0,1,-1)$. The normal forms
compatible with Jacobi's ordering are $x_{1}=-x_{2}'$,
$x_{2}''=x_{3}/2$, $x_{3}'=0$; $x_{1}=-x_{2}'$, $x_{3}=2x_{2}$,
$x_{2}'''=0$; $x_{2}'=x_{1}$, $x_{1}'=-x_{3}/2$, $x_{3}'=0$ and
$x_{2}'=x_{1}$, $x_{3}=2x_{1}'$, $x_{1}''=0$. 
\end{example}

One may use th.~\ref{th:BJ} for systems $\Sigma$ of $s$ equations in $n>s$
variables, by choosing, when it is possible, a subset $Y\subset X$,
with $\nabla_{Y,\Sigma}\neq0$.

\subsection{The saddle Jacobi number}\label{subsec:saddle}

For systems with less equations than variables, one may define the
\emph{saddle Jacobi number}.

\begin{definition}\label{def:saddle}
  Using the same notations as in def.~\ref{def:jac_num}, we define the
  \emph{saddle Jacobi number} of the system $\Sigma$ as being
$$
\hat\cO_{\Sigma}:=\min_{Y\subset X, \sharp Y=s,\cO(Y,\Sigma)\neq-\infty}\cO_{Y,\Sigma}.
$$
Recall that $\cO_{Y,\Sigma}$ is the tropical determinant of the square system obtained by restricting our attention to the variables in $Y$. 

By convention if $\cO_{Y,\Sigma}=-\infty$ for all $Y$, or if $s>n$, we set
$\hat\cO_{\Sigma}=-\infty$. If $A$ is a matrix with entries in
$\N\cup\{-\infty\}$, we define $\hat\cO_{A}$ accordingly.

Systems $\Sigma$ such that $\hat\cO_{\Sigma}=0$ are called \emph{oudephippical
  systems} or \emph{\=o-systems}. A \=o-system is called
\emph{regular} if there exists $Y\subset X$ such that
$\hat\cO_{\Sigma}=\cO_{Y,\Sigma}$ and $\nabla_{Y,\Sigma}$ does not identically
vanish. It is said to be \emph{regular at point $\eta$} if there
exists $Y\subset X$ such that $\hat\cO_{\Sigma}=\cO_{Y,\Sigma}$ and
$\nabla_{Y,\Sigma}$ does not vanish at $\eta$.
\end{definition}

\begin{proposition}\label{prop:defect}
  If $\hat\cO_{\Sigma}=\cO_{Y,\Sigma}=d$ and
  $\nabla_{Y,\Sigma}(\eta)\neq0$, then the defect of $\Sigma$ at
  $\eta$ is at most $d$.
\end{proposition}
\begin{proof} This is a straightforward consequence of th.~\ref{th:BJ}.
\end{proof}

We do not know an algorithm to compute the saddle Jacobi number faster
than by testing all possible subsets $Y\subset X$, but we will see
that it is possible to test in polynomial time if it is $0$.

\begin{definition}\label{def:lazy}
We say that a system $\Sigma\subset O(\J(\R,\R^{n}))$ of $s$ differential
equations in $n$ variables $x_{1}$, \dots, $x_{n}$ admits a lazy flat
parametrization at $\eta\in\J(\R,\R^{n})$ with flat output $Z$ if
there exists a partition $X=\{x_{1}, \ldots,
x_{n}\}=\bigcup_{h=0}^{r}\Xi_{h}$, with $\Xi_{0}=Z$, and an open
neighborhood $V$ of $\eta$, such that for all $0<h\le r$ and all
$x_{i_{0}}\in\Xi_{h}$, there exists an equation
$x_{i_{0}}-H_{i_{0}}(\Xi_{0}, \ldots, \Xi_{r-1})$, where $H_{i_{0}}$
is a differential function defined on $V$ that belongs to the algebraic
ideal\footnote{The algebraic ideal is a proper subset of the differential ideal.}  generated by $\Sigma$ in $O(V)$.
\end{definition}

\begin{remark}\label{rem:param}
It is easily checked that a system $\Sigma$ admitting a lazy flat
parametrization with flat output $Z=\Xi_{0}$ is flat.

We may indeed rewrite the parametrization $\Xi_{h}=\tilde H_{h}(\Xi_{0},
\ldots, \Xi_{h-1})$, for $1\le h\le r$. So, for $h=1$ we have an expression $\Xi_{1}=\hat H_{1}(Z):=\tilde H_{1}(Z)$. We may then recursively define $\hat H_{h}$, for $2\le h\le r$ by setting $\Xi_{h}=\hat H_{h}(Z):=\tilde H_{h}(Z,\hat H_{1}(Z), \ldots, \hat{H}_{h-1}(Z))$.
\end{remark}

Using th.~\ref{th:BJ}, one is able to bound the order of the
computations required to get a full parametrization.

\begin{proposition}
\label{prop::2ndreason}
In order to compute explicitly the full flat parametrization, we need
to differentiate equation $P_{i}$ at most $\lambda_{i}$ times, if
$\lambda$ is the minimal canon of the order matrix
$A_{\bigcup_{h=1}^{r}\Xi_{h},\Sigma}$.  
Then for a flat output $\zeta\in Z=\Xi_{0}$,
assuming that $\ord_{\zeta}H_{i}=e_{i}$, the maximal order of
$\zeta$ in the full flat parametrization is at most 
$\max_{P_{i}\in\bigcup_{h=1}^{r}\Sigma_{h}}\lambda_{i}+e_{i}$.
\end{proposition}
\begin{proof} This is a straightforward consequence of th.~\ref{th:BJ}.
\end{proof}

\begin{remark}
This result exhibits the second main reason for which the use of a canon
in our context has a very important impact.
\end{remark}

The next example will help to understand the situation.
\begin{example}
Consider the system $P_{1}:=x_{5}-(x_{4}'+x_{3}'')=0$,
$P_{2}:=x_{6}-(x_{4}''+x_{3}'+x_{1}^{(4)})=0$,
$P_{3}:=x_{3}-(x_{1}'+x_{2})=0$, $P_{4}:=x_{4}-(x_{2}'+x_{1}')=0$. We
have then a full flat parametrization, with flat outputs
$Z=\{x_1,x_2\}$, $x_3 = x_{1}'+x_{2}$, $x_4 = x_{2}'+x_{1}'$,
$x_{5}=x_{1}'''+x_{2}''+x_{1}''+ 2x_{2}''$ and $x_{6}=x_{1}^{(4)} +
x_{1}'''+x_{2}'''+x_{1}''+ x_{2}^{(3)} + x_{2}'$, that may be computed
using derivatives of $P_{3}$ up to order $2$ and $P_{4}$ up to order
$2$. The vector $(0,0,2,2)$ is indeed the minimal canon of the order
matrix
$$ A_{\{x_{3},x_{4},x_{5},x_{6}\},\Sigma}=
\left(
\begin{array}{cccc}
  2&1&0&-\infty\\
  1&2&-\infty&0\\
  0&-\infty&-\infty&-\infty\\
  -\infty&0&-\infty&-\infty
\end{array}
\right)
$$
\end{example}

  One may remark that such expressions may be much bigger and much
  harder to compute, as shown by the next example, so that we have
  advantage to achieve numerical computations with lazy
  parametrizations.

\begin{example}\label{ex:size} Consider the system $x_{3}=x_{1}x_{2}$,
  $x_{4}=\left(x_{3}^{(k)}\right)^{d}$. It is a lazy flat
  parametrization with flat outputs $x_{1}$ and $x_{2}$. If we develop
  $\left((x_{1}x_{2})^{(k)}\right)^{d}$, we get a expression with
  ${d+k\choose k}$ monomials.

  Instead of computing the flat parametrization itself, one can
  first choose for all flat outputs $z_{j}$ a function of the time
  $\zeta_{j}(t)$. If the function $\hat H_{r}$ is of order $e_{j}$ in
  $z_{j}$, then, the best is to substitute to $z_{i}$ in $Z$ the sum
  $\rj\zeta_{i}(t) = \sum_{k=0}^{e_{i}}\zeta_{i}^{(k)}(t)\tau^{k}/k!$
  before achieving the substitutions of rem.~\ref{rem:param}. On may
  then use \eqref{eq:ddS} to compute any differential expression.

  So, the size of intermediate results is only proportional to $\max
  e_{i}$, allowing much faster computations. In fact, as we will see
  is \PATCH{subsec.~\ref{subsec:Maple}}, it is enough to have a
  nonvanishing system determinant to work with series, \textit{e.g.}
  by using Newton's method, without actually computing the lazy flat
  parametrization.
\end{example}

We can now conclude this section with the following theorem that
characterizes flat $\bar o$-systems. As a lazy flat parametrization is a
special kind of
regular \=o-system, systems that admit a lazy flat parametrization are
equivalent to a regular \=o-system using simple elimination tools,
such as Gr{\"o}bner bases or characteristic set computations, without
differentiation and without solving PDE systems.

\begin{theorem}\label{th:reg_o_syst_flat}
With the notations of def.~\ref{def:lazy}, we have the
following propositions.

i) A \=o-system $\Sigma$, which is regular at point $\eta$, admits a
\emph{lazy flat parametrization} at point $\eta$.

ii) A system $\Sigma$ that admits a
\emph{lazy flat parametrization} at point $\eta$ with flat output $Z$
and such that $\nabla_{X\setminus Z,\Sigma}(\eta)\neq0$ is a regular
\=o-system at point $\eta$.

iii) If the system $\Sigma$ is a \=o-system, which is regular at point
$\eta$, it is flat at $\eta$.
\end{theorem}
\begin{proof}
i) Let $Y$ be such that $\cO_{Y,\Sigma}=0$ and
$\nabla_{Y,\Sigma}(\eta)\neq0$ and $\lambda$ be the minimal canon of
the order matrix $A_{\Sigma}$ restricted to the columns of $Y$. Then,
there is a partition $\Sigma=\bigcup_{h=1}^{r}\Sigma_{h}$, such that
$\lambda_{i}=\lambda_{i'}$ iff $P_{i}$ and $P_{i'}$ belong to the same
subset $\Sigma_{h}$. We further assume that the sets $\Sigma_{h}$ are
indexed so that the corresponding $\lambda_{i}$ for $P_{i}\in\Xi_{h}$
are decreasing. Let $\sigma\in S_{s,n}$ be such that its values are
the columns defined by $Y$ and such that
$\sum_{i=1}^{s}a_{i,\sigma(i)}=0$. Then we define $\Xi_{0}$ to be
$X\setminus Y$ and $\Xi_{h}$ to be the variables with indexes
$\sigma(i)$, where $i$ runs over the indexes of the equations in
$\Sigma_h$.

As $a_{i,\sigma(i)}+\lambda_{i}\ge a_{i',\sigma(i)}+\lambda_{i'}$ by
the definition of a canon (def.~\ref{def:canon}), the equations of
$\Sigma_{h}$ do not depend on the variables in $\Xi_{h'}$ if $h<h'$,
so that the system is block triangular\footnote{\PATCH{Considering
    the system as a system in the variables of $Y$ only, and the
    remaining variables as parametric variables, we reduce to a system
    of differential dimension $0$ that is indeed block triangular,
    according to the definition in~\cite[4.3]{Ollivier2022b}.}} and
$\nabla_{Y,\Sigma}(\eta)=\prod_{h=1}^{r}D_{h}(\eta)$, where $D_{h}$ is
the Jacobian determinant of $\Sigma_{h}$ with respect to variables in
$\Xi_{h}$, so that for all $1\le h\le r$ $D_{h}(\eta)\neq0$. We only
have to use the implicit function theorem to get the requested lazy
flat parametrization, with flat output $Z:=\Xi_{0}$.

ii) As $\nabla_{X\setminus Z,\Sigma}(\eta)\neq0$, the order of $\Sigma$
considered as a system in the subset of variables $X\setminus Z$ is
equal to $\cO_{X\setminus Z,\Sigma}$ by th.~\ref{th:BJ}. If a lazy flat
parametrization exists with flat output $Z$, then this order must be $0$.

iii) This is a consequence of i) and also a special case of th.~\ref{th:BJ}. 
\end{proof}

This is particularly important for the complexity as the size of a
nonlinear expression grows exponentially with the order of derivation,
as shown by ex.~\ref{ex:size}. This result provides a fast flat
parametrization.

\begin{remark} Any flat parametrization is a \=o-system, which shows
  that flat parametrization can be very far from the usual state space
  representation of control theory. It is known that a flat
  parametrization may be sometimes easier to compute from the physical
  equations.

  \textit{E.g.} the system $x_{1}=x_{4}$, $x_{2}=x_{5}'/x_{4}'$,
  $x_{3}=x_{4}'x_{5}''/x_{4}''-x_{5}$ is a flat parametrization, and
  so a \=o-system, that corresponds to the system of Rouchon
  $x_{3}'=x_{1}'x_{2}'$ with flat outputs $x_{5}=x_{1}$ and
  $x_{6}=x_{1}'x_{2}-x_{3}$. See~\cite{Ollivier1998,Ollivier1999} for
  more details on this classical example.
\end{remark}

\section{\=O-systems and flatness}\label{sec:bock-diag}

In this section, we give efficient criteria to test if a
system is a \=o-system or a regular \=o-system.

\begin{remark} In this section, we consider a submatrix $B$ of a
  matrix $A=(a_{i,j})$ to be defined by a set $R$ of rows and a set
  $C$ of columns, together with the values $a_{i,j}$, for
  $(i,j)\in R\times C$. The empty submatrix corresponds to
  $R=C=\emptyset$.

  By abuse of notation, we identify subsets of equations $P_{i}$ or of
  variables $x_{j}$ with the corresponding subsets of indices.
\end{remark}

\subsection{An algorithmic criterion for
  \=o-systems}\label{subsubsec:char_o-sys} 

In this section, we consider a matrix $A$ of positive integers and
$-\infty$ elements and provide an algorithm to test if $\hat\cO=0$.

We may first make some obvious simplification to spare useless computations.

\begin{remark}\label{rem:obvious}
  If $s>n$ or if $A$ contains a row of $-\infty$
  elements, then $\hat\cO_{A}=-\infty$.
  One may remove from $A$ all columns that contain only $-\infty$
  elements.
\end{remark}

The basic idea of the algorithm relies then on the following lemma.

\begin{lemma}\label{lem:basic_idea}
  Assume that $A$ is a $s\times n$ matrix with $s\le n$ such that
  $\hat\cO_{A}=0$.

  \PATCH{
  i) All the row of $A$ contain at least one element equal to $0$.}

  Let $B$ denote the submatrix formed of the
  columns $C$ of $A$ that contain only $0$ or $-\infty$ entries, $R_{0}$
  its maximal set of rows and $C_{0}$ its corresponding minimal set of
  columns, according to prop.~\ref{prop:maximum_R}.  With these
  hypotheses, we have the following propositions.

  ii) For any subset of columns $Y$ such that $\cO_{Y,A}=0$, let
  $\lambda$ be the minimal canon of $A$ restricted to the columns of
  $Y$. The set $R_{0}$ contains the rows $i$ with $\lambda_{i}=0$, so
  that it is nonempty.

  iii) The matrix $A_{2}$ formed of rows $R_{2}$ not in $R_{0}$ and
  columns $C_{2}$ not
  in $C\setminus C_{0}$ is such that
  there exists a set of columns $Y_{2}$ that satisfies
  $\hat\cO_{A_{2}}=\cO_{Y_{2},A_{2}}=0$.

  iv) The matrix $B'$ formed of rows of $B$ in $R_{0}$ and columns
  not in $C_{0}$ is such that there exists
  $Y_{1}$ with
  $\cO_{Y_{1},B'}=0$, which implies $\cO_{Y_{2}\cup Y_{1},A}=0$.
\end{lemma}
\begin{proof}
  \PATCH{ i) By definition, $\hat\cO_{A}=\cO_{Y,A}=0$ and there exists
    an injection $[1,s]\mapsto Y$ such that $a_{i,\sigma(i)}=0$, for
    $1\le i\le s$.}
  ii) With $\sigma$ as in the proof of i), if $\lambda_{i}=0$, as
  $a_{i,\sigma(i)}+\lambda_{i}\ge a_{i',\sigma(i)}+\lambda_{i'}$
  according to the canon definition, we need have
\begin{equation}\label{eq:canon}
  a_{i',\sigma(i)}=-\infty\> \hbox{if}\> \lambda_{i'}>0\> \hbox{and}\>
  a_{i',\sigma(i)}=0\> \hbox{if}\> \lambda_{i'}=0.
\end{equation}
So, if $\lambda_{i}=0$ the column $\sigma(i)$ belongs to
  the columns of $B$.

  Let $R_{0}':=\{i|\lambda_{i}=0\> \hbox{and}\> i\notin R_{0}\}$. 
Then, by \eqref{eq:canon}, the columns of $\sigma(R_{0}')$ cannot
  contain elements equal to $0$ in rows $i$ with $\lambda_{i}>0$, so
  not in $R_{0}\cup R_{0}'$. This means that all $0$ elements are located in rows
  $R_{0}\cup R_{0}'$ and columns $C_{0}\setminus\sigma(R_{0}')$, which
  contradicts the maximality of $R_{0}$ unless $R_{0}'=\emptyset$, so
  that all rows with $\lambda_{i}=0$ belong to $R_{0}$.

  iii) With the same notations as in the proof of ii), as the columns of
  $C\setminus C_{0}$ contain no element equal to $0$ outside the rows
  of $R_{0}$, $Y_{2}:=\sigma(R_{2})\subset C_{2}$ and
  $\hat\cO_{A_{2}}=\cO_{Y_{2},A_{2}}=0$.

  iv) By the definition of the sets $R_{0}$ and $C_{0}$, one may find
  a family of $r$ transversal elements of $B$
  $a_{i_{k},j_{k}}=0$, for $1\le k\le r$. Let $Y_{1}:=\{j_{k}|
  j_{k}\notin C_{0}\}$, we have $\sharp Y_{1}=\sharp R_{0}$, so that
  $\cO_{Y_{1},B'}=0$. 
\end{proof}

This provides the following recursive algorithm~\ref{algo:o-system}, denoted
\textsc{\=O-test}.  We use here the subroutine \textsc{HK} of
rem.~\ref{rem:HK} that implements the algorithm described in
th.~\ref{th:HK}. We assume moreover that it also returns the set
$Y_{1}$, with the notations of the above lemma.

If $B$ contains only $-\infty$ entries, $R_{0}$ and $Y_{1}$ are
defined to be $\emptyset$ by convention. We denote by $A\setminus R$
the matrix $A$ where the rows in the set $R$ have been suppressed.

\begin{algorithm}
\caption{\textsc{\=O-test}}\label{algo:o-system}

\textbf{Input:} A $s\times n$ matrix $A$ with entries
$a_{i,j}\in\N\cup\{-\infty\}$.

\textbf{Output:} ``failed'' or a set $Y$ of rows such that
$\cO_{Y,A}=0$.

\begin{algorithmic}[0]
\Function{\=O-test}{$A$}
\State\textbf{if}{$s>n$ or \PATCH{$\exists 1\le i\le s\>\forall 1\le j\le
    n>a_{i,j}\neq0$}} \textbf{then return} ``failed''\textbf{fi};
\State Suppress from $A$ all columns of $-\infty$;
\State Build $B$ the submatrix of $A$ of columns containing only $-\infty$ or
$0$ elements.
\State $(R_{0},C_{0},Y_{1}):=$\Call{HK}{$(B)$};
\State\textbf{if}{$R_{0}=\emptyset$} \textbf{then return} ``failed''\textbf{fi};
\State $Y_{2}:=$\Call{\=O-test}{$(A\setminus R_{0})$}; 
\State\textbf{if}{$Y_{2}=$ ``failed''} \textbf{then return} ``failed''\textbf{fi};
\State\textbf{return} $Y_{2}\cup Y_{1}$;
\EndFunction
\end{algorithmic}
\end{algorithm}

\PATCH{
\begin{example}\label{ex:O-test} We illustrate
  algo.~\ref{algo:o-system} 
  with the following example.
{\footnotesize\arraycolsep=1.4pt
  $$
A_{1}=  \left(
  \begin{array}{cccc|c||cc}
    5&2&7&3&-\infty&0&\AAAst{0}\\
    \hline
    \hline
     9&0&0&0&\AAAst{0}&-\infty&-\infty\\
     7&0&0&0&0&-\infty&-\infty\\
     0&-\infty&-\infty&-\infty&-\infty&-\infty&-\infty\\
     -\infty&0&-\infty&-\infty&-\infty&-\infty&-\infty
  \end{array}
  \right)\>\>\>\>
A_{2}=  \left(
  \begin{array}{c|c||ccc}
     9&0&0&0&\AAAst{0}\\
     7&0&0&\AAAst{0}&0\\
    \hline
    \hline
     0&-\infty&-\infty&-\infty&-\infty\\
     -\infty&0&-\infty&-\infty&-\infty
  \end{array}
  \right)\>\>\>\>
A_{3}=  \left(
  \begin{array}{cc}
     \AAAst{0}&-\infty\\
     -\infty&\AAAst{0}
  \end{array}
  \right)
  $$
}
In order to compute \textsc{\=O-test}$(A_{1})$, we need first to
apply \textsc{HK} to $(B_{1})$, which is the submatrix matrix defined by the
last $3$ columns. The zeros in $B_{1}$ belong to the union of 
column $5$ of $A_{1}$ and row $1$. A maximal set of transversal zeros in $B$
corresponds to the starred $0$. So, \textsc{HK}$(B_{1})$ returns the
triplet of sets $(R_{0},C_{0},Y_{1})=(\{1\},\{5\},\{7\})$\footnote{We
  denote for brevity sets of rows or columns by the sets of
  corresponding indices.}.

So, we call \textsc{\=O-test} again on the matrix $A_{1}\setminus
R_{0}$, from which we can supress the two last columns of $-\infty$
elements, which produces $A_{2}$. The matrix $B_{2}$ contains the last
four columns and its zeros are contained in its two first rows and its
first column, the element in a maximal set of transversal zeros being
again starred. So, \textsc{HK}$(B_{2})=(\{1,2\},\{3\},\{4,5\})$.

To conclude, we apply \textsc{\=O-test} on the matrix $A_{2}\setminus
\{1,2\}$, from which we remove columns of $-\infty$, producing
$A_{3}$. Then, the matrix is square with two (starred) transversal
zeros. So \textsc{HK}$(B_{3})=(\{1,2\},\emptyset,\{1,2\})$.

Then \textsc{\=O-test}$(A_{3})=\{1,2\}$,
\textsc{\=O-test}$(A_{2})=\{1,2\}\cup\{4,5\}=\{1,2,4,5\}$ and
\textsc{\=O-test}$(A_{1})=\{1,2,4,5\}\cup\{7\}=\{1,2,4,5,7\}$.

If we had set $a_{4,1}=-\infty$, $A_{3}$ would have contained a full
row of $-\infty$, so that
\textsc{\=O-test}$(A_{3})=$
\textsc{\=O-test}$(A_{2})=$
\textsc{\=O-test}$(A_{1})=$ ``failed''.
\end{example}
}

In the previous algorithm, finding the columns of $B$ can be made
faster using balanced trees.

\begin{remark}\label{rem:AVL} We may sort the elements in the columns of a 
  $s\times n$ matrix $A$ and store the results in \emph{balanced (or
    AVL) trees} (Adel'son-Vel'skii and Landis~\cite{AVL} or
  Knuth~\cite[sec.~6.2.3]{Knuth1998}), with complexity $ns\ln(s)$,
  which allows to delete an element in a column with cost $\ln(s)$,
  preserving the order, and to get the greatest element or a column
  with cost $\ln(s)$.
\end{remark}

The following theorem provides an evaluation of the complexity.

\begin{theorem}\label{th:complexity_o-syst} i) This algorithm tests if
  $\hat\cO_{A}=0$ and if yes returns 
  $Y$ such that $\hat\cO_{A}=\cO_{Y,A}$.

  ii) It works in $O(d^{1/2}psn)$
  elementary operations, where $d$ is the maximal number of
  transversal $0$ in the sets $B$ built at each call of the algorithm
  \textsc{HK} and $p$ the number of recursive calls of
  the main algorithm \textsc{\=O-test}.

  iii) If at each step $C_{0}=\emptyset$, then the algorithm works in
  $O((d^{1/2}+\ln(s))sn)$ operations, using balanced trees as in
  rem.~\ref{rem:AVL}.
\end{theorem}
\begin{proof}
i) The algorithms produces the correct result as a consequence of
lem.~\ref{lem:basic_idea}. Indeed, by
lem.~\ref{lem:basic_idea}~ii) if \textsc{HK}$(B)$ returns $R_{0}=\emptyset$,
we know that $\hat\cO_{A}\neq0$. In the same way if
\textsc{\=O-test}$(A\setminus R_{0})$ returns ``failed'',
$\hat\cO_{A}$ is not $0$ by
lem.~\ref{lem:basic_idea}~iii). Furthermore, except in those two
cases, we know that $\cO_{Y_{2}\cup Y_{1}',A}$ is $0$ by
lem.~\ref{lem:basic_idea}~iv).

ii) For the complexity, at each call the computation of the set $B$
requires $O(sn)$ operations, which is proportional to the size of the
matrix. Then, computing a maximal set of diagonal $0$ requires at most
$O(d^{1/2}sn)$ operations, using Hopcroft and Karp~\cite{HK1973}
algorithm, by th.~\ref{th:HK}. So the total cost is $O(d^{1/2}psn)$,
where $p$ is the number of recursive call to \textsc{\=O-test}.

  iii) Using balanced trees, the cost of the first inspection of the
  matrix becomes $O(\ln(s)sn)$. Let $b_{i}$ and $c_{i}$ respectively
  denote the cardinal of $R_{0}$ and the number of columns in $B$ at
  step~$i$ and let $p$ be the total number of steps. At step
  $i$, we just
  have to remove at most $b_{i} n$ elements from the AVL trees with
  cost $O(\ln(s)b_{i}n)$ and from the matrix with cost
  $O(b_{i}n)$. The detection of columns in $B$ at each step can be made
  with cost $O(\ln(s)n)$. This provides a total cost $O(\ln(s)sn)$.

  The remaining costs lie in the \textsc{HK} routine, for which the
  cost is $O(b_{i}^{1/2}b_{i}c_{i})$ at step $i$ by th.~\ref{th:HK}. So, the
  remaining total cost is
  $O(\sum_{i=1}^{p}b_{i}^{1/2}b_{i}c_{i})\le O(d^{1/2}sn)$, as $d=\max_{i=1}^{p}
  b_{i}$ and $s=\sum_{i}^{p} b_{i}$, which
  concludes the proof.
\end{proof}

More examples will be given in the section~\ref{subsec:examples_chained}.

\subsection{Sufficient condition for regularity}\label{subsec:suf_reg_cond}

The next theorem is an easy sufficient criterion for the
regularity of \=o-systems.

\begin{theorem}\label{th:suff_reg_cond}
  Let $\Sigma$ be a \=o-system defining a diffiety $V$ in some
  neighborhood of a point $\eta$. Using algo.~\ref{algo:o-system}, one
  may assume that we have a partition of the equations
  $\Sigma=\bigcup_{h=1}^{p}\Sigma_{h}$, where $\Sigma_{p-h}$ corresponds to
  rows in $R_{0}$ at step $h$ of the algorithm, and a partition of the
  variables $X=\bigcup_{h=1}^{p}\Xi_{h}$, where $\Xi_{h}$ corresponds
  to columns in $B$ and not in $C_{0}$ at step $h$.

  With these notations, if for all $1\le h\le p$,
  \begin{equation}\label{eq:suf_reg_cond}
      \rank\,\left(\frac{\partial P_{i}}{\partial x}(\eta)
      \mid P_{i}\in\Sigma_{h},\ x\in\Xi_{h}\right)=\sharp\Sigma_{i},
  \end{equation}
then $\Xi$ is \=o-regular at point $\eta$.
\end{theorem}
\begin{proof} Eq.~\ref{eq:suf_reg_cond} implies that for all $1\le
  i\le p$ there exists a subset $Y_{h}\subset\Xi_{h}$ such that
  $\nabla_{Y_{h},\Sigma_{h}}(\eta)\neq0$. By construction, setting
  $Y:=\bigcup_{h=1}^{p}Y_{h}$, we have then
$$
\cO_{Y_{h},\Sigma_{h}}=0,
$$
 so that
$$
\cO_{Y,\Sigma}(\eta)=\sum_{h=1}^{p}\cO_{Y_{h},\Sigma_{h}}=0,
$$
as the equations in $\Sigma_{h}$ do not depend on the variables in
$\Xi_{\ell}$ for $\ell>h$. Furthermore, we have
$$
\nabla_{Y,\Sigma}(\eta)=\prod_{h=1}^{p}\nabla_{Y_{h},\Sigma_{h}}(\eta)\neq0,
$$
so that $\Sigma$ is \=o-regular at point $\eta$.
\end{proof}

\PATCH{So, possible sets of flat outputs in this setting are built by
  choosing $Y_{h}\subset\Xi_{i}$ with $\sharp Y_{h}=\sharp\Sigma_{h}$,
  for $1\le h\le p$, so that $\nabla_{Y_{h},\Sigma_{h}}$ does not
  vanish. This is what will be done on the aircraft example at
  subsec.~\ref{subsec:block_triang} and
  subsec.~\ref{subsec:restricted}.}

This condition is not necessary, as shown by the next example.

\begin{example}\label{ex:not_necessary_reg}
    Let be the system $\Sigma:=\{P_{1}, \ldots, P_{4}\}$ with
    $P_{1}:=x_{1}+x_{3}'+x_{6}^{2}$, $P_{2}=x_{1}+x_{4}'+x_{5}^{3}$,
    $P_{3}:=x_{1}'+x_{3}$, $P_{4}:=x_{2}'+x_{4}$. \PATCH{The order
      matrix of $\Sigma$ is $A_{1}$ below. The matrix $B_{1}$ of rows
      containing only $0$ and $-\infty$ elements correspond to the two
      right columns. It contains two zeros, that are transversal and
      belong to the two upper rows, so that
      \textsc{HK}$(B_{1})=(\{1,2\},\emptyset,\{5,6\})$. Then,
      \textsc{\=O-test} is called on the matrix $A_{2}$.
$$
      A_{1}=\left(
      \begin{array}{cccc||cc}
        0&-\infty&1&-\infty&-\infty&\Ast{0}\\
        0&-\infty&-\infty&1&\Ast{0}&-\infty\\
        \hline\hline
        1&-\infty&0&-\infty&-\infty&-\infty\\
        -\infty&1&-\infty&0&-\infty&-\infty
      \end{array}
      \right)
      A_{2}=\left(
      \begin{array}{cc||cc}
        1&-\infty&\Ast{0}&-\infty\\
        -\infty&1&-\infty&\Ast{0}
      \end{array}
      \right)
      $$
The matrix $B_{2}$ contains the two right columns and we have
\textsc{HK}$(B_{2})=(\{1,2\},\emptyset,\{3,4\})$, so that
\textsc{\=O-test}$(A_{1})=\{3,4,5,6\}$. 
      }
We have thus
    $\Xi_{1}=\{x_{5},x_{6}\}$, $\Xi_{2}=\{x_{3},x_{4}\}$ and
    $\Sigma_{1}=\{P_{1},P_{2}\}$, $\Sigma_{2}=\{P_{3},P_{4}\}$. The
    only set $Y_{i}$ that may be extracted from the $\Xi_{i}$ are the
    $\Xi_{i}$ themselves, as $\sharp \Xi_{i}=\sharp \Sigma_{i}$, for
    $i=1,2$. The corresponding flat outputs are then $x_{1}$
    and $x_{2}$. With this choice,
    $\nabla_{Y,P}=\nabla_{\Xi_{1},\Sigma_{1}}\nabla_{\Xi_{2},\Sigma_{2}}=\nabla_{\Xi_{1},\Sigma_{1}}=3x_{5}^{2}$, 
    which vanishes when $x_{5}=0$.

    But at such a point, we can use the alternative values
    $Y=\{x_{1},x_{3},x_{4},x_{6}\}$, with $\nabla_{Y,\Sigma}=2x_{6}$,
    \PATCH{when $x_{6}\neq0$}.
\end{example}

So we need a more precise criterion of \=o-regularity, that will be
described below.

\subsection{Characterization of regular
  \=o-systems}\label{subsec:char_reg_o-sys}

\PATCH{In this subsection, we will use ex.~\ref{ex:not_necessary_reg}
  and the following example~\ref{ex:p-q-3} as a running examples to
  illustrate all the algorithms.

\begin{example}\label{ex:p-q-3} We will consider the system defined
  by $P=0$, with $P_{1}:=x_{1}^{2}+x_{6}+x_{7}'$,
  $P_{2}:=x_{1}+x_{2}^{2}+x_{5}+x_{6}'$ and $P_{3}:=x_{2}+x_{3}^{2}+x_{4}+x_{5}'$
  at a point where $x_{1}=x_{2}=x_{3}=0$.
\end{example}}

We work at a point $\eta$ of the
diffiety. We assume that the coordinates of this point belong to some
effective subfield $\K$ of $\R$ and that at the point $\eta$ of the
diffiety, for all $P\in\Sigma$ and all $x\in\Xi$, $(\partial
P/\partial x)(\eta)$ can be computed. We neglect the cost of this
computation, the Jacobian matrix being assumed to be given.

If the equations of $\Sigma$ are algebraic,
$\K$ can be $\Q$ or an algebraic extension of $\Q$. We do not consider
here the size of the elements of $\K$ and bit complexity and only
evaluate the number of elementary operations in $\K$, counting in them
operations in $\N$ and all other faster elementary operations.

The basic principle looks much like that of
sec.~\ref{subsubsec:char_o-sys}.  We will need the following easy
technical lemma.

\begin{lemma}\label{lem:vec_sup}
Let $v_{i}\in\R^{n}$ with $v_{i}=(c_{i,1}, \ldots, c_{i,n})$, for
$1\le i\le s$. Let $K$ be the vector space generated by the linear
relations between the vectors $v_i$. Let
$e_{k}=\sum_{i=1}^{s}b_{k,i}v_{i}$, for $1\le k\le r$, be a basis of
$K$ and $R_{\rm [K]}:=\{i\mid \exists k, 1\le k\le r\>,
b_{k,i}\neq0\}$.

The set $R_{\rm [K]}$ does not depend on the choice of the basis $e$.
\end{lemma}
\begin{proof} Let $\bar e$, with
  $\bar e_{k}=\sum_{\ell=1}^{r}\gamma_{k,\ell}e_{\ell}$ be another
  basis, defining another set of rows $\bar R_{\rm [K]}$. If $\bar
  b_{k,i}=\sum_{\ell=1}^{r}\gamma_{k,\ell}b_{\ell,i}$ is nonzero,
  then $b_{\ell,i}$ is nonzero for some $1\le\ell\le r$. So $\bar
  R_{\rm [K]}\subset R_{\rm [K]}$. We prove in the same way the
  reciprocal inclusion.
\end{proof}

The following lemma will allow us to design a recursive process.

\begin{lemma}\label{lem:basic_idea_reg} Assume that a system $\Sigma$ of
  $s$ equations
  in the variables $X:=\{x_{1}, \ldots, x_{n}\}$ is \=o-regular at
  point $\eta$, with $\nabla_{Y,\Sigma}(\eta)\neq0$. We use the same
  notations and hypotheses as in lemma~\ref{lem:basic_idea}, with $A:=A_{\Sigma}$,
  the order matrix of $\Sigma$.

  We consider the minimal canon $\lambda$ of the order matrix $A_{Y,\Sigma}$
  restricted to the columns in $Y$ and a bijection $\sigma:[1,s]\mapsto
  Y$ such that $\sum_{i=1}^{s}a_{i,\sigma(i)}=0$.

  Let $B$ be defined as in lem.~\ref{lem:basic_idea} and $B'$ be a
  submatrix of $B$ restricted to a set of rows $R$ that contains
  $\bar R:=\{i\mid\lambda_{i}=0\}$ and a set of columns $C$ that contains
  $\sigma(\bar R)$.

  i) Let $R_{0}$ and $C_{0}$ be respectively sets of
  rows and columns that contain all entries of $B'$ equal to $0$, with
  $\sharp R_{0}+\sharp C_{0}=\cO_{B'}$ and $R_{0}$ maximal (as in
  th.~\ref{th:Konig}), then $\bar R\subset R_{0}$ and
  $\sigma(\bar R)\subset C\setminus C_{0}$.

  ii) Let $J$ be the value of the Jacobian matrix $(\partial
  P_{i}/\partial x_{j}\mid a_{i,j}\in B')$ at point $\eta$ and
  $R_{[K]}$ be the set of rows associated to a basis of linear
  relations between the rows of $J$ as in lem.~\ref{lem:vec_sup}. We
  also define the set of columns $C_{\rm [K]}:=\{j\mid \exists i\in
  R_{\rm [K]}\> a_{i,j}=0$. Then, $\bar R\subset R\setminus R_{\rm
    [K]}$ and $\sigma(\bar R)\subset C\setminus C_{\rm [K]}$.
\end{lemma}
\begin{proof}
  i) We proceed as in the proof of lem.~\ref{lem:basic_idea} ii).  Let
  $R_{0}':=\{i|\lambda_{i}=0\> \hbox{and}\> i\notin R_{0}\}$.  Then,
  by \eqref{eq:canon}, the columns of $\sigma(R_{0}')$ cannot contain
  elements equal to $0$ in rows $i$ with $\lambda_{i}>0$, which means
  that all $0$ elements are located in rows $R_{0}\cup R_{0}'$ and
  columns $C_{0}\setminus\sigma(R_{0})$, which contradicts the
  maximality of $R_{0}$ unless $R_{0}'=\emptyset$, so that all rows
  with $\lambda_{i}=0$ belong to $R_{0}$.

  ii) Assume that $\bar R\cap R_{\rm [K]}$ is nonempty and contains
  row $i$. Then, by \eqref{eq:canon}, the columns of $\sigma(\bar R)$
  cannot contain entries equal to $0$ and not located in the rows of
  $\bar R$. So, $\nabla_{\sigma(\bar R),\bar \Sigma}$, where $\bar
  \Sigma$ is the subset of equations that corresponds to the rows of
  $\bar R$, is a factor of $\nabla_{Y,\Sigma}$. As the rows of $\bar
  R\cap K_{\rm K}$ in $J$ are involved in a nontrivial linear
  relation, their restrictions to the columns of $\sigma(\bar R)$ are
  linearly dependent and $\nabla_{\sigma(\bar R),\bar \Sigma}$ must vanish at
  $\eta$, which contradicts $\nabla_{Y,\Sigma}(\eta)\neq0$. So $\bar
  R\cap R_{0}=\emptyset$.

  Using again \eqref{eq:canon}, the columns of $\sigma(\bar
  R)$ can only contain entries equal to $-\infty$ in the rows of $R_{\rm
    [K]}$, so that $\sigma(\bar R)\cap C_{\rm [K]}$ is equal to $\emptyset$.
\end{proof} 

\begin{definition} With the notations of
  lem.~\ref{lem:basic_idea_reg}, we define the \emph{row kernel
    support} of $B'$ at point $\eta$ to be $R_{\rm [K]}$ and
  \emph{column kernel support} of $B'$ at point $\eta$ to be $C_{\rm
    [K]}$. \emph{Nontrivial rows} of $B'$ are those that contain
  entries equal to $0$.
\end{definition}

We assume that a function \textsc{Kernel-Support} implements the
computation of the row and
column kernel support of the Jacobian matrix $J$ at point $\eta$ and
returns $B''$ retricted to columns not in $C_{\rm [K]}$ and rows not in
$R_{\rm [K]}$, with the notations of lem.~\ref{lem:basic_idea_reg}.

\begin{remark}\label{rem:complex-KS}
It is easily seen that the complexity of this algorithm is $O(s^{2}n)$
elementary operations in $\K$, using Gaussian elimination. In fact, we
may only consider nontrivial rows of $B'$ (containing elements
different from $-\infty$) and the corresponding rows of $J'$. If
  their number is $q$, the complexity is $O(q^{2}n+sn)$.
\end{remark}

\begin{algorithm}\label{algo:Kernel-Support}
\caption{\textsc{Kernel-Support}}\label{algo:KS}

\textbf{Input:} A matrix $B'$ of $0$ and $-\infty$
elements and a matrix $J$ of real elements with equal numbers of
rows and columns.

\textbf{Output:} The submatrices $B''$ of $B'$ and $J'$ of $J$ where
rows in $R_{\rm [K]}$ and columns in $C_{\rm [K]}$ have been
suppressed.

\begin{algorithmic}[0]
\Function{Kernel-Support}{$A,J$}
\State Compute the kernel $K$ of the linear mapping $L$ defined by the
rows of $J'$.
\State Compute the row support $R_{\rm [K]}$ and column support
$C_{\rm [K]}$ of $K$.
\State\textbf{return}
$(B'':=B'\setminus R_{\rm [K]}\setminus C_{\rm [K]},J':=J\setminus
R_{\rm [K]}\setminus C_{\rm [K]})$.
\EndFunction
\end{algorithmic}
\end{algorithm}

\PATCH{
  \begin{example}\label{ex:not_necessary_reg2}
    We start with ex.~\ref{ex:not_necessary_reg} at a
point where $x_{5}=0$ and $x_{6}\neq0$. We have
$$
      B'=\left(
      \begin{array}{c|c}
        -\infty&\Ast{0}\\
        \hline
        \Ast{0}&-\infty
      \end{array}
      \right)
      \quad\hbox{and}\quad
            J=\left(
      \begin{array}{c|c}
        -\infty&1\\
        \hline
        0&-\infty
      \end{array}
      \right),
$$ so that $R_{\rm [K]}$ is reduced to the last row of $B'$ and $C_{\rm [K]}$
      to its first column. The function \textsc{Kernel-Support} returns
      $(B'',J')=((0),(1))$.

      When $x_{5}=x_{6}=0$, $R_{\rm [K]}$
      contains the two rows and $B''=\emptyset$.
\end{example}

We continue with ex.~\ref{ex:p-q-3}.

\begin{example}\label{ex:p-q-3_2} We have:
$$
B'=
  \left(
  \begin{array}{c|ccc}
    0&-\infty&-\infty&-\infty\\
    \hline
    0&0&-\infty&-\infty\\
    -\infty&0&0&0
  \end{array}
  \right)
  \quad\hbox{and}\quad
J=  \left(
  \begin{array}{c|ccc}
    0&0&0&0\\
    \hline
    1&0&0&0\\
    0&1&0&1
  \end{array}
  \right),
  $$
so that $R_{\rm [K]}$ is reduced to the first row of $B'$ and $C_{\rm [K]}$
      to its first column. The function \textsc{Kernel-Support} returns
$$
B''=
  \left(
  \begin{array}{ccc}
    0&-\infty&-\infty\\
    0&0&0
  \end{array}
  \right)
  \quad\hbox{and}\quad
J'=  \left(
  \begin{array}{ccc}
    0&0&0\\
    1&0&1
  \end{array}
  \right).
  $$
\end{example}
}

We may now iterate \textsc{HK} and \textsc{Kernel-Support} until the returned
matrix $B''$ is equal to $B'$. We call the process \textsc{Seq}. We do
not need any more to assume that \textsc{HK} returns a set $Y_{1}$ of
columns, but only the sets $R_{0}$ and $C_{0}$.

\begin{algorithm}
  \caption{\textsc{Seq}}\label{algo:Seq}

\textbf{Input:} A matrix $B'$ of $0$ and $-\infty$
elements and a matrix $J$ of real elements with the same number of
rows and columns.

\textbf{Output:} Submatrices $(B'',J')$ of $B'$ and $J$ such that
\textsc{HK}$(B'')=(R_{0},C_{0})$ with $C_{0}=\emptyset$ and
\textsc{Kernel-Support}$(B'',J')=(B'',J')$.

\begin{algorithmic}[0]
\Function{Seq}{$B,J'$}
\State$(R_{0},C_{0}):=$\textsc{HK}$(B')$;
\State$(B'',J'):=$\textsc{Kernel-Support}$((B'\setminus C_{0})\cap R_{0},(J\setminus C_{0})\cap R_{0})$
\State\textbf{if} $B''=B'$ \textbf{then} \textbf{return} $(B'',
J')$ \textbf{else} \textbf{return} \textsc{Seq}$(B'',J')$ \textbf{fi}
\EndFunction
\end{algorithmic}
\end{algorithm}

To alleviate the presentation, we avoid going too deeply in
computational details. The following remark should be enough for our
purpose.

\begin{remark}\label{rem:compl_q}
As the number of rows of $B'$ decreases at each recursive call, the
number of iterations is bounded by $s$, so that the complexity is
$O(s^{3}n)$. We can obviously neglect the trivial rows of $B'$ and the
corresponding rows of $J'$ that are rows of $0$ elements. If $d$ is
the number of nontrivial rows in $B'$ and $q$ the number of
iterations, then the computations only require $O(qd^{2}n')$
operations, where $n'$ is the number of columns of $B$.
\end{remark}

\begin{example}\label{ex:not_necessary_reg3} The matrix $J'$ computed at
  ex.~\ref{ex:not_necessary_reg2} has full rank, so that the next call
  to \textsc{HK} and \textsc{Kernel-Support} will return again
  $(B'',J')$ and the iterations in \textsc{Seq} will stop.
\end{example}

\begin{example} Applying \textsc{Kernel-Support} to the matrices
  $(B',J)$ in ex.~\ref{ex:p-q-3_2}, we get the sequence of results
$(B'',J')$, $B''':=(0\>0)$, $J'':=(0\>1)$, and then $(B''',J'')$
  again, as $J''$ has full rank, so that the process stops.
\end{example}

In the following lemma, ii) and iii) are analogs of
lem.~\ref{lem:basic_idea}~iii) and iv).

\begin{lemma}\label{lem:Seq} i) a) With the hypotheses and  notations of
  lem.~\ref{lem:basic_idea_reg}, the process \textsc{Seq} returns
  $(B'',J')$, where 
  $B''$ is a submatrix of $B$ containing at least
  one entry equal to $0$ and such that \textsc{HK}$(B'')=(R_{0},C_{0})$
  with $C_{0}=\emptyset$ and \textsc{Kernel-Support}$(B'',J')=(B'',J')$. b)~It
  contains the intersection of the rows of $\bar R$ and the columns of
  $\sigma(\bar R)$. c)~The only elements different from $-\infty$ in the
  columns of $B''$ are contained in the rows of $B''$.

  ii) Let $R_{1}$ denote the rows of $B''$ and $C_{1}$ its columns.
  The matrix $A_{2}$ formed of rows $R_{2}$ not in $R_{1}$ and columns
  $C_{2}$ not in $C_{1}$ is such that there exists a set of columns
  $Y_{2}$ that satisfies $\hat\cO_{A_{2}}=\cO_{Y_{2},A_{2}}=0$ and
  $\nabla_{Y_{2},\Sigma_{2}}(\eta)\neq0$, where $\Sigma_{2}$ is the
  subset of $\Sigma$ $\Sigma_{2}:=\{P_{i}\mid i\in R_{2}\}$.

  iii) The matrix $B''$
  is such that there exists
  $Y_{1}$ with
  $\cO_{Y_{1},B''}=0$ and $\nabla_{Y_{1},B''}(\eta)\neq0$, which implies
  $\cO_{Y_{2}\cup Y_{1},\Sigma}=0$ and $\nabla_{Y_{2}\cup Y_{1},\Sigma}(\eta)\neq0$.
\end{lemma}
\begin{proof} i) a) It is 
  straightforward as when the process \textsc{Seq} stops, it must
  return a $B''$ such that \textsc{HK}$(B)=(R_{0},C_{0})$ with
  $C_{0}=\emptyset$ and \textsc{Kernel-Support}$(B'',J')=(B'',J')$.  b)~By
  lem.~\ref{lem:basic_idea_reg} i) and ii), the intersection of the
  rows of $\bar R$ and the columns of $\sigma(\bar R)$ belong to $B''$
  and $J''$ at each iteration of algo.~\ref{algo:Seq}, so that they
  are contained in the matrices of the output. c)~In the same way, at
  each iteration, the only elements different from $-\infty$ in the
  columns of $B''$ are contained in the set of rows $R_{1}$, so that this
  property also stands for the output.
  
  ii) and iii) We know that there exists $Y$ such that
  $\cO_{Y,\Sigma}=0$ and $\nabla_{Y,\Sigma}(\eta)\neq0$. So, one may
  find $Y_{2}\subset Y$, such that
  $\nabla_{Y_{2},\Sigma_{2}}(\eta)\neq0$ and
  $\cO_{Y_{2},\Sigma_{2}}=0$. As $J'$ has full rank at
  $\eta$, there exists $Y_{1}$ such that
  $\nabla_{Y_{1},\Sigma_{1}}(\eta)\neq0$ and
  $\cO_{Y_{1},\Sigma_{1}}=0$, where $\Sigma_{1}$
  corresponds to rows in $R_{1}$. Then,
  i)~c)~implies that
$$
\cO_{Y_{1}\cup Y_{2},\Sigma}= \cO_{Y_{1},\Sigma_{1}}+\cO_{Y_{2},\Sigma_{2}}=0,
$$
and
$$
\nabla_{Y\cup Y_{2},\Sigma}(\eta)=\nabla_{Y_{1},\Sigma_{1}}(\eta)
                      \nabla_{Y_{2},\Sigma_{2}}(\eta)\neq0.
$$
\end{proof}

It is then easy to design an algorithm \textsc{\=O-reg} to test
\=o-regularity at point $\eta$. The routine \textsc{Seq2} is assumed
to be a variant that returns the set of columns $Y_{1}$ and the set of
equations $\Sigma_{2}$ of lem.~\ref{lem:Seq}~iii) or ``failed'' if $B''$
is empty.

\begin{algorithm}
  \caption{\textsc{\=O-reg}}\label{algo:o-reg}

  \textbf{Input:} A differential system $\Sigma$ of $s$ equations in $n$
variables $x_{1}$, \dots, $x_{n}$, defining a diffiety in the
neighborhood of $\eta\in\J(\R,\R^{n})$.

\textbf{Output:} ``failed'' or a set $Y$ of rows such that
$\cO_{Y,\Sigma}=0$ and $\nabla_{Y,\Sigma}(\eta)\neq0$.

\begin{algorithmic}[0]
\Function{\=O-reg}{$\Sigma,\eta$}
\State $A:=A_{\Sigma}$;
\State \textbf{if} $s>n$ or $A$ contains a full row of $-\infty$ \textbf{then}
\textbf{return} ``failed'';
\State Suppress from $A$ all columns of $-\infty$;
\State Build $B$ the submatrix of $A$ of columns $C$ containing only
$-\infty$ or $0$ elements and
$J:=(\partial P_{i}/\partial x_{j}\mid 1\le i\le s,\> j\in C)$
\State \textbf{if} \textsc{Seq2}$(B,J)=$ ``failed'' \textbf{then}
\textbf{return} ``failed''
\textbf{else} $(Y_{1},\Sigma_{2}):=$\textsc{Seq2}$(B,J)$ \textbf{fi}
\State \textbf{if} \textsc{\=O-reg}$(\Sigma_{2},\eta)=$ ``failed''
\textbf{return} ``failed'' \textbf{else} $Y_{2}:=$\textsc{\=O-reg}$(\Sigma_{2},\eta)$ \textbf{fi}
\State\textbf{return} $Y_{1}\cup Y_{2}$
\EndFunction
\end{algorithmic}
\end{algorithm}

The following theorem provides an evaluation of the complexity.

\begin{theorem}\label{th:complexity_reg_o-syst} i) The
  algorithm \textsc{\=o-reg}~\ref{algo:o-reg} tests if $\Sigma$ is
  \=o-regular at point $\eta$ and then returns $Y$ such that
  $\cO_{Y,\Sigma}=0$ and $\nabla_{Y,\Sigma}(\eta)\neq0$.

  ii) It works in $O(p(qd^{2}+s)n)$ elementary operations in $\K$, where

  --- $p$ is the number of recursive calls of \textsc{\=o-reg};

  --- $q$ and $d$ are respectively the
  maximal number of iterations of \textsc{Seq} and the maximal number
  of nontrivial rows in $B$ at each recursive call of  \textsc{\=o-reg}.

  iii) If $q$ is equal to $1$ and $C_{0}=\emptyset$ at each recursive
  call of \textsc{\=o-reg}, then using balanced trees
  (rem.~\ref{rem:AVL}) as in algo.~\ref{algo:o-system}, we can achieve
  a complexity of $O((d^{2}+\ln(s)s)n$.
\end{theorem}
\begin{proof} i) Proceeding as in the proof of
  th.~\ref{th:complexity_o-syst}, it is a straightforward consequence of
  lem.~\ref{lem:Seq}.

  ii) The construction of $B$ at each recursive call requires $O(sn)$
  operations and \textsc{Seq} $O(qd^{2}n)$ operations, using
  rem.~\ref{rem:compl_q}.

  iii) Using balanced trees, the first construction of $B$ requires
  $O(\ln(s) sn)$ operations and further constructions only $O(\ln(s)
  n)$. Let $n_{h}$ be the number of columns in $B$ at recursive call
  $h$ of \textsc{\=o-reg}. By rem.~\ref{rem:compl_q}, \textsc{Seq}
  requires then $O(d^{2}n_{h})$ operations, so that the total cost is
  $O(\ln(s)sn)+\sum_{h=1}^{p}O(d^{2}n_{h})=O((\ln(s)s+d^{2})n)$
  operations. Indeed, if $C_{0}=\emptyset$, all columns of $B$ only
  contain $-\infty$ elements in $A_{2}$ and can be removed.
\end{proof}

\PATCH{
\begin{example}
Going back to ex.~\ref{ex:not_necessary_reg} in the case $x_{5}=0$ and
$x_{6}\neq 0$, the computations in
ex.~\ref{ex:not_necessary_reg3} show that \textsc{Seq2} will return
$Y_{1}=\{x_{6}\}$ and $\Sigma_{2}$ that corresponds to
$\{P_{2},P_{3},P_{4}\}$. In the following recursive call, all Jacobian
matrices have full rank, so that the first iteration in \textsc{Seq}
returns the good result. The successive values returned by
\textsc{\=O-reg} are
\textsc{\=O-reg}$(\Sigma_{2})=\{x_{1},x_{3},x_{4}\}$,
\textsc{\=O-reg}$(\Sigma_{3}:=\{P_{2},P_{4}\})=\{x_{1},x_{4}\}$ and
\textsc{\=O-reg}$(\Sigma_{4}:=\{P_{4}\})=\{x_{4}\}$, so that
\textsc{\=O-reg}$(\Sigma_{2})=\{x_{1},x_{3},x_{4},x_{6}\}$.

When $x_{5}=x_{6}=0$, $B''=\emptyset$ and \textsc{Seq2} and
\textsc{\=O-reg} returns ``failed''.
\end{example}
}

The following example shows that $p$ and $q$ can both be equal to $s$.

\begin{example}\label{ex:p-q}
We consider the system $\Sigma_{s}$ defined by
$P_{1}:=x_{1}^{2}+x_{2s}+x_{2s+1}'$ and
$P_{i}:=x_{i-1}+x_{i}^{2}+x_{2s-i+1}+x_{2s-i+2}'$, for $1<i\le s$
in $2s+1$ 
variables $x_{j}$ at a point $\eta=0$ where $x_{i}=0$, for $1\le i\le
s$. \PATCH{Ex.~\ref{ex:p-q-3} corresponds to $s=3$. Again in the
  general case, the first row and the first column are removed at each
  iteration of \textsc{Kernel-Support}, until the last, for which the
  final value of the Jacobian matrix has full rank and the sequence
  stops.

Starting with $\Sigma_{s}$, the system considered at the $h^{th}$
recursive call is $\Sigma_{s-h}(x_{1}, \ldots, x_{s-h}, x_{s+h+1},
\ldots, x_{2s+1})$.} So, we have then $p=s$ and at iteration $h$ of
\textsc{\=o-reg}, we have $q=s-h$ iterations in \textsc{Seq} with
columns of $B'$ reducing from $\{1, \ldots, h, 2s-h+1\}$ to $\{h,
2s-h+1\}$.
\end{example}

\subsection{Examples}\label{subsec:examples_chained}

We will consider here some examples of classes of ``chained'' or
``triangular'' flat systems in the literature. We stress on the fact
that in the papers quoted here, the main issue is to test the
existence of a change of variables that may reduce a given system to
such a form, whereas our problem here is to test if a system is already in
such a form, up to a permutation of indices.

\subsubsection{Goursat normal form}

It is known that all driftless systems with two controls of the
general form
\begin{equation}\label{eq:driftless_two_controls}
  x_{i}'=f_{i}(x)u+g_{i}(x)v,\>\hbox{for}\> 1\le i\le n
\end{equation}
can be reduced to the Goursat normal form
\begin{equation}\label{eq:Goursat}
  \begin{array}{lll}
    z_{0}'&=&v_{0};\\
    z_{i}'&=&z_{i+1}v_{0}\>\hbox{for}\> 1\le i\le n-2;\\
    z_{n-1}'&=&v_{1},
  \end{array}
\end{equation}
iff it is flat and we may use the following flatness criterion.
\begin{theorem}\label{th:Cartan}
  A driftless systems with two
  controls~\eqref{eq:driftless_two_controls} is flat iff the vector
  spaces $E_{i}$, $0\le i\le n-2$ defined by $E_{0}:=\{f,g\}$ and
  $E_{i+1}=E_{i}+\{[E_{i},E_{i}]\}$ satisfy $\dim E_{i}=i+2$ for all
    $0\le i\le n-2$. 
\end{theorem}
This result goes back to the work of
Cartan~\cite{Cartan1914,Cartan1915} and has been adapted to control by
Martin and Rouchon~\cite{Martin1994}. See also Li \textit{et
  al.}~\cite{Li2012}.

The system \eqref{eq:Goursat} is obviously a \=o-system, with a single
set of possible flat outputs; $\{z_{0},z_{1}\}$, which is regular iff
$v_{0}\neq0$. 

\subsubsection{Complexity issues}\label{subsubsec:complexity_issues}

Without going into useless details, for which we refer to the
references quoted above, we need to give some idea of the complexity
of computations involved to work out a flat parametrization after
having proved the existence of a suitable change of variables.  We
assume here for simplicity that the fields $f$, $g$, \dots\ are
defined by rational functions of the state variables and that vector
spaces are $\R(x)$-vector spaces. For simplicity, we identify the
fields $f$, $g$, \dots\ with the associated derivations.

The next theorem is the basis of a step
by step reduction of a two inputs driftless system in Goursat normal
form.

\begin{proposition}\label{th:devissage_Goursat}
  Assume that a two inputs driftless system
  \eqref{eq:driftless_two_controls}, with $n>3$ states, admits a
  change of variable $y_{i}=Y_{i}(x)$ such that the system becomes
\begin{equation}\label{eq:devissage}
\begin{array}{lll}
  y_{i}'&=&(\bar f_{i}(y_{1}, \ldots, y_{n-1})
  +\bar g(y_{1}, \ldots, y_{n-1})y_{n})\bar u;\>\hbox{for}\> 1\le i\le n-1;\\
  y_{n}'&=&\bar{w},
\end{array}
\end{equation}
with $[\bar f,\bar g]\notin \langle\bar f,\bar g\rangle$. Then $\dim
E_{1}=3$ and $\dim E_{2}=4$.
\end{proposition}
\begin{proof}
  Using the new coordinates, we have have $E_{3}=\langle \partial_{y_{n}},
  \bar f, \bar g\rangle$ and $E_{4}=\langle \partial_{y_{n}},
  \bar f, \bar g, [\bar f, \bar g]\rangle$, that must have respective
  dimensions $3$ and $4$, according to our independence hypothesis.
\end{proof}

Easy computations imply the following corollary.

\begin{corollary} Under the hypotheses of the proposition, there
  exists a couple of functions of the state variables $x$
  $(a,b) \neq (0,0)$ such that $a[f,[f,g]]+b[g,[f,g]]=0$ modulo
  $E_{3}$. Then, $af+bg=c\partial_{y_{n}}$ and the $y_{i}$, $1\le i\le
  n-1$ are functionally independent first integrals common to the
  field $af+bg$.
\end{corollary}

Using this lemma, we are reduced to a new two inputs system
\begin{equation}\label{eq:devissage2}
\begin{array}{lll}
  y_{i}'&=&(\bar f_{i}(y_{1}, \ldots, y_{n-1})
  +\bar g(y_{1}, \ldots, y_{n-1})y_{n})\bar u;\>\hbox{for}\> 1\le i\le n-1;\\
  &=&\bar f_{i}(y_{1}, \ldots, y_{n-1})\bar u
     +\bar g(y_{1}, \ldots, y_{n-1})\bar v,
\end{array}
\end{equation}
  with $\bar v=y_{n}\bar u$.

  Successive applications of this process reduces the state dimension
  and produces a Goursat normal form and a flat parametrization. One
  may notice that for $n=3$, all combinations $a(x)f+b(x)g$ work.

  The main issue then is to look for first integrals. There can exist
  no rational solutions and there is no general method to test if a
  rational solution exists, even when looking for first integrals of a
  field in the affine plane. Already looking for the existence of
  such an integral up to to a given degree is computationally
  difficult. See \textit{e.g.} Ch{\`e}ze and Combot~\cite{Cheze2020} and
  the references therein for more details.

  One may also look for closed forms solutions as did Rouchon for the
  car with one trailer in the general case~\cite{Rouchon1993b}, but
  again there might not exist any. This does not mean that the task is
  hopeless, but justifies some special interest to situations where
  the computations are much easier, although not completely trivial,
  mostly when the size of initial equations is already
  appreciable. Such situations are obviously nongeneric, but may
  often be encountered in practice with the help of some
  simplifications. This is not uncommon with flat systems, that are
  themselves nongeneric but quite ubiquitous in engineering practice.

\begin{example}
An affine generalization with two inputs has been considered by
Silveira~\cite{Silveira-PHD} and Silveira \textit{et
  al.}~\cite{Silveira2015}:
\begin{equation}\label{eq:Silveira}
  \begin{array}{lll}
    z_{0}'&=&v_{0};\\
    z_{i}'&=&f_{i}(z_{0},z_{1}, \ldots, z_{i+1})+
    z_{i+1}v_{0}\>\hbox{for}\> 1\le i\le n-2;\\
    z_{n-1}'&=&v_{1}.
  \end{array}
\end{equation}
They provide necessary and sufficient conditions to reduce a system of the form
$x'=f(x)+g_{1}(x)u_{1}+g_{2}(x)u_{2}$ to the form~\eqref{eq:Silveira}.

Such a system is oudephippical. Using algo.~\ref{algo:o-system}, one
can conclude that all matrices $B$ are such that
$C_{0}=\emptyset$. The sets $\Xi_{h}$ of th.~\ref{th:suff_reg_cond}
are $\Xi_{1}=\{v_{0},v_{1}\}$ and $\Xi_{h}=\{z_{h}\}$, for
$1<h<n-1$. For best efficiency, a sparse version of the algorithm
should be designed for such sparse systems.  The only possible flat
outputs set in the setting of th.~\ref{th:suff_reg_cond} is
$\{z_{0},z_{1}\}$ and the regularity condition is $v_{0}+\partial
f_{i}/\partial z_{i+1}\neq0$, for all $1\le i\le n-2$.
\end{example}
\vfill\eject
  
\subsubsection{Multi input chained forms}

Some notions of chained forms may be found in the literature for
systems with many inputs.

\begin{example}\label{ex:Li2016}
A multi-input generalization, the ``$m$-chained form'', has been
proposed by Li \textit{et al.}~\cite{Li2016}:
\begin{equation}\label{eq:m-chained}
  \begin{array}{lll}
    z_{0}'&=&v_{0};\\
    z_{i,\ell}'&=&f_{i,\ell}(z_{0},\bar z_{\ell})
    +z_{i,\ell+1,}v_{0}\>\hbox{for}\> 1\le i\le
    m\>\hbox{and}\> 1\le\ell<k;\\
    z_{i,k}'&=&v_{i}, \>\hbox{for}\> 1\le i\le m
  \end{array}
\end{equation}
  where $\bar z_{\ell}:=(z_{1,1}, \ldots, z_{1,\ell}, \ldots, z_{m,1},
  \ldots, z_{m,\ell})$.
  
This system is also oudephippical. All matrices $B$ of
algo.~\ref{algo:o-system} are again such that $C_{0}=\emptyset$. The
sets $\Xi_{h}$ of th.~\ref{th:suff_reg_cond} are
$\Xi_{1}=\{v_{0}, \ldots, v_{m}\}$ and $\Xi_{h}=\{z_{1,h},
\ldots, z_{m,h}\}$,
for
$1< h\le k$. The only possible flat
outputs set in the setting of th.~\ref{th:suff_reg_cond} is
$\{z_{0},z_{1,1}, \ldots, z_{m,1}\}$
and the regularity condition is $v_{0}\neq0$.
\end{example}

The authors consider the case of a rolling coin on a moving table.
$$\left(\begin{array}{c} x'\\ y'\\ \theta'\\ \phi'\end{array}\right)=
  \left(\begin{array}{c} \cos\theta(\alpha\cos\theta+\beta\sin\theta)\\
    \sin\theta(\alpha\cos\theta+\beta\sin\theta)\\
    0\\0\end{array}\right)
+\left(\begin{array}{c} 0\\0\\1\\0\end{array}\right)u_{1}
+\left(\begin{array}{c} R\cos\theta\\R\sin\theta\\0\\1\end{array}\right)u_{2}
$$ where $\alpha$ and $\beta$ are known functions of the time, that
  describe the motion of the table and $R$ is a constant.  They show
  that it can be reduced to the form \eqref{ex:Li2016} only in the
  case of a constant speed rotation of the table.  Nevertheless, one
  may choose $\theta$ as a flat output. The system becomes linear in
  the remaining variables and is flat when $\theta'\neq0$, with flat
  output $R\phi-\cos(\theta) x-\sin(\theta)y$.

This is very close to the example
$x'=\cos(\theta)\phi'$, $y'=\sin(\theta)\phi'$, which may be traced
back to Monge~\cite{Monge1787} under the form $\dd\phi^{2}=\dd
x^{2}+\dd y^{2}$. It has been rediscovered independently by
Petitot~\cite{Petitot1992}.

We see that it is not obvious to design meaningful examples
for which the changes of coordinates described by PDE remain
tractable. Working with such functions, when they cannot be expressed
by closed form formulas remains a challenge.

All these examples enter in more general notions of block diagonal or
chained systems.

\subsubsection{Block diagonal and chained systems}

``Almost chained systems'', as described
in~\cite{Ollivier2022},
\begin{equation}
  (Z'_{h}, X_{h}')=G_{h}(Z_{1}, \ldots Z_{h+1}, X_{1}, \ldots , X_{h+1})
+H_{h}(X_{h+2}, \ldots , X_{h+\ell_{h}} ), 1 \le h \le r, 
\end{equation}
are chained when the extra functions $H_{h}$ are $0$. We propose here
some more precise definition.

\begin{definition}\label{def:block_triang}
An \emph{order $1$ block triangular systems} is a system $\Sigma$ in
the variables $\Xi$, for which there exist partitions
$\Sigma=\bigcup_{h=1}^{p}\Sigma_{h}$ and
$\Xi=\bigcup_{h=0}^{p}\Xi_{h}$ such that all equations in $\Sigma_{i}$
depend only in variables in $\bigcup_{k=0}^{i}\Xi_{k}$ and are of
order $1$ in variables of $\Xi_{i-1}$ and $0$ at most in the other
variables, with

i) $\sharp\Xi_{h-1}=\sharp\Sigma_{h}$;

ii) $\cO_{\Xi_{h-1},\Sigma_{h}}=\sharp\Xi_{h-1}$;

iii) $\cO_{\Xi_{h},\Sigma_{h}}=0$.

It is said to be \emph{dense} if moreover we have

iv) $a_{i,j}=0$ in $A_{\Sigma_{h}}$ for all $1\le h\le p$ and all
$(i,j)\in\Sigma_{h}\times \Xi_{h}$.

An order $1$ block triangular system is said to be \emph{chained} at
level $h>0$ if all equations in $\Sigma_{h}$ depend only in variables
in $\Xi_{h}\cap\Xi_{h-1}$. It is said to be \emph{stricly chained} if
is chained at level $h$ and the equation in $\Sigma_{h}$ depend only
on derivatives of order $1$ of the variables in $\Xi_{h-1}$ and not of
those variables themselves.
\end{definition}

\begin{remark}
Condition iii) means that there is an
  injection $\sigma:\Sigma_{h}\mapsto\Xi_{h}$ with
  $\sum_{P\in\Sigma_{h}}\ord_{\sigma(P)}P=0$. Then, considered as a
  system in the variables of $\bigcup_{h=1}^{p}\sigma(\Sigma_{h})$, 
  the system $\Sigma$ is block triangular according to the definition
  of~\cite[4.3]{Ollivier2022b}.
\end{remark}

We have the following proposition, of which the easy proof is left to
the reader.

\begin{proposition} Order $1$ block triangular systems are \=o-systems
  such that $C_{0}=\emptyset$ at each step of
  algo.~\ref{algo:o-system}.

  With the notations of this algorithm and of the previous definition,
  the set $\Xi_{p-h+1}$ (resp.~$\Sigma_{p-h+1}$), for $1\le h\le p$
  corresponds to the columns (resp.~the rows) of $B$ at step $h$ of
  the algorithm.
\end{proposition}

The next subsection will provide some sufficient conditions of
regularity and singularity in the block diagonal case.

\subsection{Some special results for block
  triangular systems}\label{subsubsec:suf_sing_cond}

\subsubsection{A sufficient condition for \=o-regularity}\label{subseubsec:o-reg_block}

We have seen with th.~\ref{th:suff_reg_cond} a sufficient condition
for regularity. Example~\ref{ex:not_necessary_reg} shows that it is
not a necessary solution for \=o-regularity for all \=o-systems. The
following corollary shows that the condition of
th.~\ref{th:suff_reg_cond} is indeed a necessary and sufficient
condition of \=o-regularity in the case of generic block diagonal systems.

\begin{corollary}\label{cor:o-reg_block}
  With the hypotheses of th.~\ref{th:suff_reg_cond}, assume that

  a) $\Xi$
  is a dense order $1$ block triangular system or that

  b) for all $1\le h\le p$ and all $\Sigma'\subset\Sigma_{h}$ with
  $\sharp\Sigma'=\sharp\Sigma_{h}-1$, the Jacobian matrix $(\partial
  P/\partial x \mid (P,x)\in \Sigma'\times\Xi_{h})$ has rank equal to
  $\sharp\Sigma_{h}-1$.

  Then, it is
  \=o-regular at point $\eta$ iff for all $1\le h\le p$, we have
  eq.~\ref{eq:suf_reg_cond}, \textit{i.e.}
  \begin{equation}\label{eq:suf_reg_cond_h}
      \rank\,\left(\frac{\partial P}{\partial x}(\eta)
      \mid P \in \Sigma_{h},\ x\in\Xi_{h}\right)=\sharp\Sigma_{h}.
  \end{equation}
\end{corollary}
\begin{proof} The sufficient part is th.~\ref{th:suff_reg_cond}.

  The necessary part is a consequence of the correctness of
  algo.~\ref{algo:o-reg} th.~\ref{th:complexity_reg_o-syst}~i). If a),
  then if the rank of \eqref{eq:suf_reg_cond_h} is not
  $\sharp\Sigma_{h}$ at level $h$, then at stage $p-h$ of the
  algorithm, all columns of $B$ are suppressed, so that the system is
  not regular. The same applies using hypothesis b), as then relations
  between the rows of the Jacobian matrix must imply all
  $\sharp\Sigma_{h}$ rows, so that again we need to suppress all
  columns of $B$. 
\end{proof}

\subsubsection{A sufficient condition for singularity}
We will need some technical lemma about linear systems.

\begin{lemma}\label{lem:sub_sys} If a block diagonal linear system
  $\Sigma$, chained at level 
  $i_{0}$ is flat, then, using 
  the notations of def.~\ref{def:block_triang}, for all $1\le
  i_{0}<i_{1}\le p$, the block diagonal system
  $\bar\Sigma:=\bigcup_{i=i_{0}}^{i_{1}}\Sigma_{i}$ in the variables
  $\bar\Xi:=\bigcup_{i=i_{0}-1}^{i_{1}}\Xi_{i}$ is flat.
\end{lemma}
\begin{proof} By the structure theorem, the linear system $\Sigma$ is
  flat iff the module
  $$\left [\sum_{x\in\Xi}\R(t)[\dd/\dd t]x \right ]/\R(t)[\dd/\dd t]\Sigma$$
  contains no torsion element, which implies that the module
  $\left [ \sum_{x\in\bar\Xi}\R(t)[\dd/\dd t]x \right ]/\R(t)[\dd/\dd t]\bar\Sigma$
  contains no torsion element.
\end{proof}

We can now state the following sufficient condition of singularity. 

\begin{theorem}\label{th:suff_sing}
  Let $\Sigma=\bigcup_{i=1}^{p}\Sigma_{h}$ be a block diagonal system
  in the variables $\bigcup_{h=0}^{p}\Xi_{h}$,
  chained at level $h_{0}$ and strictly chained at level $h_{0}-1$,
  and such that all variables in $\Xi_{h_{0}-1}\cup\Xi_{h_{0}}$
  appearing nonlinearly in $\Sigma_{h_{0}-1}\cup\Sigma_{h_{0}}$ are
  constants along a given trajectory. Let $\eta$ denote a point of
  this trajectory.

  Assume moreover that the Jacobian matrix
\begin{equation}\label{eq:Rank1}
\left(\frac{\partial P}{\partial x}\mid
(P,x)\in\Sigma_{h_{0}}\times\Xi_{h_{0}}\right) 
\end{equation}
has rank $m_{0}<n_{0}:=\sharp\Xi_{h_{0}-1}$ and that the Jacobian
determinants
\begin{equation}\label{eq:Rank2}
\left|\frac{\partial P}{\partial x'}\mid
(P,x)\in\Sigma_{h_{0}}\times\Xi_{h_{0}-1}\right| 
\end{equation}
do not vanish at $\eta$ for all $1\le h\le p$.

With these hypotheses, $\Sigma$ is not flat at $\eta$.
\end{theorem}
\begin{proof} The nonvanishing of the determinants \eqref{eq:Rank2}
  implies the existence of explicit differential equations
  \begin{equation}\label{eq:th_sing}
    x'=f_{x}(\Xi_{0}, \ldots, \Xi_{h}),\> x\in\Xi_{h}\>\hbox{for}\>
    1\le h\le p.
  \end{equation}
  We will consider the linearized system (as defined in
  def.~\ref{def:power-series})
  \begin{equation}\label{eq:th_sing2}
    \dd_{\eta} x'=\dd_{\eta} f_{x}(\dd \Xi_{0}, \ldots, \dd \Xi_{h-1}),\>
    x\in\Xi_{h}\>\hbox{for}\> 0\le h<p.
  \end{equation}
As the system is chained at level $h_{0}$, for $x\in\Xi_{h_{0}-1}$,
$f_{x}$ and $\dd f_{x}$ only depend respectively on
$\Xi_{h_{0}}\cup\Xi_{h_{0}-1}$ and $\dd(\Xi_{h}\cup\Xi_{h_{0}-1})$. In
the same way, as the system is strictly chained at level $h_{0}-1$,
for $x\in\Xi_{h_{0}-2}$, $f_{x}$ and $\dd f_{x}$ only depend
respectively on $\Xi_{h-1}$ and $\dd\Xi_{h-1}$. Moreover, as the
$x\in\Xi_{h_{0}-1}\cup\Xi_{h_{0}}$ appearing nonlinearly are
constants, the coefficients in the $f_{x}$, for
$x\in\Xi_{h_{0}-1}\cup\Xi_{h_{0}}$, are constants too.
\medskip

We will show that the linearized system, truncated at level $h_{0}$,
admits torsion elements, so that the system $\Xi$ is not flat by
lem.~\ref{lem:sub_sys}. Such torsion elements are first integrals of the
Lie algebra generated by $\partial_{x}$, for $x\in\Xi_{h}$ and 
$$
\uptau:=\partial_{t}+\sum_{x\in\bigcup_{h=0}^{h_{0}-1}\Xi_{h}}f_{x}\partial_{x},
$$
where $\partial_{x}$ denotes $\partial/\partial x$.
See \cite[4.3]{Kaminski-et-al-2020} for more details on such
constructions. 
\medskip

As the rank of the Jacobian matrix \eqref{eq:Rank1} is $n_{0}<m$,
$$
\Lie(\uptau;\partial_{x}\mid x\in\Xi_{h_{0}})
=\langle\partial_{x}\mid x\in\Xi_{h_{0}}\rangle
+\Lie(\uptau;[\uptau,\partial_{w}]\mid w\in W\subset\Xi_{h_{0}}),
$$
for some subset of variables $W$, such that $\sharp W=n_{0}$. We
denote $[\uptau,\partial w]$ by 
$\hat\uptau\partial_{w}$ and $\hat\uptau^{k+1}\partial
w:=[\uptau,\hat\uptau^{k}\partial w]$. One may find integers $e_{w}$,
$w\in W$
such that the $\pi_{h_{0}-1}\hat\uptau^{k}\partial_{x}$, for $x\in W$ and
$1\le k\le e_{w}$ are linearly independent, where $\pi_{h_{0}-1}$ denotes
the projection on the vector space $\langle\partial_{x}\mid
x\in\Xi_{h_{0}-1}\rangle$. Those integers may be chosen so that
$\sum_{w\in W}e_{w}$ is maximal.
\medskip

Assume now that
$$
\pi_{h_{0}-1}\hat\uptau^{e_{w}+1}\partial_{w}=
\sum_{x\in W}\sum_{k=0}^{e_{x}}c_{w,x,k}\pi_{h_{0}-1}\hat\uptau^{k}\partial_{x}
$$
and define
$$
\varpi(w):=\hat\uptau^{e_{w}+1}\partial_{w}
-\sum_{x\in W}\sum_{k=0}^{e_{x}}c_{w,x,k}\hat\uptau^{k}\partial_{x},
$$
for all $w\in W$.
Then, we have
$$
\begin{array}{ll}
\Lie(\uptau;\partial_{x}\mid x\in\Xi_{h_{0}})
=&\langle\partial_{x}\mid x\in\Xi_{h_{0}}\rangle
+\langle\hat\uptau^{k}\delta_{x}\mid w\in W,\> 0\le k\le e_{w}\rangle\\
&+\Lie(\uptau;\varpi(w)\mid w\in W).
\end{array}
$$ The brackets $[\uptau,\varpi(w)]$, for $w\in W$, do not involve
derivations $\partial_{x}$, for $x\in\Xi_{h_{0}-1}\cup\Xi_{h_{0}-2}$,
as the system is chained at level $h_{0}$ and strictly chained at
level $h_{0}-1$. All variables that appear nonlinearly are constant,
so that all coefficients are constants. This implies that the
projection $\pi_{h_{0}-2}\Lie(\uptau;\varpi(w)\mid w\in W)$, where
$\pi_{h_{0}-2}$ denotes the projection on the vector space
$\langle\partial_{x}\mid x\in\Xi_{h_{0}-2}\rangle$, is equal to
$\pi_{h_{0}-2}\langle\varpi(w)\mid w\in W\rangle$, so that the
intersection of $\Lie(\uptau,\varpi(w)\mid w\in W)$ with
$\langle\partial_{x}\mid x\in \Xi_{h_{0}-2} \rangle$ has dimension at
most $\sharp W=n_{0}<\sharp\Xi_{h_{0}-2}$. This implies that the
dimension of $\Lie(\uptau,\partial_{x}\mid x\in\Xi_{h_{0}})$ is not
the maximal dimension $\sum_{h=0}^{h_{0}}\sharp\Xi_{h}$ and that
nontrivial torsion elements exist, which concludes the proof.
\end{proof}
\vfill\eject

\section{The simplified aircraft. A block diagonal system}\label{sec:aircraft}

\subsection{The aircraft model}\label{subsec:model}

\subsubsection{Nomenclature}\label{sec:nomenclature}

In this section, we have collected all the notations to make easier
reading the sequel\PATCH{\footnote{Wing span $a$ and mean aerodynamic
    chord $b$ are respectively denoted by $b$ and $c$ in
    \cite{Grauer-Morelli}.}}.

{\small
  
\begin{multicols}{3}
{\parindent=0pt\parskip=1mm  

  \textbf{Roman}
  
  $a$: wing span
  
  $b$: mean aerodynamic chord

  $C_{x}$, $C_{y}$, $C_{z}$: aerodynamic force coefficients, wind frame

  $C_{D}$, $C_{Y}$, $C_{L}$: \PATCH{aerodynamic force coefficients
  in~\cite{Grauer-Morelli}}

  $C_{l}$, $C_{m}$, $C_{n}$: aerodynamic moment coefficients
  
    $F$: thrust

  $L$, $M$, $N$: aerodynamic moments

  $m$: mass

  $p$, $q$, $r$: roll, pitch and yaw rates

  $S$: wing area

  $V$: airspeed

  $X$, $Y$, $Z$: aerodynamic forces

  $y_{p}$: distance of the engines to the plane of symmetry

  \textbf{Greek}

  $\alpha$ angle of attack

  $\beta$ sideslip angle

  $\gamma$: flight path angle

  $\eta$: differential thrust ratio 

  $\vartheta$: pitch angle

  $\theta$: parameters

  $\mu$: bank angle

  $\phi$: roll angle

  $\chi$: aerodynamic azimuth or heading angle

  $\psi$: yaw angle

  $\delta_{l}$,   $\delta_{m}$,   $\delta_{n}$: aileron, elevator,
  rudder deflexion

  $\theta$: model parameters

}

\end{multicols}

}

The model presented here relies on Martin~\cite{Martin-PHD,Martin-96}. On may
also refer to Asselin~\cite{Asselin1997},
Gudmundsson~\cite{Gudmundsson} or McLean~\cite{McLean} for more details.
%================================

\subsubsection{Earth frame, wind frame and body frame}

\begin{figure}[h!]
\begin{center}
\hbox to
\hsize{\hss
  \includegraphicsh[height=5cm]{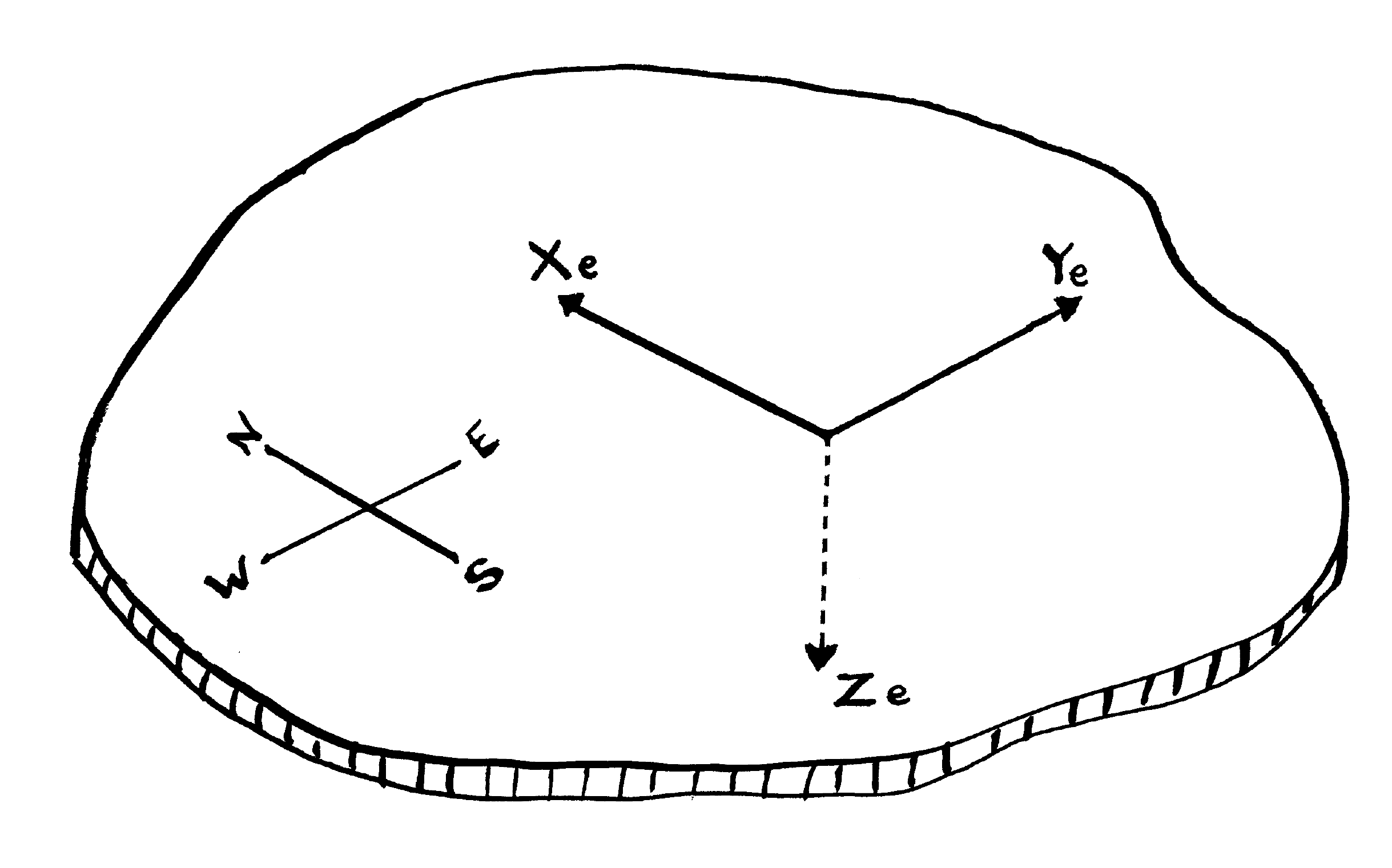}{a)}
  \hss
\includegraphicsh[height=5cm]{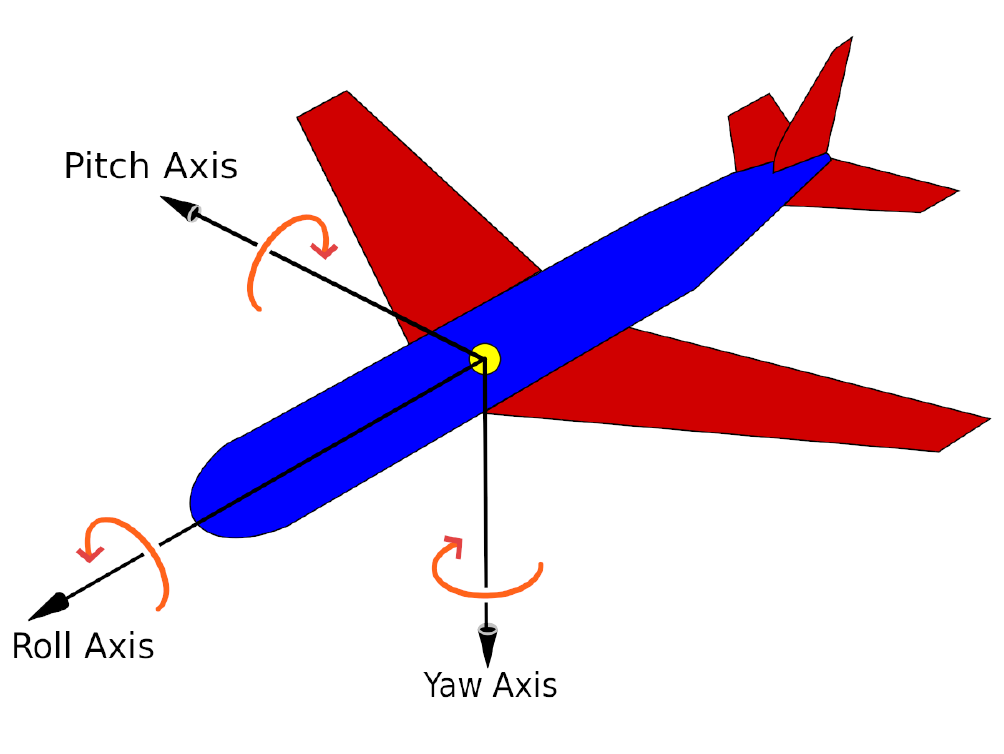}{b)}\hss}
\caption{a)~Earth frame and b)~body frame.} 
\label{fig::earth_frame}
\end{center}
\end{figure}

We use an earth frame with origin at ground level, with $z$-axis
pointing downward, as in the figure~\ref{fig::earth_frame}~(a). The
coordinates of the gravity center of the aircraft are given in this
referential.

%\begin{figure}[h!]
%\begin{center}
%\hbox to
%\hsize{\hss
%  \includegraphics[height=5cm]{./Earth_ref_ft.pdf}{a)}
%  \hss
%\includegraphics[height=5cm]{./body_frame.pdf}{b)}\hss}
%\caption{a)~Earth frame and b)~body frame.} 
%\label{fig::earth_frame}
%\end{center}
%\end{figure}

The body frame or aircraft referential is defined as in
figure~\ref{fig::earth_frame}~(b), where $x_{b}$ corresponds to the
roll axis, $y_{b}$ to the pitch axes and $z_{b}$ to the yaw axes,
oriented downward. The angular velocity vector $(p,q,r)$ is given in
this referential, or to be more precise, at each time, in the Galilean
referential that is tangent to this referential.

\begin{figure}[h!]
\begin{center}
  \hbox to\hsize{\hss\includegraphicsh[height=5cm]{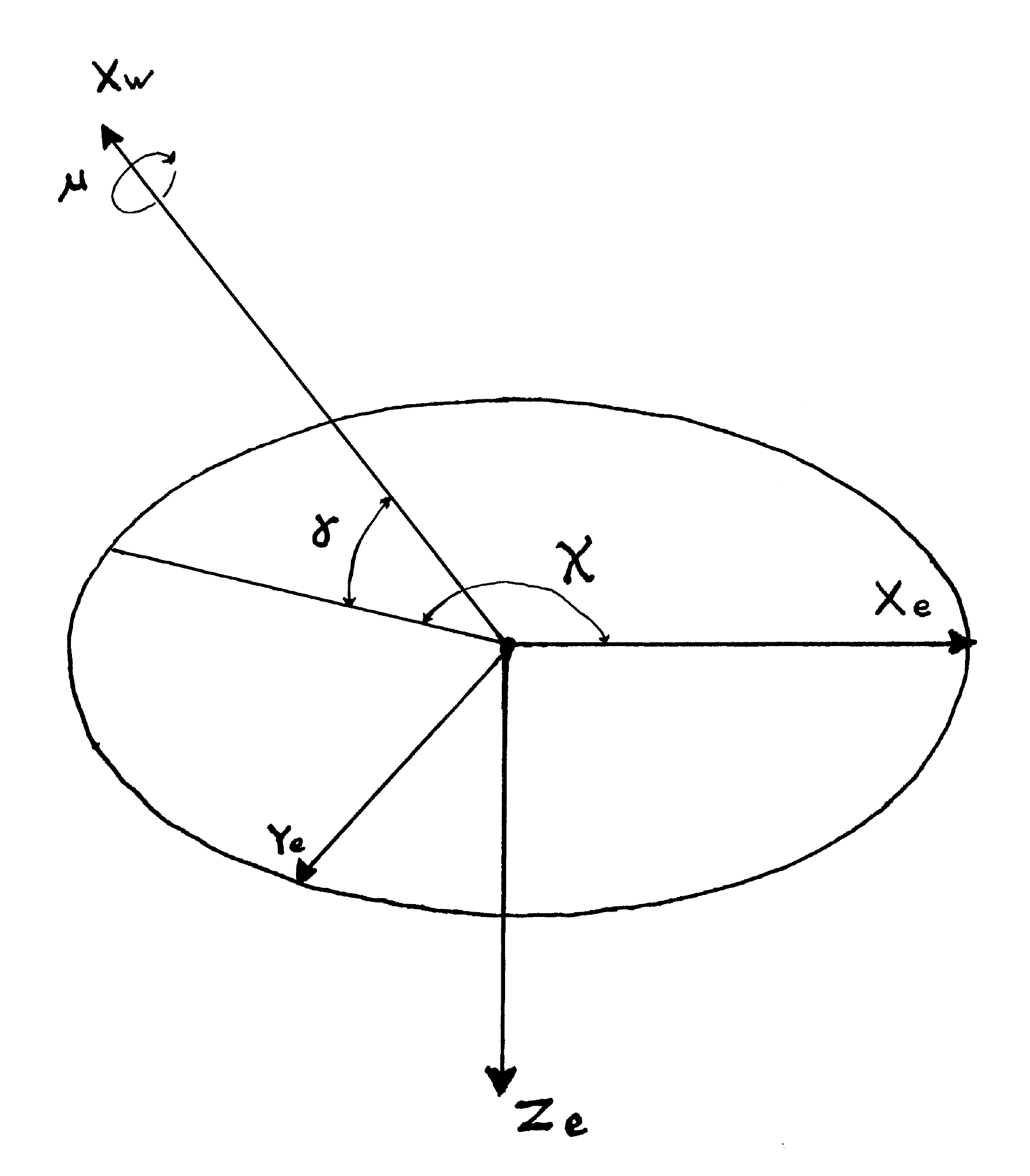}{a)}\hss
\includegraphicsh[height=5cm]{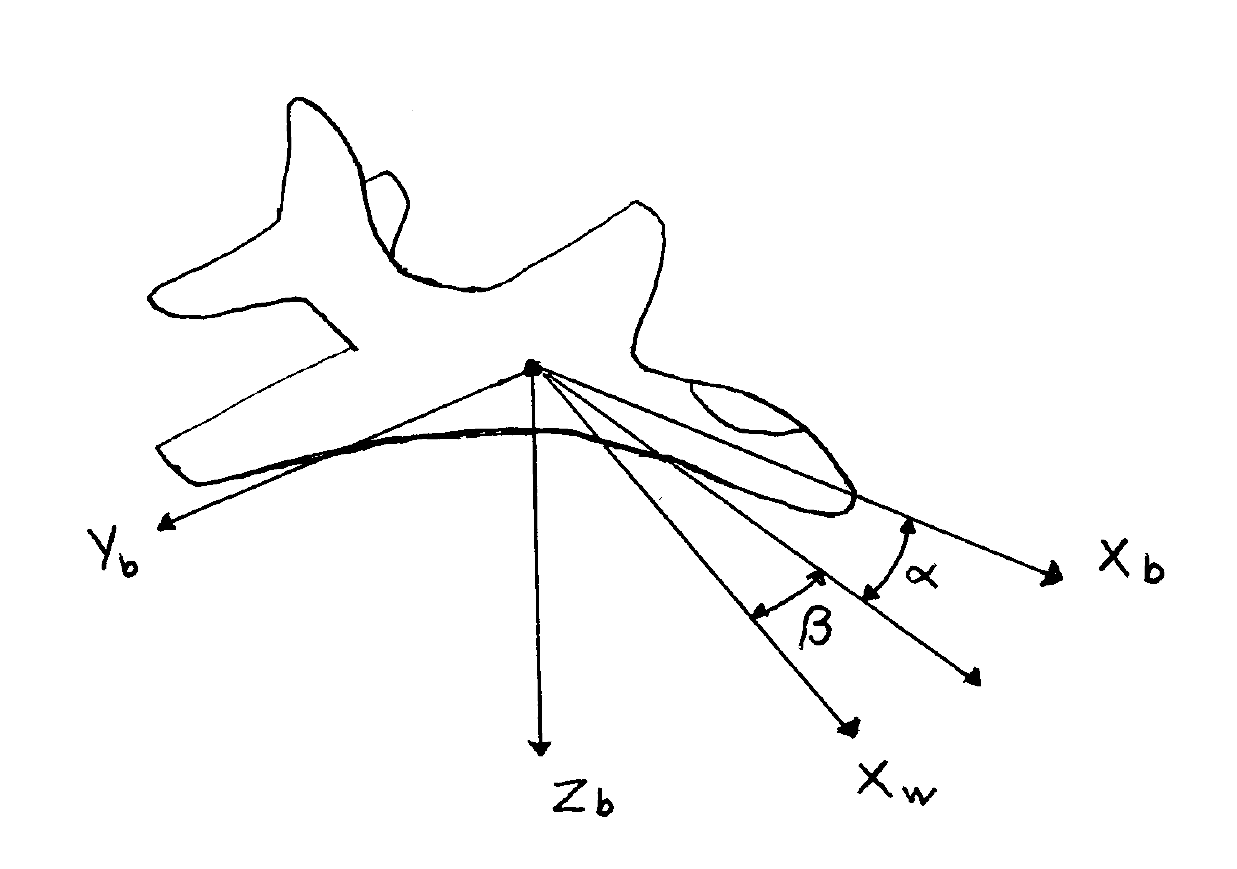}{b)}\hss}
    \caption{a)~Wind frame and b)~From wind to body frame.} 
\label{fig::wind_frame}
\end{center}
\end{figure}

The wind frame, with origin the center of gravity of the aircraft has
an axis $x_{w}$, in the direction of the velocity of the aircraft, the
axis $z_{w}$ being in the plane of symmetry of the aircraft.  The
Euler angles\footnote{\PATCH{More precisely, such angles are known as
    Tait-Bryan angles.}} that define the orientation of the wind frame
in the earth frame are denoted $\chi(t),\gamma(t),\mu(t)$, and are
respectively the \textit{aerodynamic azimuth} or \textit{heading
  angle}, the \textit{flight path angle} and the \textit{aerodynamic
  bank angle}, positive if the port side of the aircraft is higher
than the starboard side. See figure~\ref{fig::wind_frame}~(a). We go
from earth referential to wind referential using first a rotation with
respect to $z$ axis by the \textit{heading angle} $\chi$, then a
rotation with respect to $y$ axis by the \textit{flight path angle}
$\gamma$, and last a rotation with respect to $x$ axis by the
\textit{bank angle} $\mu$.

%\begin{figure}[h!]
%\begin{center}
%  \hbox to\hsize{\hss\includegraphics[height=5cm]{./Earth2wind.pdf}{a)}\hss
%\includegraphics[height=5cm]{./Wind2Body_ft.pdf}{b)}\hss}
%    \caption{a)~Wind frame and b)~From wind to body frame.} 
%\label{fig::wind_frame}
%\end{center}
%\end{figure}

The orientation of the wind frame with respect to the body frame is
defined by two angles: the \textit{angle of attack} $\alpha(t)$ and the
\textit{sideslip angle} $\beta(t)$, which is positive when the wind is on the
starboard side of the aircraft, as in figure~\ref{fig::wind_frame}~(b).
We go from the wind referential to the body referential using first a
rotation with respect to $z$ axis by the side slip angle
$\beta$ and then a rotation with respect to $y$ axis by the
angle of attack $\alpha$.

%%%%%%%%%%%%%%%%%%%%% 
\subsubsection{Dynamics}

In the sequel, we shall write $p(t),q(t),r(t)$ the coordinates of the
rotation vector of the body frame with respect to the the earth frame
expressed \PATCH{in a frame attached to a Galilean referential,} that
coincides at time $t$ with the body frame, and $L(t)$, $M(t)$, $N(t)$
the corresponding torques. In the same way $X(t)$, $Y(t)$ and $Z(t)$
denote the forces applied on the aircraft, expressed \PATCH{in a frame
attached to a Galilean referential,} that coincide at time $t$ with the
the wind referential.

\subsubsection{Aircraft geometry}\label{subsubsec:aircr_geom}

The mass of the aircraft is denoted by $m$, $S$ is the surface of the
wings. In the body frame, we assume that the
aircraft is symmetrical with respect to the $xz$-plane, so that the
tensor of inertia has the following form:

 \begin{equation}
 J := 
\left[\begin{matrix}I_{xx} & 0 & - I_{xz} \\
0 & I_{yy} & 0 \\
- I_{xz} & 0 & I_{zz}
\end{matrix}\right]. 
\end{equation}
 
In the standard equations~\eqref{eq::GNA}, we also need $a$, that
stands for the \emph{wing span} and $b$ for the \emph{mean aerodynamic
  chord}. 

\subsubsection{Forces and torques}\label{subsubsec::forces}

The force $(X,Y,Z)$ in the wind frame and the torque $(L,M,N)$ are
expressed by the following formulas:

\begin{subequations}
\begin{align}
\label{eq:X}
X&=F(t)\cos(\alpha+\epsilon)\cos(\beta(t))-\frac{\rho}{2}SV(t)^{2}C_{x}- gm\sin{\left(\gamma(t) \right)};\\
\label{eq:Y}
Y&=F(t)\cos(\alpha+\epsilon)\sin(\beta(t))+\frac{\rho}{2}SV(t)^{2}C_{y}+gm\cos(\gamma(t))\sin(\mu(t));\\
\label{eq:Z}
Z&=-F(t)\sin(\alpha+\epsilon)-\frac{\rho}{2}SV(t)^{2}C_{z}+gm\cos(\gamma(t))\cos(\mu(t));\\
\label{eq:L}
L&=-y_{p}\sin(\epsilon)(F_{1}(t)-F_{2}(t))+\frac{\rho}{2}SV(t)^{2}bC_{l};\\
\label{eq:M}
M&=\frac{\rho}{2}SV(t)^{2}cC_{m};\\
\label{eq:N}
N&=y_{p}\cos(\epsilon)(F_{1}(t)-F_{2}(t))+\frac{\rho}{2}SV(t)^{2}cC_{n}.
\end{align}
\end{subequations}
The angle $\epsilon$ is related to the lack of parallelism of the
reactors with respect to the $xy-$plane in the body frame and is small.

The aerodynamic coefficients $C_x,C_y,C_z,C_{l},C_{m},C_{n}$ depend on
$\alpha$ and $\beta$ and also on the angular speeds $p$, $q$, $r$ and
the controls are virtual angles\footnote{The angles are not physical
  angles, but rather measures related to some physical angles and that
  \PATCH{are} 
  calibrated by the aircraft producer.} $\delta_{l}$, $\delta_{m}$ and
$\delta_{n}$, that respectively express the positions of the ailerons,
elevator and rudder.

\begin{remark}\label{rem:diff-thrust}
One may use $\eta=(F_{1}-F_{2})/(F_{1}+F_{2})$ as an alternative
control instead of $\delta_{n}$ in case of rudder jam, \PATCH{where $F_1$ and
$F_2$ represent the thrust of the right and left engines, the sum of
which is equal to the total thrust $F$.}
\end{remark}

\subsubsection{Equations}\label{subsubsec:equations}
 
Following Martin~\cite{Martin-PHD,Martin-96}, the dynamics of the system is
modeled by the following set of explicit differential equations
(\ref{eq::aircraft1:x}--\ref{eq::aircraft1:mu},\ref{eq::aircraft2}):
{\small
\begin{subequations}
\begin{align}
\label{eq::aircraft1:x}
  \frac{d}{d t} x(t)  & =  V(t) \cos(\chi(t)) \cos(\gamma(t));\\[2mm]
\label{eq::aircraft1:y}
\frac{d}{d t} y(t)  & =  V(t) \sin(\chi(t)) \cos(\gamma(t)); \\[2mm]
\label{eq::aircraft1:z}
\frac{d}{d t} z(t)  & =  - V(t) \sin(\gamma(t));\\[2mm]
\label{eq::aircraft1:V}
\frac{d}{d t} V(t) & =  \frac{X}{m} ;\\[2mm]  
\label{eq::aircraft1:gamma}
\frac{d}{d t} \gamma(t)& = 
                    -\frac{Y\sin(\mu(t))+Z\cos(\mu(t))}{mV(t)};\\[2mm]
\label{eq::aircraft1:chi}
\frac{d}{d t} \chi(t)& =
              \frac{Y\cos(\mu(t))-Z\sin(\mu(t))}{\cos(\gamma(t))mV(t)};\\[2mm]
\label{eq::aircraft1:alpha}
\frac{d}{d t} \alpha(t) & = \begin{array}[t]{l}
               \frac{1}{\cos(\beta(t))}
              (-p\cos(\alpha(t))\sin(\beta(t))
                     +q\cos(\beta(t))\\
                     -r\sin(\alpha(t))\sin(\beta(t))
                     +\frac{Z}{mV(t)});\end{array}\\[2mm]
\label{eq::aircraft1:beta}
\frac{d}{d t} \beta(t)& = 
                    +p\sin(\alpha(t))
                    -r\cos(\alpha(t))
                    +\frac{Y}{mV(t)};\\[2mm]
\label{eq::aircraft1:mu}
\frac{d}{d t} \mu(t)& =  \begin{array}[t]{l}\frac{1}{\cos(\beta(t))}
                    (p\cos(\alpha(t))+r\sin(\alpha(t))\\
            +\frac{1}{mV(t)}(Y\cos(\mu(t))\tan(\gamma(t))\cos(\beta(t))\\
            -Z(\sin(\mu(t))\tan(\gamma(t))\cos(\beta(t))
            +sin(\beta(t)))));
\end{array}
\end{align}
\end{subequations}

}

\begin{equation}
\label{eq::aircraft2}
\left(\begin{array}{l}\frac{d}{d t} p(t)\\
  \frac{d}{d t} q(t)\\
  \frac{d}{d t} r(t)
\end{array}\right) =
J^{-1}
\left(\begin{array}{l}
 (I_{yy}-I_{zz})qr+I_{xz}pq+L\\
 (I_{zz}-I_{xx})pr+I_{xz}(r^{2}-p^{2})+M\\
  (I_{xx}-I_{yy})pq-I_{xz}rq+N
\end{array}\right).
\end{equation}
In the last expressions, the terms depending on gravity have been
incorporated to the expressions $X$, $Y$ and $Z$, as
in~\cite{Martin-96}.

We notice with Martin that this set of equations imply
$\cos(\beta)\cos(\gamma)V\neq0$. The nonvanishing of $V$ and $\cos(\beta)$
seems granted in most situations; the vanishing of $V$ may occur
with aircrafts equipped with vectorial thrust, which means a larger set
of controls, that we won't consider here. The vanishing of
$\cos(\gamma)$ can occur with loopings etc.\ and would require the
choice of a second chart with other sets of Euler angles. 
\PATCH{According to \eqref{eq::aircraft1:z}, the value of $z$ is
  negative, as in fig.~\ref{fig::earth_frame}~a), the axis $z$ points
  to the ground.}

\subsection{The GNA model}\label{sub:GNA}

In the last equations, $\rho$ can depend on $z$, as the air density vary with
altitude. The expression of $C_{x}$ and $C_{z}$ could also depend on $z$
to take in account ground effect. These expressions that depend on $\alpha$,
$\beta$, $p$, $q$, $r$,
$\delta_{l}$, $\delta_{m}$ and $\delta_{n}$,
should also depend on the Mach number, but most available
formulas are given for a limited speed range and the dependency on $V$
is limited to the $V^{2}$ term in factor.
In the literature, the available expressions are often partial or
limited to linear approximations. McLean~\cite{McLean} provides such
data for various types of aircrafts; for different speed and flight
conditions, including landing conditions with gears and flaps
extended.

We have chosen here to use the Generic Nonlinear Aerodynamic (GNA)
subsonic models, given by Grauer
and Morelli~\cite{Grauer-Morelli} that cover a wider range of values,
given in the following table.
\begin{figure}[h!]
  \begin{center}
  \begin{tabular}{|l|l|l|}\hline
    $-4\Deg\le\alpha\le 30\Deg$&$-20\Deg\le\beta\le20\Deg$,& \hfill\\
    \hline
$-100\Deg/\s\le p\le100\Deg/\s$& $-50\Deg/\s\le q\le50\Deg/\s$
& $-50\Deg/\s\le r\le50\Deg/\s$\\
\hline
$-10\Deg\le\delta_{l}\le10\Deg$& $-20\Deg\le\delta_{m}\le10\Deg$&
$-30\Deg\le\delta_{n}\le30\Deg$\\
\hline
  \end{tabular}
\end{center}
  \caption{\label{eq:model_range} Range of values for the GNA model}
\end{figure}
Among the 8 aircrafts
in their database, we have made simulations with 4: fighters F-4
and F-16C, STOL utility aircraft DHC-6 Twin Otter and the sub-scale
model of a transport aircraft GTM (see~\cite{Hueschen}).

The GNA model depends on 45 coefficients:

{\small
\begin{equation}\label{eq::GNA}
\begin{array}{lll}
  C_{D} &=& \theta_{1}+\theta_{2}\alpha+\theta_{3}\alpha\tilde q
    +\theta_{4}\alpha\delta_{m}+\theta_{5}\alpha^{2}
    +\theta_{6}\alpha^{2}\tilde q+\theta_{7}\delta_{m}
    +\theta_{8}\alpha^{3}+\theta_{9}\alpha^{3}\tilde q
    +\theta_{10}\alpha^{4},\\
  C_{y} &=& \theta_{11}\beta+\theta_{12}\tilde p
    +\theta_{13}\tilde r +\theta_{14}\delta_{l}+\theta_{15}\delta_{n},\\
  C_{L} &=& \theta_{16}+\theta_{17}\alpha+\theta_{18}\tilde q
    +\theta_{19}\delta_{n}+\theta_{20}\alpha\tilde q+\theta_{21}\alpha^{2}
    +\theta_{22}\alpha^{3}+\theta_{23}\alpha^{4},\\
  C_{l} &=& \theta_{24}\beta+\theta_{25}\tilde p+\theta_{26}\tilde r
    +\theta_{27}\delta_{l}+\theta_{28}\delta_{n}\\
  C_{m} &=& \theta_{29}+\theta_{30}\alpha+\theta_{31}\tilde q
    +\theta_{32}\delta_{e}+\theta_{33}\alpha\tilde q
    +\theta_{34}\alpha^{2}\tilde q+\theta_{35}\alpha^{2}\delta_{e}
    +\theta_{36}\alpha^{3}\tilde q+\theta_{37}\alpha^{3}\delta_{e}
    +\theta_{38}\alpha^{4},\\
  C_{n} &=& \theta_{39}\beta+\theta_{40}\tilde p+\theta_{41}\tilde r
    +\theta_{42}\delta_{l}+\theta_{43}\delta_{n}
    +\theta_{44}\beta^{2}+\theta_{45}\beta^{3},
\end{array}
\end{equation}

}
\PATCH{where $\tilde p=ap$, $\tilde r=ar$, $\tilde q=bq$
(see~\ref{subsubsec:aircr_geom} for the meaning of $a$ and $b$).
The aerodynamic coefficient
$C_{D}$, $C_{Y}$ and $C_{L}$ are defined
by~\cite[(1)]{Grauer-Morelli}.
The coefficients $C_{x}$ $C_{y}$ and $C_{z}$ in the wind
frame are then given by the formulas:}
\begin{equation}\label{eq::LDtoXY}
\begin{array}{lll}
C_{x}&=&\cos(\beta)C_{D}-\sin(\beta)C_{Y},\\
C_{y}&=&\sin(\beta)C_{D}+\cos(\beta)C_{Y},\\
C_{z}&=&C_{L}
\end{array}
\end{equation}

\begin{definition}\label{def:simplified-model}
  The simplified model is obtained by replacing $p,q,r$,
  $\delta_{\ell}$, $\delta_{m}$ an $\delta_{n}$ by $0$ (or some known
  constants) in the expressions of $C_{x}$, $C_{y}$ and $C_{z}$.
\end{definition}

\subsection{Block triangular structure of the simplified model}\label{sub:parametrization}\label{subsec:block_triang}

\def\HB#1{\hbox to 5pt{\hss{\footnotesize $#1$}\hss}}
Using the simplified model, \PATCH{the order matrix, where $-\infty$
  terms do not appear for better readability, is the following:
$$
\begin{array}{l}
    \phantom{\eqref{eq::aircraft1:x}}
  \begin{array}{cccc||ccc||cccc||ccc||ccc}
    \phantom{x\hskip-0.48ptx}&\HB{x}&\HB{y}&\HB{z}&\HB{V}&\HB{\gamma}&\HB{\chi}&\HB{\alpha}&\HB{\beta}&\HB{\mu}&\HB{F}&\HB{p}&\HB{q}&\HB{r}&\HB{\delta_{\ell}}&\HB{\delta_{m}}&\HB{\delta_{n}}
  \end{array}\\
    \begin{array}{c}
    \eqref{eq::aircraft1:x}\\
    \eqref{eq::aircraft1:y}\\
    \eqref{eq::aircraft1:z}\\
     \hline\hline
    \eqref{eq::aircraft1:V}\\
    \eqref{eq::aircraft1:gamma}\\
    \eqref{eq::aircraft1:chi}\\
     \hline\hline
    \eqref{eq::aircraft1:alpha}\\
    \eqref{eq::aircraft1:beta}\\
    \eqref{eq::aircraft1:mu}\\
     \hline\hline
    (\ref{eq::aircraft2}\hbox{p})\\
    (\ref{eq::aircraft2}\hbox{q})\\
    (\ref{eq::aircraft2}\hbox{r})
  \end{array}  
 \left(
  \begin{array}{ccc||ccc||cccc||ccc||ccc}
    1& & &0&0&0& & & & & & & & & & \\
     &1& &0&0&0& & & & & & & & & & \\
     & &1&0&0& & & & & & & & & & & \\
     \hline\hline
     & &0&1& & &0&0& &0& & & & & & \\
     & &0&0&1& &0&0&0&0& & & & & & \\
     & &0&0&0&1&0&0&0&0& & & & & & \\
     \hline\hline
     & &0&0& & &1&0& &0&0&0&& & & \\
     & &0&0& & &0&1& &0&0& &0& & & \\
     & &0&0&0& &0&0&1&0&0& &0& & & \\
     \hline\hline
     & &0&0& &0& & & &0&1&0&0&0&0&0\\
     & &0&0& &0& & & &0&0&1&0&0&0&0\\
     & &0&0& &0& & & &0&0&0&1&0&0&0
  \end{array}
 \right).
\end{array}
$$
So, algorithm \textsc{\=O-test} provides the following partition,
proceeding as in sec.~\ref{subsec:suf_reg_cond}:}
with $\Xi_{0}:=\{x,y,z\}$,
$\Xi_{1}:=\{V,\gamma,\chi\}$, $\Xi_{2}:=\{\alpha,\beta,\mu,F\}$,
$\Xi_{3}:=\{p,q,r\}$ and
$\Xi_{4}:=\{\delta_{\ell},\delta_{m},\delta_{n}\}$ or
$\Xi_{4}:=\{\delta_{\ell},\delta_{m},\eta\}$ when differential thrust is
used (see~rem.~\ref{rem:diff-thrust}), we see that, for $1\le k\le 3$,
$\dot \Xi_{i}$ only depends on $\bigcup_{\kappa=0}^{k+1}\Xi_{\kappa}$,
so that the model is block triangular with $\Sigma_{1}$ corresponding to
(\ref{eq::aircraft1:x}--\ref{eq::aircraft1:z}), $\Sigma_{2}$
corresponding to (\ref{eq::aircraft1:V}--\ref{eq::aircraft1:chi}),
$\Sigma_{3}$ corresponding to
(\ref{eq::aircraft1:alpha}--\ref{eq::aircraft1:mu}) and $\Sigma_{4}$
to \eqref{eq::aircraft2}.

Simple computations show that
\begin{equation}\label{eq:rank1}
\left|\frac{\partial P}{\partial \xi}| P\in\Sigma_{1};\> \xi\in
  \Xi_{1}\right|=-V^{2}\cos(\gamma); \quad 
\left|\frac{\partial P}{\partial \xi}| P\in\Sigma_{3};\> \xi\in
  \Xi_{3}\right|=\frac{1}{\cos(\beta)},
\end{equation}
that do not vanish. In the same way,
\begin{equation}\label{eq:rank2}
\left|\frac{\partial P}{\partial \xi}| P\in\Sigma_{4};\> \xi\in
  \Xi_{3}\right|=\frac{\rho^{3}}{8}S^{3}V^{6}a^{2}b|J|^{-1}
  \left|\frac{\partial C_{i}}{\partial \delta_{j}}| (i,j)\in\{\ell,m, n\}^{2}\right|,
\end{equation}
does not vanish as the diagonal terms ${\partial
  C_{i}}/{\partial \delta_{i}}$, $i=\ell,m,n$, are much bigger than
the others.

So, to apply th.~\ref{th:suff_reg_cond}, we only need to consider the rank of the Jacobian matrix
$$
\left(\frac{\partial P}{\partial \xi}| P\in\Sigma_{2};\> \xi\in
  \Xi_{2}\right),
$$
which  is equal to the rank of the Jacobian matrix
\begin{equation}\label{eq::Delta}
   \Delta:= \left(\frac{\partial Q}{\partial \xi}|
      Q\in\{X,\sin(\mu)Y+\cos(\mu)Z,-\cos(\mu)Y+\sin(\mu)Z\};\> \xi\in
      \Xi_{2}\right).
\end{equation}

\begin{proposition}\label{prop:simpl-aicr-flat}
  The simplified aircraft is flat when the matrix $\Delta$ has full
  rank.

  Let $\Delta_{\xi}$ be the matrix $\Delta$ where the column
  corresponding to $\xi\in\Xi_{2}$ has been suppressed. If
  $|\Delta_{\xi}|\neq 0$ at some point, then $x,y,z,\xi$ is a regular
  flat output at that point.
  \end{proposition}
\begin{proof} This is a direct consequence of th.~\ref{th:suff_reg_cond}~ii).
\end{proof}

\begin{example}\label{ex::plane} We can associate a diffiety to our aircraft
  model, which is defined by $\R^{12}\times\left(\R^{\N}\right)^{4}$,
  and a derivation $\meth$ defined by $\meth:=\meth_{0}+\meth_{1}$, where
  $\meth_{1}$ is the trivial derivation on $\left(\R^{\N}\right)^{4}$:
  $\meth_{1}:=\sum_{z\in\{F,\delta_{l},\delta_{m},\delta_{n}\}}^{m}\sum_{k\in\N}
  z_{i}^{(k+1)}\partial/\partial z_{i}^{(k+1)}$ and $\meth_{0}$ is
  defined on $\R^{12}$ by the differential equations
  (\ref{eq::aircraft1:x}--\ref{eq::aircraft1:mu},\ref{eq::aircraft2}):
$$
\delta_{0}:= V(t) \cos{\left(\chi(t) \right)} \cos{\left(\gamma(t)
  \right)}\frac{\partial}{\partial x}+V(t) \sin{\left(\chi(t) \right)} \cos{\left(\gamma(t)
  \right)}\frac{\partial}{\partial y}+\cdots
$$

In practice, the diffiety is defined by a smaller open set, because of
the bounds that exist on the values of the variables. The values of
the controls are bounded and one wishes to restrict the values of the
angle of attack $\alpha$ or side-slip angle $\beta$ for safety
reasons. The maximal values for the GNA model are given below in
table~\ref{eq:model_range}. Other limitations must be included, such
as the maximal value of the thrust. The speed $V$ should also be
greater than the stalling speed (see~\ref{subsec:stalling}).
\end{example}

We will now consider more closely the possible choices of flat outputs.

\subsection{Choices of flat outputs}\label{subsec:choices-FO}

\subsubsection{The side-slip angle choice}\label{subsubsec:classical-FO}

Martin~\cite{Martin-PHD,Martin-96} has used the set of flat outputs:
$x,y,z,\beta$. We need to explicit under which condition such a flat
output may be chosen, \textit{i.e.} when the Jacobian determinant
$$
\Delta_{\beta}=\left|
\begin{array}{ccc}
\frac{\partial X}{\partial \alpha}&
\frac{\partial X}{\partial \mu}&
\frac{\partial X}{\partial F}\\
\frac{\partial(\sin\mu\,Y+\cos\mu\, Z)}{\partial \alpha}
&\frac{\partial(\sin\mu\,Y+\cos\mu\, Z)}{\partial \mu}
&\frac{\partial(\sin\mu\,Y+\cos\mu\, Z)}{\partial F}\\
\frac{\partial(\cos\mu\,Y-\sin\mu\, Z)}{\partial \alpha}
&\frac{\partial(\cos\mu\,Y-\sin\mu\, Z)}{\partial \mu}
&\frac{\partial(\cos\mu\,Y-\sin\mu\, Z)}{\partial F}
\end{array}
\right|
$$ does not vanish, according to prop.~\ref{prop:simpl-aicr-flat}.
First, we remark, following Martin~\cite[p.~80]{Martin-PHD} that when
$Y=gm\cos(\gamma)\sin(\mu)$ and $Z=gm\sin(\gamma)\sin(\mu)$,
\textit{i.e.} when the lift is zero, $\sin\mu\,Y+\cos\mu\,
Z=gm\cos(\gamma)$ and $\cos\mu\,Y-\sin\mu\, Z=0$, so that the second
column of $\Delta_{\beta}$ is zero and the determinant vanishes. This
means that $0$-g flight trajectories are singular for this flat
output. On the other hand, when $\beta$ and $\alpha$ are close to $0$,
which is the case in straight and level flight, easy computations
using eq.~(\ref{eq:X}--\ref{eq:Z}) allow Martin to conclude that
$$
|\Delta_{\beta}| \approx-Z\left(\frac{\rho}{2}SV^{2}\frac{\partial
  C_{z}}{\partial\alpha}+F\right)\gg 0.
$$

To go further, one may use the expression of $X$ eq.~\eqref{eq:X} and
deduce from it
\begin{equation}\label{eq:F}
  F=\frac{X+\frac{\rho}{2}SV^{2}C_{x}+gm\sin(\gamma)}
  {\cos(\alpha+\epsilon)\cos(\beta)}, 
\end{equation}
assuming $\cos(\alpha+\epsilon)\cos(\beta)$.
Substituting this expression in $Y$ and $Z$, we define $\tilde Y$ and $\tilde
Z$ and further define $\hat
Y:=\cos(\mu)\tilde Y-\sin(\mu)\tilde Z$ and $\hat Z:=\sin(\mu)\tilde
Y+\cos(\mu)\tilde Z$. Then, $|\Delta_{\beta}|\neq0$ when
\begin{equation}\label{eq::cond_beta}
\left|\begin{array}{cc}
\frac{\partial\hat Y}{\partial\alpha}&
\frac{\partial\hat Y}{\partial\mu}\\
\frac{\partial\hat Z}{\partial\alpha}&
\frac{\partial\hat Z}{\partial\mu}\end{array}\right|\neq0.
\end{equation}
The main interest of this choice is to be able to
impose easily $\beta=0$, which is almost always required.

\subsubsection{The bank angle choice}\label{subsec:Bank-FO}

As the angle $\mu$ is known, we may compute from $\Xi_{1}'$ and
$\Xi_{1}''$ the values $X$, $Y$ and $Z$. So, singularities for this
flat outputs are such that
\begin{equation}\label{eq::stalling_cond}
|\Delta_{\mu}|=\left|
\begin{array}{ccc}
\frac{\partial X}{\partial \alpha}&
\frac{\partial X}{\partial \beta}&
\frac{\partial X}{\partial F}\\
\frac{\partial Y}{\partial \alpha}
&\frac{\partial Y}{\partial \beta}
&\frac{\partial Y}{\partial F}\\
\frac{\partial Z}{\partial \alpha}
&\frac{\partial Z}{\partial \beta}
&\frac{\partial Z}{\partial F}
\end{array}
\right|\neq0
\end{equation}
Using $\tilde Y$ and $\tilde Z$, as defined in
subsec.~\ref{subsubsec:classical-FO}, we see that this is equivalent
to
\begin{equation}\label{eq::cond_mu}
\left|\begin{array}{cc}\frac{\partial\tilde Y}{\partial\alpha}&
\frac{\partial\tilde Y}{\partial\beta}\\
\frac{\partial\tilde Z}{\partial\alpha}&
\frac{\partial\tilde Z}{\partial\beta}\end{array}\right|\neq0.
\end{equation}
When $\beta$ is $0$, $\partial\tilde Z/\partial\beta$ is also $0$, due
to the aircraft symmetry. Using the GNA model (see~\ref{sub:GNA}), we
have $\partial C_{x}/\beta=0$ and $\partial C_{z}/\partial\beta=0$, so
that $\partial\tilde Z/\partial\beta=0$. The value of the
determinant~\eqref{eq::cond_mu} is then
\begin{equation}\label{eq::cond_mu_bis}
-\frac{\partial\tilde Z}{\partial\alpha}\frac{\partial\tilde Y}{\partial\beta}.
\end{equation}
For most aircrafts, $\partial
C_{y}/\partial \beta$ is negative at $\beta=0$, with values in the
range $[-1., -0.5]$.
Delta wing aircrafts seem to be a exception, with smaller absolute
values ($-0.014$ for the X-31) or even negative ones ($+0.099$ for the
F-16XL). It seems granted that for regular transport planes, $\partial
C_{y}/\partial \beta$ is negative, so that the determinant vanishes
only when $\partial\tilde Z/\partial\alpha$ is $0$. We will see
in~\ref{subsec:stalling} that this may be interpreted as stalling
condition and that the vanishing of~\eqref{eq::cond_mu_bis} on a
trajectory with constant controls means that the points of this
trajectory are flat singularities, so that no other flat outputs could
work. 

This choice is the best to impose $\mu=0$ and is natural for decrabe
maneuver, that is when landing with a lateral wind, which implies
$\beta\neq0$. We then need to maintain
$\mu$ close to $0$ to avoid the wings hitting the runway.

It is also a good choice when
$Y=Z=0$, a situation that may be encountered in aerobatics or when
training for space condition with $0$-g flights (see
subsec~\ref{subsec:parabolic-F16C}). The choices $\beta$ and $\mu$ are
compared in~\cite[7.1]{Ollivier2022}
with the simulation of a twin otter
flying with one engine.

\subsubsection{The thrust choice}\label{subsubsec:thrust-FO}

The choice of thrust $F$ has one main interest: to set $F=0$ and
consider the case of a aircraft having lost all its engines.  The
aircraft must land by gliding when all engines are lost. This is a
rare situation, but many successful examples are known, including the
famous US Airways Flight 1549~\cite{NTSB2010}.  The singularities of
this flat output are such that
$$
|\Delta_{F}|=\left|
\begin{array}{ccc}
\frac{\partial X}{\partial \alpha}&
\frac{\partial X}{\partial \mu}&
\frac{\partial X}{\partial \beta}\\
\frac{\partial(\sin\mu\,Y+\cos\mu\, Z)}{\partial \alpha}
&\frac{\partial(\sin\mu\,Y+\cos\mu\, Z)}{\partial \mu}
&\frac{\partial(\sin\mu\,Y+\cos\mu\, Z)}{\partial \beta}\\
\frac{\partial(\cos\mu\,Y-\sin\mu\, Z)}{\partial \alpha}
&\frac{\partial(\cos\mu\,Y-\sin\mu\, Z)}{\partial \mu}
&\frac{\partial(\cos\mu\,Y-\sin\mu\, Z)}{\partial \beta}
\end{array}
\right|
$$
vanishes. When $F=0$, by eq.~(\ref{eq:X}--\ref{eq:Z}), the
vanishing of $\Delta_{F}$ is
equivalent to
$$
\left|
\begin{array}{ccc}
\frac{\partial C_{x}}{\partial \alpha}&
0&
\frac{\partial C_{x}}{\partial \beta}\\
\frac{\partial(\cos\mu\, C_{z})}{\partial \alpha}
&\frac{\partial(\sin\mu\,C_{y}+\cos\mu\, C_{z})}{\partial \mu}
&\frac{\partial(\sin\mu\,C_{y})}{\partial \beta}\\
\frac{\partial(-\sin\mu\, C_{z})}{\partial \alpha}
&\frac{\partial(\cos\mu\,C_{y}-\sin\mu\, C_{z})}{\partial \mu}
&\frac{\partial(\cos\mu\,C_{y})}{\partial \beta}
\end{array}
\right|=\frac{\partial C_{x}}{\partial \alpha}\left(C_{y}\frac{\partial C_{y}}{\partial\beta}-C_{z}\frac{\partial C_{x}}{\partial \beta}\right).
$$ When $\beta$ vanishes, $C_{y}$ and $\partial C_{z}/\partial\beta$
also vanish, due to the aircraft symmetry with respect to the
$xz$-plane. So, we need have $\beta\neq0$ to use those flat outputs.
Using the GNA model, $C_{x}$ and $C_{z}$ 
depend only on $\alpha$ and $C_{y}$ depends linearly on $\beta$.

In the case of a gliding aircraft, situations with $\beta\neq0$ could
indeed be useful to achieve the forward slip maneuver.
When the aircraft is too high, combining nonzero $\beta$ and
$\mu$ precisely allows a fast descent while keeping a moderate
speed. This is very useful when gliding, as there is no option for a
go around when approaching the landing strip too high or too
fast. This maneuver was performed with success by the ``Gimli
Glider''~\cite{Gimli}, Air Canada Flight 143, that ran out of fuel on
July 23, 1983, which could land safely in Gimli former \PATCH{Air Force
base.}  A simulation of the forward slip may be found
in~\cite[7.2]{Ollivier2022}.

\subsubsection{Other sets of flat outputs}\label{subsec:Other-FO}

Among the other possible choices for completing the
set $\Xi_{1}$ in order to get flat outputs, the interest of $\alpha$
seems mostly academic. Indeed,
$$
\Delta_{\alpha}=\left|
\begin{array}{ccc}
\frac{\partial X}{\partial \beta}&
\frac{\partial X}{\partial \mu}&
\frac{\partial X}{\partial F}\\
\frac{\partial(\sin\mu\,Y+\cos\mu\, Z)}{\partial \beta}
&\frac{\partial(\sin\mu\,Y+\cos\mu\, Z)}{\partial \mu}
&\frac{\partial(\sin\mu\,Y+\cos\mu\, Z)}{\partial F}\\
\frac{\partial(\cos\mu\,Y-\sin\mu\, Z)}{\partial \beta}
&\frac{\partial(\cos\mu\,Y-\sin\mu\, Z)}{\partial \mu}
&\frac{\partial(\cos\mu\,Y-\sin\mu\, Z)}{\partial F}
\end{array}
\right|
$$

Then, $|\Delta_{\alpha}|=0$ when
\begin{equation}\label{eq::cond_alpha}
\left|\begin{array}{cc}
\frac{\partial\hat Y}{\partial\beta}&
\frac{\partial\hat Y}{\partial\mu}\\
\frac{\partial\hat Z}{\partial\beta}&
\frac{\partial\hat Z}{\partial\mu}\end{array}\right|=0.
\end{equation}
Easy computations show that it is the case when $\mu=\beta=0$, so that
$\alpha$ is not a suitable alternative input near stalling conditions,
except possibly when $\mu\neq0$. But mostly, our aerodynamic model cannot
fully reflect the real behavior near stalling so that it seems safe
to exclude stalling a.o.a from the domain of definition.
\bigskip

One may also consider time-varying
expressions, e.g. linear combinations of $\beta$ and $\mu$, to
smoothly go from one choice to another.

\subsection{Stalling conditions}\label{subsec:stalling}

It is known that the lift of a wing reaches a maximum at a critical
angle of attack, due to flow separation. This phenomenon can be hardly
reversible and creates a sudden drop of the lift force $Z$ from its peack
value. Our mathematical model is too poor to fully reflect such
behavior, but a maximum for the lift can still be computed. 

We need to take also in account the contribution of the thrust in the
expression of $Z$, and simple computations show that the critical
angle of attack corresponds in our setting to a maximum of $\tilde Z$,
that is $\partial\tilde Z/\partial\alpha=0$, which corresponds to the
singularity for flat output $\mu$ already observed
above~\eqref{eq::cond_mu_bis}.

Three cases may appear with stalling:

\noindent 1) to reach an extremal value of $\tilde Z$, meaning that
$\partial\tilde Z/\partial\alpha$ vanishes;

\noindent 2) to reach the maximum thrust $F_{\mbox{max}}$ before
reaching a maximum of $\tilde Z$;

\noindent 3) reaching no maximum of $\tilde Z$ with an aircraft with trust/weight
ratio greater than $1$: in such a case, there is no stalling.

For horizontal straight line trajectories, we may compute the speed
$V$ and the thrust $F$, depending on $\alpha$, for $\beta=\mu=0$,
using the simplified model. We may also use the full model. As
$\alpha$, $\beta$ and $\mu$ are constants, $p=q=r=0$, which further
allows to express $\delta_{\ell}$, $\delta_{m}$ and $\delta_{n}$
depending on $\alpha$, so that $C_{\ell}=C_{m}=C_{n}=0$, which is easy
with the GNA model that is linear in those quantities.

\begin{figure}[h!]
\hbox to \hsize{\hss
\includegraphicsh[width=4.5cm]{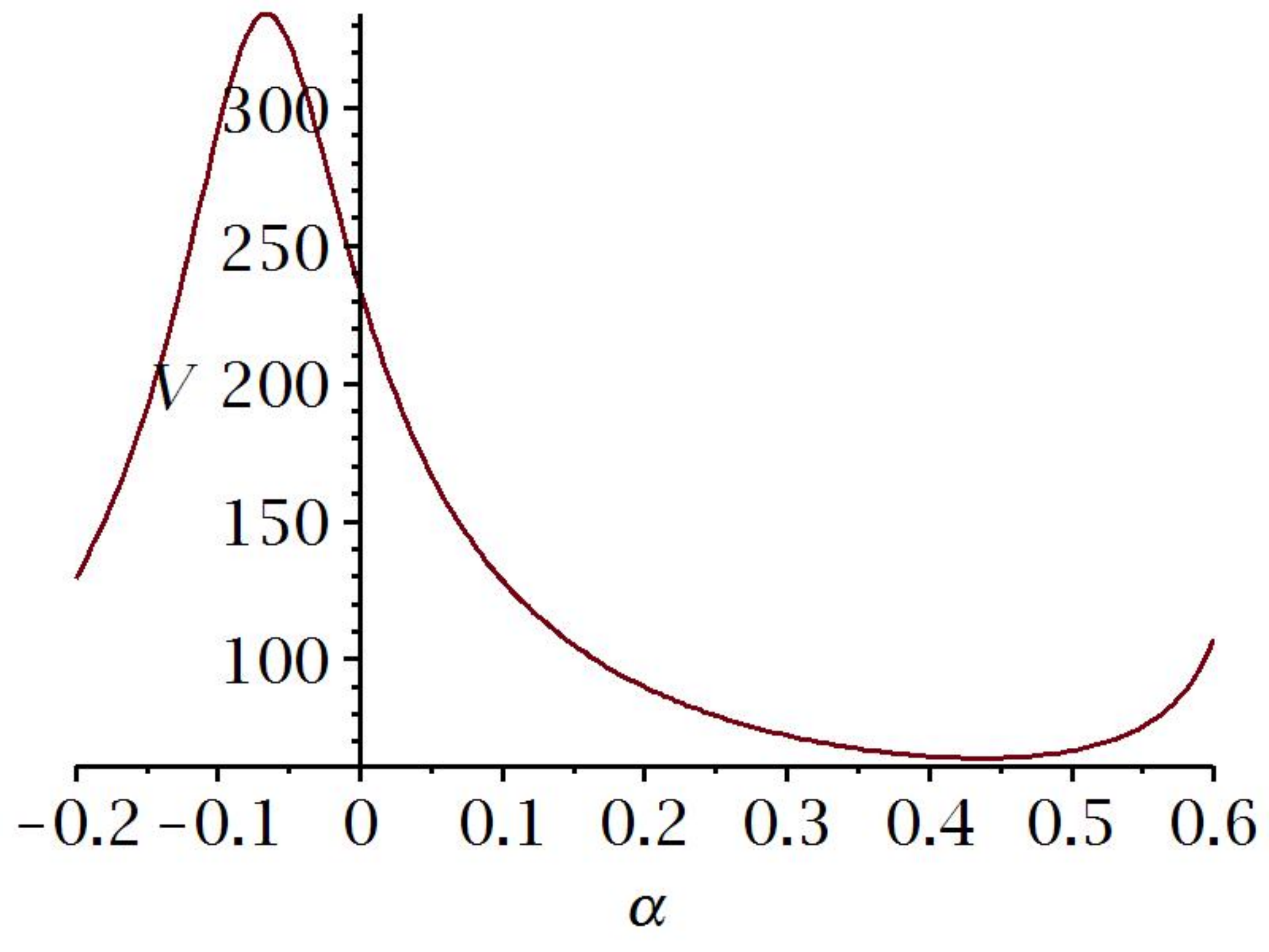}{a)}\hss
\includegraphicsh[width=4.5cm]{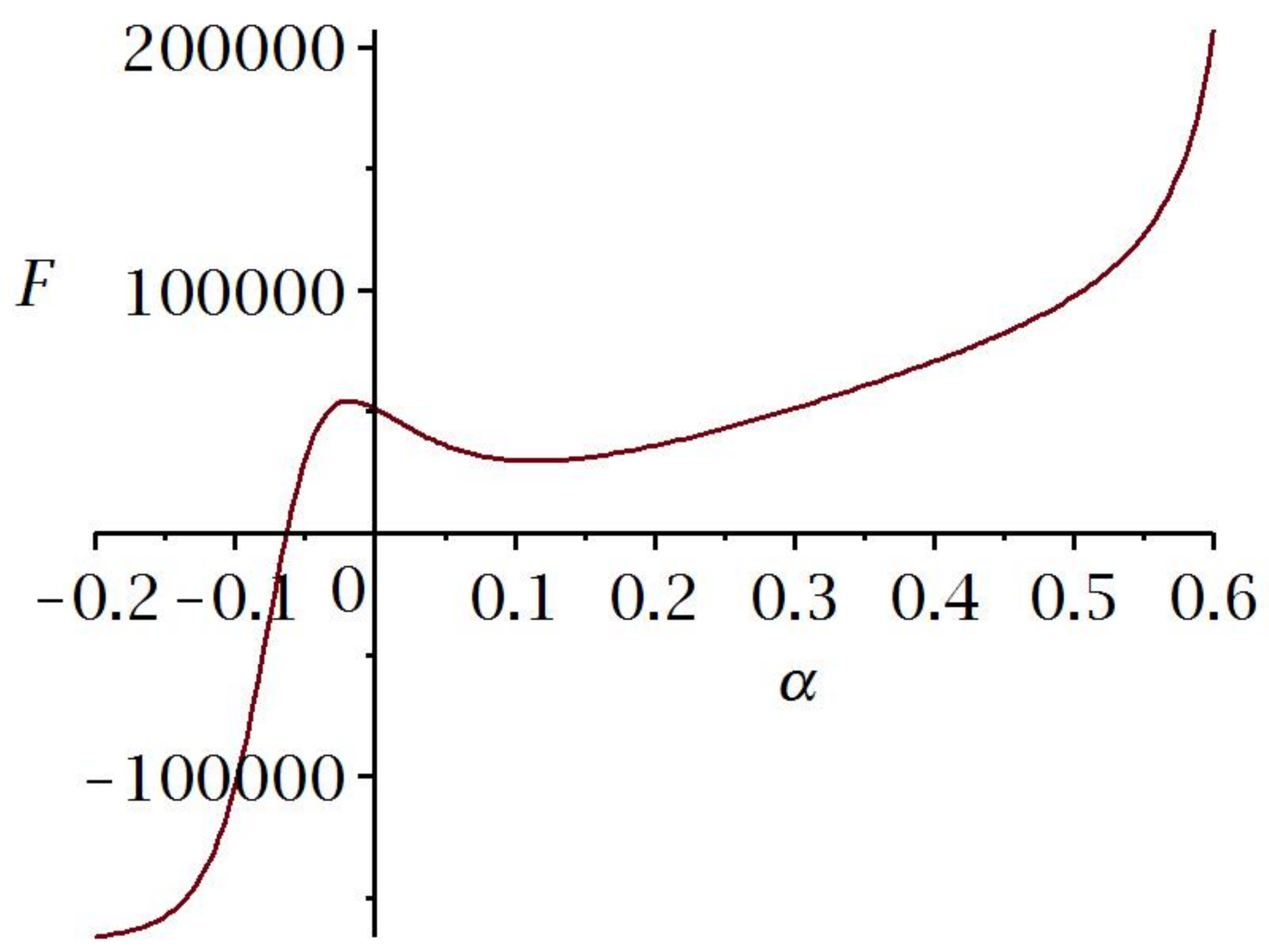}{b)}\hss
}
\hbox to \hsize{\hss
\includegraphicsh[width=4.5cm]{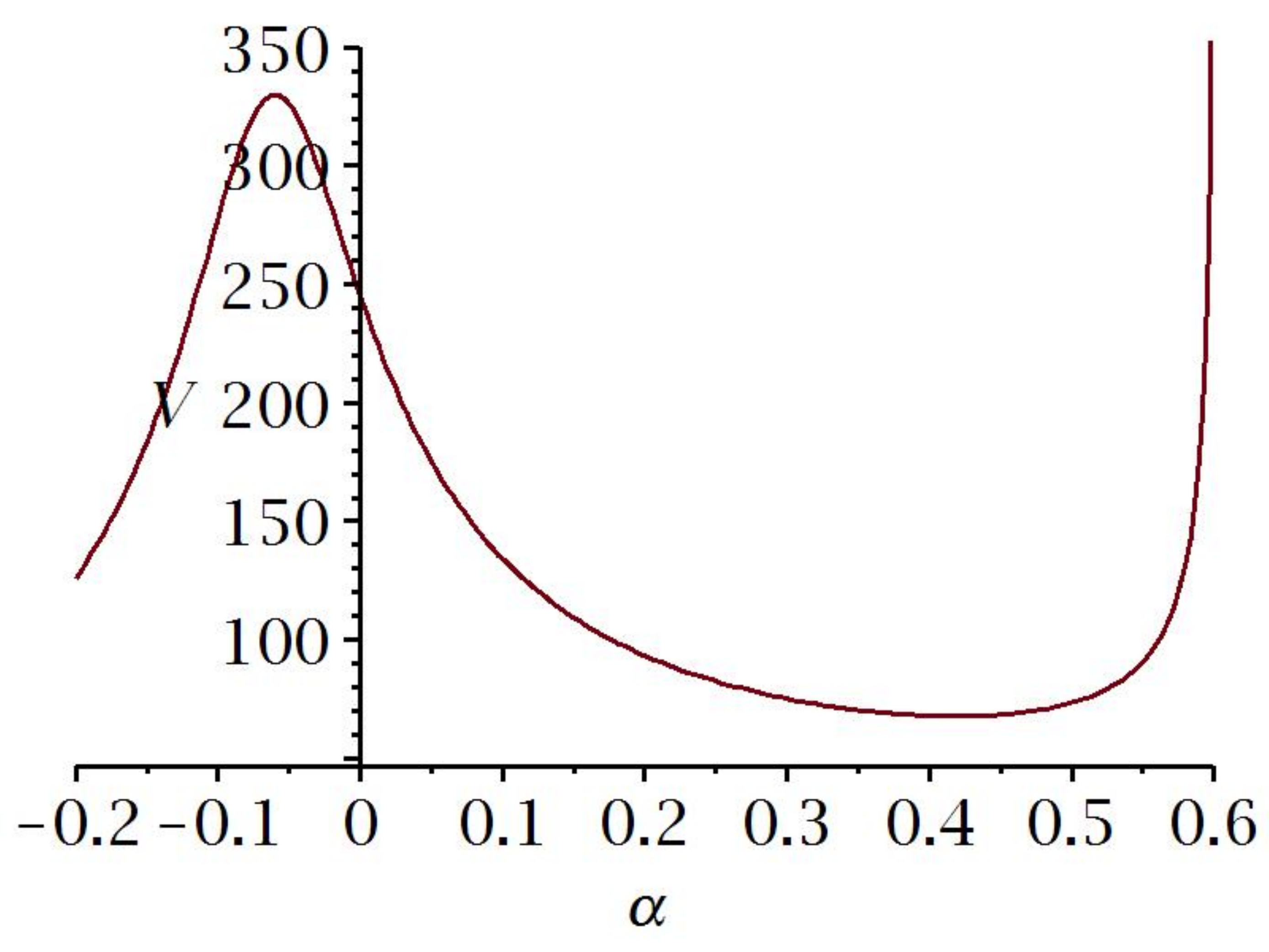}{c)}\hss
\includegraphicsh[width=4.5cm]{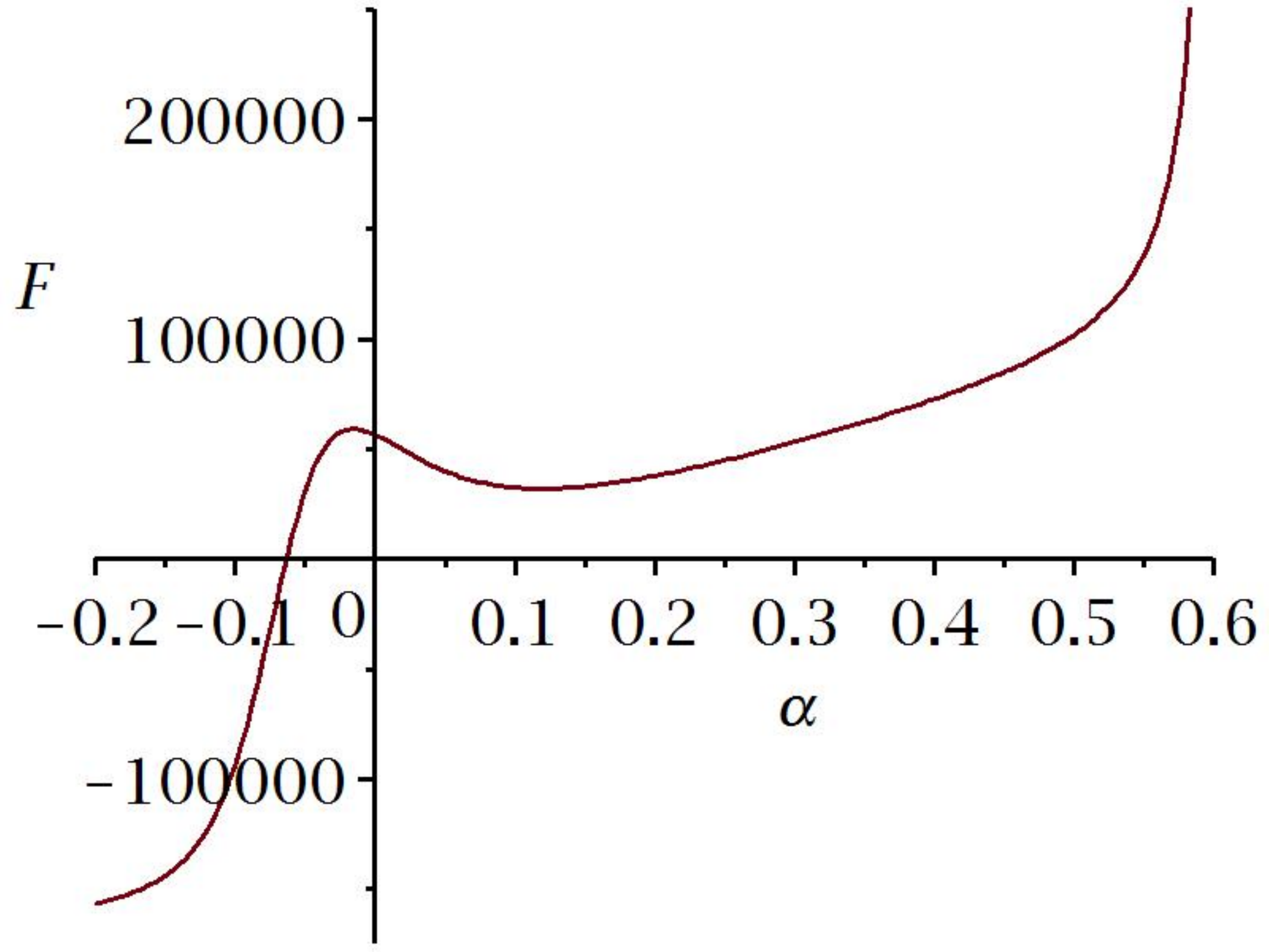}{d)}\hss
}
\caption{\label{fig::F4_F_V} F-4: values of $V$ and $F$ depending on
  $\alpha$. a) \&\ b)
  full generic nonlinear aerodynamic model; c) \&\ d)
  simplified model}\end{figure}

\textit{E.g.}, For the F4, setting the weight to
$38924\mbox{lb}$~\cite{Grauer-Morelli}, the evaluated stall speed,
angle of attack and thrusts are respectively $67.6789$m/s
($131.56$kn), $0.4200\rad$ ($24.07\Deg$) and $77.0436$ for the full
model and $64.0904$m/s, $0.4366\rad$ and $78.8806$ for the simplified
one, without thrust limitations. These thrust values are below the
thrust of J79-GE-17A engines of later versions ($79.38$kN with
afterburning). Assuming a maximal thrust of $2\times 71.8\mbox{kN}$,
that corresponds to the J79-GE-2 engines of the first production
aircrafts, the stall speed and angle of attack are $67.9835$m/s and
$0.3969\rad$ with the full model, $64.5515$m/s and $0.4057\rad$ with
the simplified one. The full model stall speed of about $132$kn agrees
more or so with the stall speed values of $146\mbox{KCAS}$ or
$148\mbox{KCAS}$, according to the models, computed with the NATOPS
manual~\cite[fig.~4.1 and 4.2]{NATOPS-F4} at $10000$~ft and below and
the computed a.o.a of $23\Deg$
with the indicated stall a.o.a of $27$ to $28$ ``units'', keeping in
mind that those units are not exactly degrees and that our
mathematical definition of stalling cannot fully match the actual
behavior. See fig.~\ref{fig::F4_F_V}.

\PATCH{
We notice in the figure that the curves become unrealistic outside the
range of values of the GNA model given in fig.~\ref{eq:model_range},
so that $F$ even takes negative values for $\alpha<0.065\rad$, which is
very close to the $-4\Deg$ limit.} 

%\begin{figure}[h!]
%\hbox to \hsize{\hss
%\includegraphics[width=3.8cm]{./plot_V_F4.pdf}{a)}
%\includegraphics[width=3.8cm]{./plot_F_F4.pdf}{b)}
%\includegraphics[width=3.8cm]{./plot_V_simp_F4.pdf}{c)}
%\includegraphics[width=3.8cm]{./plot_F_simp_F4.pdf}{d)}\hss
%}
%\caption{\label{fig::F4_F_V} F-4: values of $V$ and $F$ depending on
%  $\alpha$. a) \&\ b)
%  real model; c) \&\ d)
%  simplified model}\end{figure}

We have the following theorem, that shows that stalling condition
means that the system is not flat for some trajectories with $\beta=\mu=0$.

\begin{theorem} Let a trajectory be such that $\alpha$, $\beta$,
  $\mu$, $F$, $\gamma$ and $V$ are constants, with moreover
  $\beta=0$, $\alpha$ and $V$ respectively equal to the stall
  a.o.a. and stall speed. There is no flat point in this trajectory.
\end{theorem}
\begin{proof} Simple computations show that when $\beta=0$,
  $\alpha$ is equal to stalling a.o.a and $V$ the stall speed, the rank
  of the Jacobian matrix $(\partial P/\partial x\mid P\in\Sigma_{2},\>
  x\in\Xi_{3})$ is at most $3$.  This is then a straightforward
  consequence of th.~\ref{th:suff_sing}, noticing that $\chi'$ is the
  only term involving $\chi$ in the equations and that it appears
  linearly.
\end{proof}
\vfill\eject

\PATCH{
\subsection{Restricted model with 9 state functions}\label{subsec:restricted}

When foccussing on flight handling quality, one may restrict to $9$
state variables by neglecting $x$, $y$ and $\chi$ and suppressing
equations
(\ref{eq::aircraft1:x},\ref{eq::aircraft1:y}\ref{eq::aircraft1:chi}),
as those variables do not appear in the remaining $9$ equations. We
get then the order matrix:
$$
\begin{array}{l}
    \phantom{\eqref{eq::aircraft1:x}}
  \begin{array}{cc||cc||cccc||ccc||ccc}
    \phantom{x\hskip-0.48ptx}&\HB{z}&\HB{V}&\HB{\gamma}&\HB{\alpha}&\HB{\beta}&\HB{\mu}&\HB{F}&\HB{p}&\HB{q}&\HB{r}&\HB{\delta_{\ell}}&\HB{\delta_{m}}&\HB{\delta_{n}}
  \end{array}\\
    \begin{array}{c}
    \eqref{eq::aircraft1:z}\\
     \hline\hline
    \eqref{eq::aircraft1:V}\\
    \eqref{eq::aircraft1:gamma}\\
     \hline\hline
    \eqref{eq::aircraft1:alpha}\\
    \eqref{eq::aircraft1:beta}\\
    \eqref{eq::aircraft1:mu}\\
     \hline\hline
    (\ref{eq::aircraft2}\hbox{p})\\
    (\ref{eq::aircraft2}\hbox{q})\\
    (\ref{eq::aircraft2}\hbox{r})
  \end{array}  
 \left(
  \begin{array}{c||cc||cccc||ccc||ccc}
    1&0&0& & & & & & & & & & \\
     \hline\hline
    0&1& &0&0& &0& & & & & & \\
    0&0&1&0&0&0&0& & & & & & \\
     \hline\hline
    0&0& &1&0& &0&0&0& & & & \\
    0&0& &0&1& &0&0& &0& & & \\
    0&0&0&0&0&1&0&0& &0& & & \\
     \hline\hline
    0&0& &0& & &0&1&0&0&0&0&0\\
    0&0& &0& & &0&0&1&0&0&0&0\\
    0&0& &0& & &0&0&0&1&0&0&0
  \end{array}
 \right).
\end{array}
$$ Following the construction of subsec.~\ref{subsec:suf_reg_cond}, as
we did in subsec.~\ref{subsec:block_triang}, we define
$\bar\Xi_{0}:=\{z\}$, $\bar\Xi_{1}:=\{V,\gamma\}$,
$\bar\Xi_{2}:=\{\alpha,\beta,\mu,F)=\Xi_{2}$, $\bar\Xi_{3}:=\{p,q,r\}=\Xi_{3}$ and
$\bar\Xi_{4}:=\{\delta_{\ell},\delta_{m},\delta_{n}\}=\Xi_{4}$, we see that, for
$1\le k\le 3$, $\dot \Xi_{i}$ only depends on
$\bigcup_{\kappa=0}^{k+1}\Xi_{\kappa}$, so that this new model is also
block
triangular with $\bar\Sigma_{1}$ corresponding to
(\ref{eq::aircraft1:z}), $\bar\Sigma_{2}$
corresponding to (\ref{eq::aircraft1:V}--\ref{eq::aircraft1:gamma}),
$\Sigma_{3}$ corresponding to
(\ref{eq::aircraft1:alpha}--\ref{eq::aircraft1:mu}) and $\Sigma_{4}$
to \eqref{eq::aircraft2}.  It is easily seen that the only Jacobian
matrix that may possibly not be of full rank is $(\partial P/\partial
\xi|(P,\xi\in\bar\Sigma_{2}\times\Xi_{2})$. As it is a submatrix of
$(\partial P/\partial \xi|(P,\xi\in\Sigma_{2}\times\Xi_{2})$, points
that do not satisfy the regularity condition of
th.~\ref{th:suff_reg_cond} for the $9$ state variables system do not
satify it for the $12$ state variables system. We can associate to
this \=o-system $12$ possible sets of linearizing outputs: first, we
take $\xi_{1}=z$, then $\xi_{2}\in\bar\Xi_{1}$ (2 possibilities) and
$\{\xi_{3},\xi_{4}\}\in\bar\Xi_{2}$ (6 possibilities).}

\section{Simulations using the simplified flat model}\label{sec:simplified}

In this section, we show simulations done with the flat approximation
of the model. These simulations are conducted with the classical set
of flat outputs described in section~\ref{subsubsec:classical-FO},
that is $x,y,z,\beta$.

\PATCH{As mentioned above, the flat approximation was obtained by
  neglecting the dependency of $C_x,C_y,C_z$ on
  $p,q,r,\delta_l,\delta_m,\delta_n$. Comparing with the full aircraft
  model considered in the next section, no significant change of
behaviour could be noticed. We show here that the model remains robust
to various perturbations in the expression of the forces.}

This tends to show that in many contexts the flat approximation is
quite sufficient.

Moreover, we show that the flat model allows a high flexibility in
trajectory planning and tracking.

All these suitable properties remain when one reactor is out of order.

\subsection{Theoretical setting for feed-back
  design}\label{subsec:simplified_feed-back}

The great advantage of flatness is that the flat motion planning makes
an open loop control immediately available. When a closed loop is
required, the feedback is designed from the difference between the
actual values of the flat outputs and their reference values, so that
this difference, being the solution of some differential equation,
tends to zero. \PATCH{The reader may refer to~\cite{Buccieri-et-al-2009} for
a more robust approach of the control of flat systems.}

In the framework of the flat aircraft model, the feedback is done is
two stages. Indeed the dependency of the system variables on $F,p,q,r$
has a slow dynamics in comparison to the rapidity of the dynamics that
controls $p,q,r$ from $\delta_l,\delta_m,\delta_n$. This allows to
construct a cascade feedback, as done in~\cite{Martin-96}. More
precisely, one can build a dynamic linearizing feedback that allows
controlling the partial state vector $\Xi =
(x,y,z,V,\alpha,\beta,\gamma,\chi,\mu,F)$ using the command $\dot
F,p,q,r$, which allows
following the reference trajectories of the flat outputs
$x,y,z,\beta$, using static linearizing feedback. More precisely,
one can compute a vector valued function $\Delta_0$ and a matrix
valued function $\Delta_1$, both depending on
$x,y,z,V,\alpha,\beta,\gamma,\chi,\mu, F$ such that:
$$
\left (\begin{array}{c}
x^{(3)} \\
y^{(3)} \\
z^{(3)} \\
\dot{\beta}
\end{array} \right) = \Delta_0 + \Delta_1 \left ( \begin{array}{c}
p \\
q \\
r \\
\dot{F}
\end{array} \right )
$$

At this stage, the variables $p,q,r,\dot{F}$ are seen as commands. In
order to make the system linear, one introduces a new vector valued
command $v$, such that: 
$$
\left ( \begin{array}{c}
p \\
q \\
r \\
\dot{F}
\end{array} \right ) = \Delta_1^{-1} (v - \Delta_0)
$$

Eventually the command $v$ is chosen of the form:
$$
v(t) = \left ( \begin{array}{c}
P_0 (x_{\hbox{\tiny ref}}(t) - x(t)) + P_1 (\dot{x}_{\hbox{\tiny ref}}(t) - \dot{x}(t)) + P_2
(\ddot{x}_{\hbox{\tiny ref}}(t) - \ddot{x}(t)) + P_3 x_{\hbox{\tiny ref}}^{(3)}(t) \\ 
P_0 (y_{\hbox{\tiny ref}}(t) - y(t)) + P_1 (\dot{y}_{\hbox{\tiny ref}}(t) - \dot{y}(t)) + P_2
(\ddot{y}_{\hbox{\tiny ref}}(t) - \ddot{y}(t)) + P_3 y_{\hbox{\tiny ref}}^{(3)}(t) \\ 
P_0 (z_{\hbox{\tiny ref}}(t) - z(t)) + P_1 (\dot{z}_{\hbox{\tiny ref}}(t) - \dot{z}(t)) + P_2
(\ddot{z}_{\hbox{\tiny ref}}(t) - \ddot{z}(t)) + P_3 z_{\hbox{\tiny ref}}^{(3)}(t) \\  
-k_1(\beta_{\hbox{\tiny ref}}(t) - \beta(t)) + \dot{\beta}_{\hbox{\tiny ref}}(t),
\end{array} \right ) 
$$ where $P(X) = P_0 + P_1X + P_2X^2 + P_X^3$ is actually the
following polynomial $P(X) = (X-k_1)^3$. Therefore the error function
$e_s(t) = s_{\hbox{\tiny ref}}(t) - s(t)$ satisfies the following differential
equations $P(e_s(t)) = 0$. In our experiments, $k_1 = -5$, so that
$e_s(t) \underset{t \rightarrow +\infty}{\longrightarrow} 0$, for each
value of $s$ in $x,y,z,\beta$. The value of $k_1$ has tuned manually and is similar to values found in~\cite{Martin-PHD}. 

In a second stage the variables $p,q,r$ are controlled through a
static linearizing feedback based on
$\delta_l,\delta_m,\delta_n$. This part of the system, as mentioned
above, is fast in comparison to the first part. More precisely, one
can compute the vector valued function $\Lambda_0$ and a matrix valued
function $\Lambda_1$, both depending on $V,\alpha,\beta,p,q,r$ such
that:
$$
\left (\begin{array}{c}
\dot{p} \\
\dot{q} \\
\dot{r} \\
\end{array} \right) = \Lambda_0 + \Lambda_1 \left ( \begin{array}{c}
\delta_l \\
\delta_m \\
\delta_n \\
\end{array} \right )
$$

Then as previously, one introduced a new command $w$ such that
$$
\left ( \begin{array}{c}
\delta_l \\
\delta_m \\
\delta_n \\
\end{array} \right ) = \Lambda_1^{-1}(w - \Lambda_0)
$$
and 
$$
w = \left (\begin{array}{c}
-k_2(p_{\hbox{\tiny ref}}(t) - p(t)) +  \dot{p}_{\hbox{\tiny ref}}(t) \\
-k_2(q_{\hbox{\tiny ref}}(t) - q(t)) +  \dot{q}_{\hbox{\tiny ref}}(t) \\
-k_2(r_{\hbox{\tiny ref}}(t) - r(t)) +  \dot{r}_{\hbox{\tiny ref}}(t) \\
\end{array} \right),
$$
where $k_2 = -15$ in our experiments. Therefore $s(t) - s_{\hbox{\tiny ref}}(t)
\underset{t \rightarrow +\infty}{\longrightarrow} 0$, for $s \in
\{p,q,r\}$.  

The rationale behind this cascade feedback is the following. The
variables $(x,y,z,V,\alpha,\beta,\gamma,\chi,\mu,F)$ are slowly
controlled through $\dot{F},p,q,r$. Once the required values of
$p,q,r$ are known, they are quickly reached through the control
performed with $\delta_l,\delta_m,\delta_n$. The respective values of
$k_1$ and $k_2$ reflect the disparity of speed between the two
dynamics. \PATCH{See e.g.~\cite{Kokotovic1986,Ren2016} and the
references therin for a singular perturbation
approach of such systems.}

\subsection{Conventions used in our simulations}\label{subsec:simu_python}

We now show a series of experiments that illustrate the strength of
the flat approximation to control the aircraft in various
situations. Those experiments were all performed with GTM extracted
from~\cite{Hueschen}. The implementation is made in Python, relying
on the symbolic library \textsc{sympy}, the numerical array library
\textsc{numpy} and the numerical integration of ODE systems from the
library \textsc{scipy}.

The experiments are all about following a reference trajectory defined
by the following expressions: 

$$
\left \{ \begin{array}{rcl}
x_{\hbox{\tiny ref}}(t) & = & V_1\cos(\pi(t-T_{\hbox{\tiny initial}})/(T_{\hbox{\tiny final}}-T_{\hbox{\tiny initial}})) \\
y_{\hbox{\tiny ref}}(t) & = & V_1\sin(\pi(t-T_{\hbox{\tiny initial}})/(T_{\hbox{\tiny final}}-T_{\hbox{\tiny initial}})) \\ 
z_{\hbox{\tiny ref}}(t) & = & -V_2t-1000 \\ 
\beta_{\hbox{\tiny ref}}(t) & = & 0
\end{array} \right.,
$$
where $T_{\hbox{\tiny initial}} = 0 s, T_{\hbox{\tiny final}} = 30 s, V_1 = 30 ms^{-1}, V_2 = 5
ms^{-1}$. This reference trajectory is an upward helix.

\subsection{Initial perturbation}\label{subsec:initial_perturbation}

We carried out experiments where the aircraft
started away from the reference trajectory and then joined it after a
few seconds. If the initial perturbation is not too big, the feedback
alone is capable to attract the aircraft to the reference
trajectory. If the initial starting point is really far away from the
reference trajectory, the flexibility of the flatness based approach
allows designing very easily a transition trajectory which can be
followed with the feedback and that brings the aircraft to the
reference upward helix trajectory.

\begin{figure}[h!]
\hbox to \hsize{\hss\renewcommand\tabcolsep{0pt}
  \lower1cm\hbox{\includegraphicsh[width=4.cm]{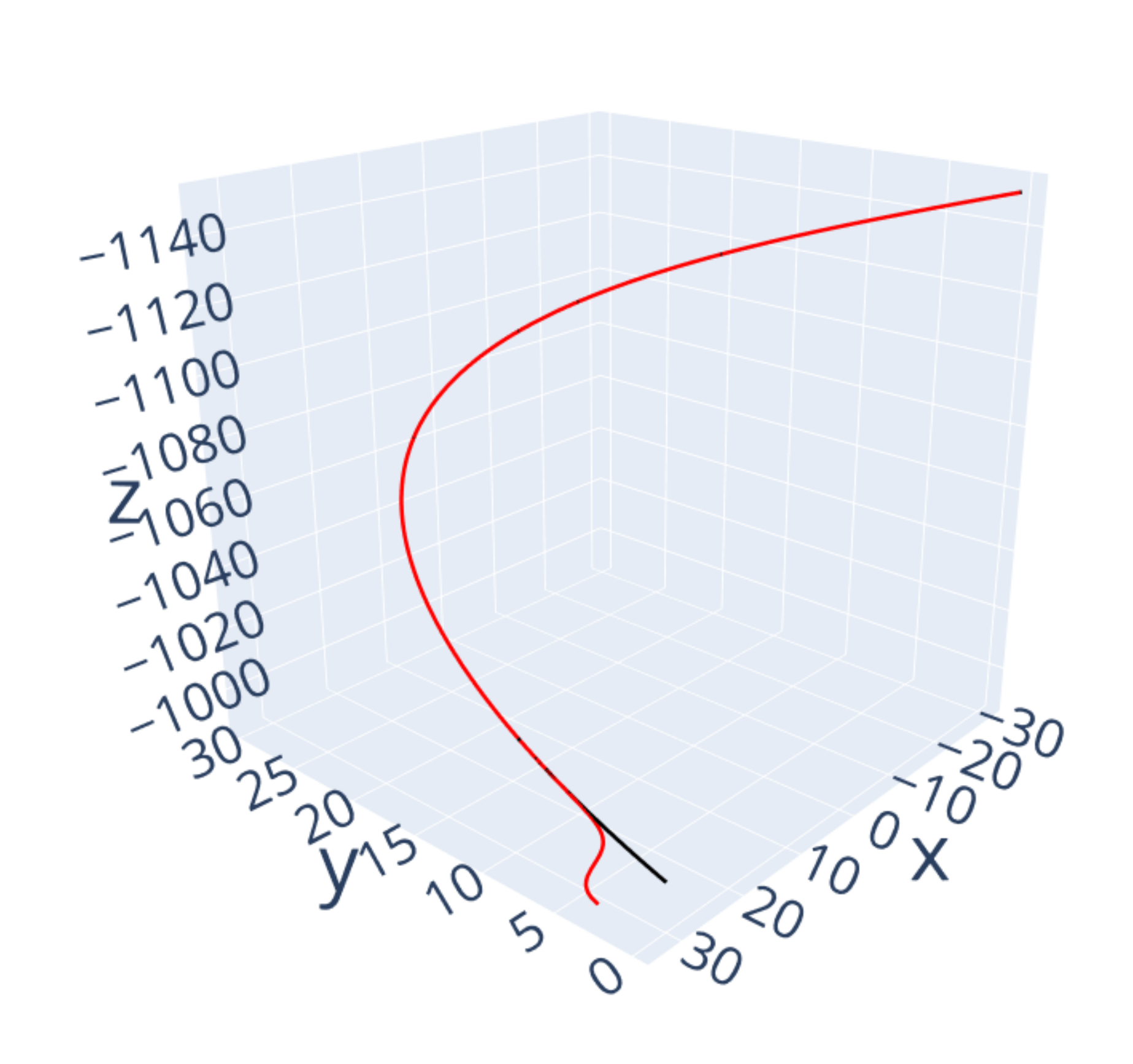}
    {a)~\begin{tabular}[b]{l}
        \hbox to 0.4cm{\color{red} \hrulefill}{\footnotesize simulated trajectory with initial perturbation} \\ 
        \hbox to 0.4cm{\color{blue} \hrulefill}{\footnotesize reference trajectory} \\
        \end{tabular}}
  \begin{tabular}[b]{cc}
  \graph[3.cm]{./OneReactorP10_X.pdf}{x\>(\m)}{b)}&
  \graph[3.cm]{./OneReactorP10_Y.pdf}{y\>(\m)}{c)}\\
  \graph[3.cm]{./OneReactorP10_Z.pdf}{z\>(\m)}{d)}&
  \graph[3.cm]{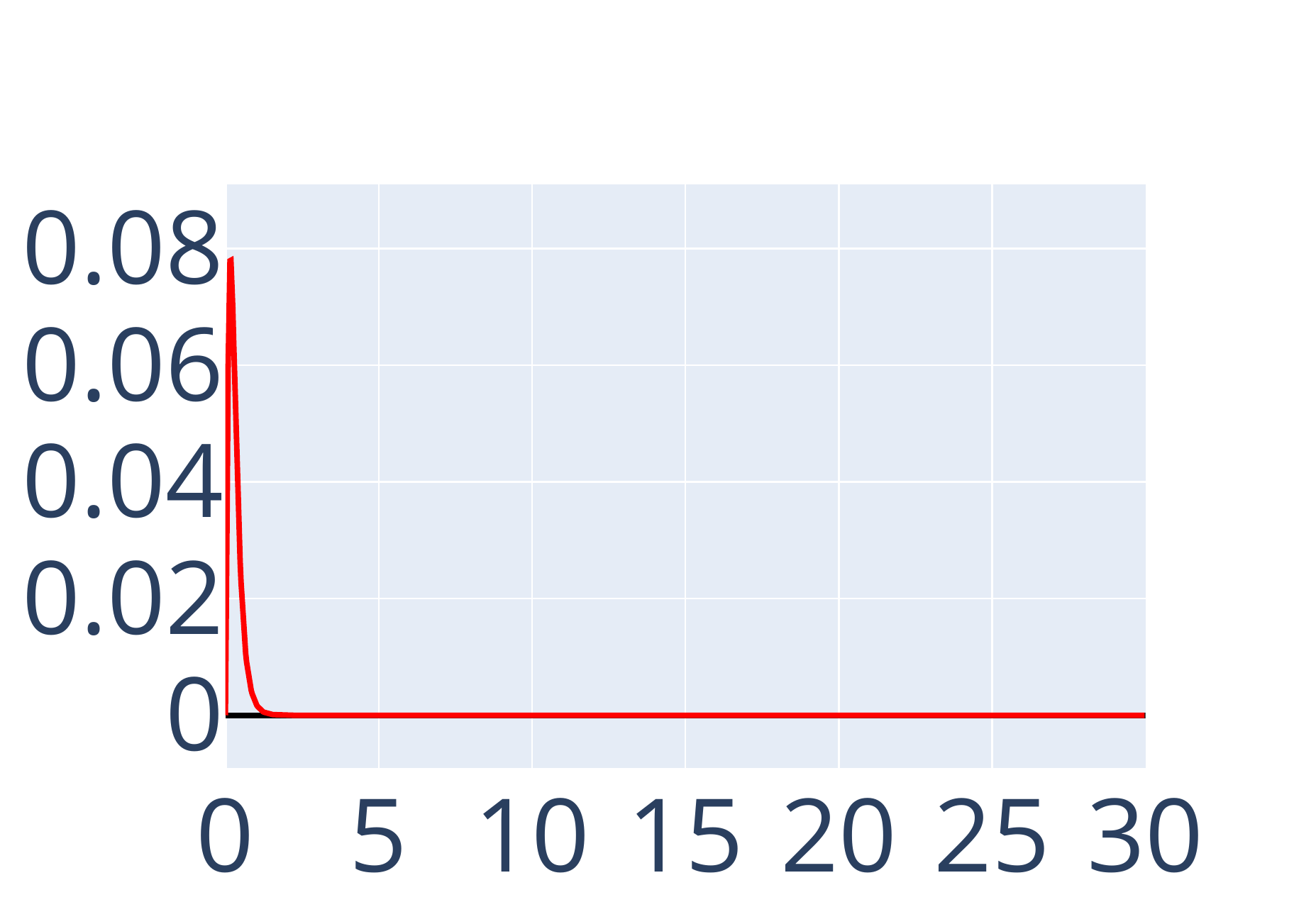}{\beta\>(\rad)}{e)}
  \end{tabular}
}\hss}
\caption{\label{fig::GTM_perturbation1_OneReactor}
The values of the GTM trajectory, one engine: the aircraft converges toward the reference trajectory. a) Trajectory, $3$D view;
  b)--e) Histories of $x,y,z,\beta$}
\end{figure}

The experiment is moreover performed when a one reactor is broken.  We
observe that the actual trajectory of the aircraft merges with the
reference one, as shown in
figure~\ref{fig::GTM_perturbation1_OneReactor}. The reference and the
actual trajectories merge perfectly, even when starting from a point
off the trajectory. \PATCH{In~\cite[fig.~2]{Ollivier2022}, the power
  of one engine is slowly reduced, so that one can see the variation
  of state functions, such as bank angle.}

\subsection{Variable wind}\label{subsec:wind}

In this section, we address the most critical problem about the flat
approximation. Since the dependency of the aerodynamic coefficient on
$p,q,r,\delta_l,\delta_m,\delta_n$ is discarded, one can wonder if the
model is robust enough to significant perturbations in the values of
the thrust. It turns out that under mild external forces, the model
remains reliable.

In the last experiment, the motion of the aircraft is perturbed by a
variable wind. This perturbation force is a sinusoidal function which
amplitude is 50lbf and frequency is 0.1Hz. This perturbation is
assumed to roughly model a quasi-periodic change in wind
direction. Doing a more realistic analysis is beyond our scope.

This setting is applied to the GTM with the reference trajectory
defined above. We observe a very robust behavior of the model, as
rendered in figures~\ref{fig::GTM_variable_wind1}
and~\ref{fig::GTM_variable_wind2}. The reference and the actual
trajectories still merge perfectly with a variable wind. For the
variables $F,\alpha,p,q,r$ variations due to the variable wind are
noticeable in the graphs.

\begin{figure}[h!]
  \hbox to \hsize{\hss
  \graph[3.5cm]{./TwoReactorsPWSinusoidalWind_X.pdf}{x\>(\m)}{a)}
  \graph[3.5cm]{./TwoReactorsPWSinusoidalWind_Y.pdf}{y\>(\m)}{b)}
  \graph[3.5cm]{./TwoReactorsPWSinusoidalWind_Z.pdf}{z\>(\m)}{c)}
%  \graph[3.5cm]{./TwoReactorsPWSinusoidalWind_beta.pdf}{\beta\>(\rad)}{d)}
\hss}
  \caption{\label{fig::GTM_variable_wind1}
    Variable wind. a)--d) Histories of $x,y,z$}
\end{figure}

\begin{figure}[h!]
\hbox to \hsize{\hss 
  \graph[3.cm]{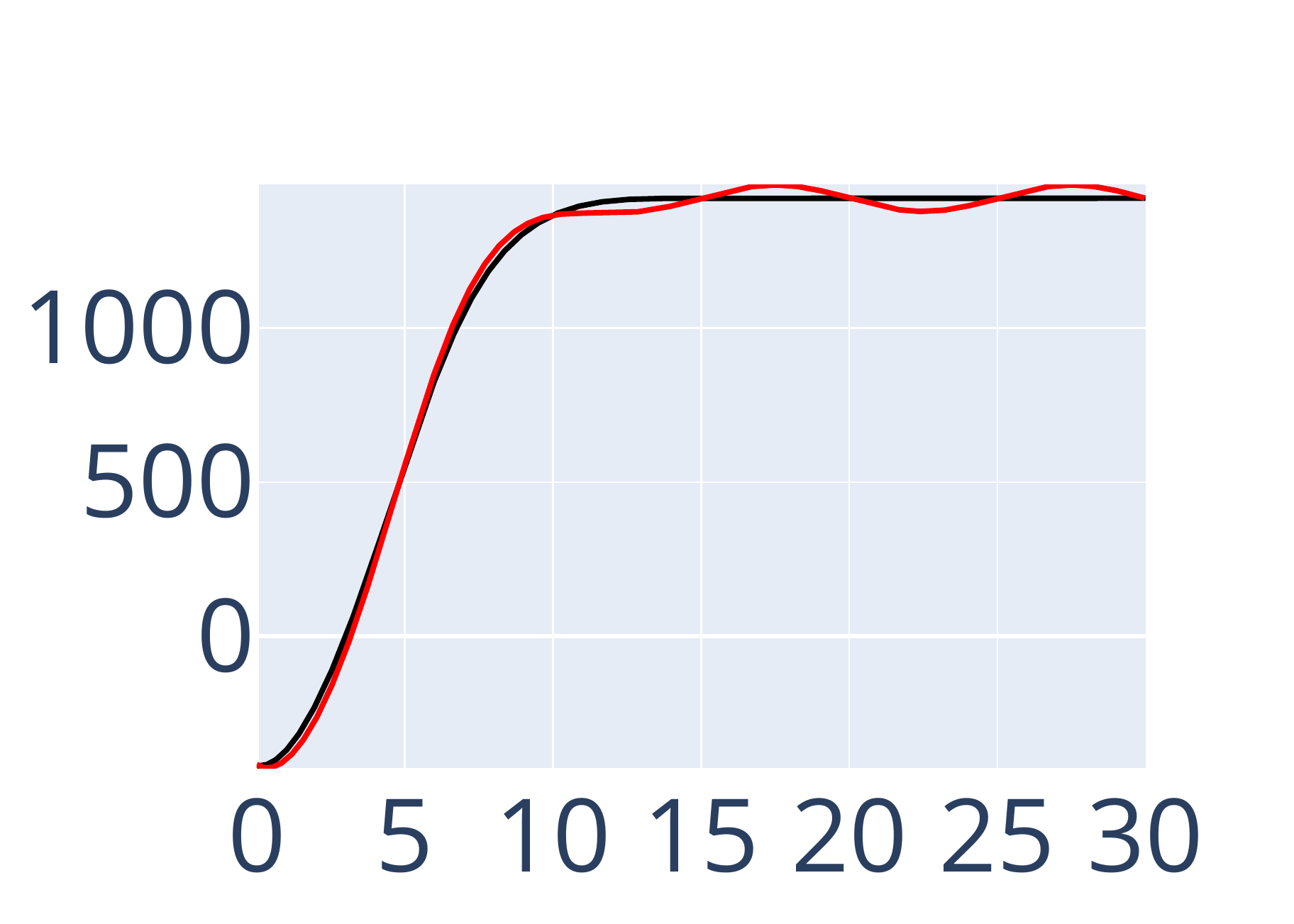}{F\>(\Newton)}{a)}\hfill
  \graph[3.cm]{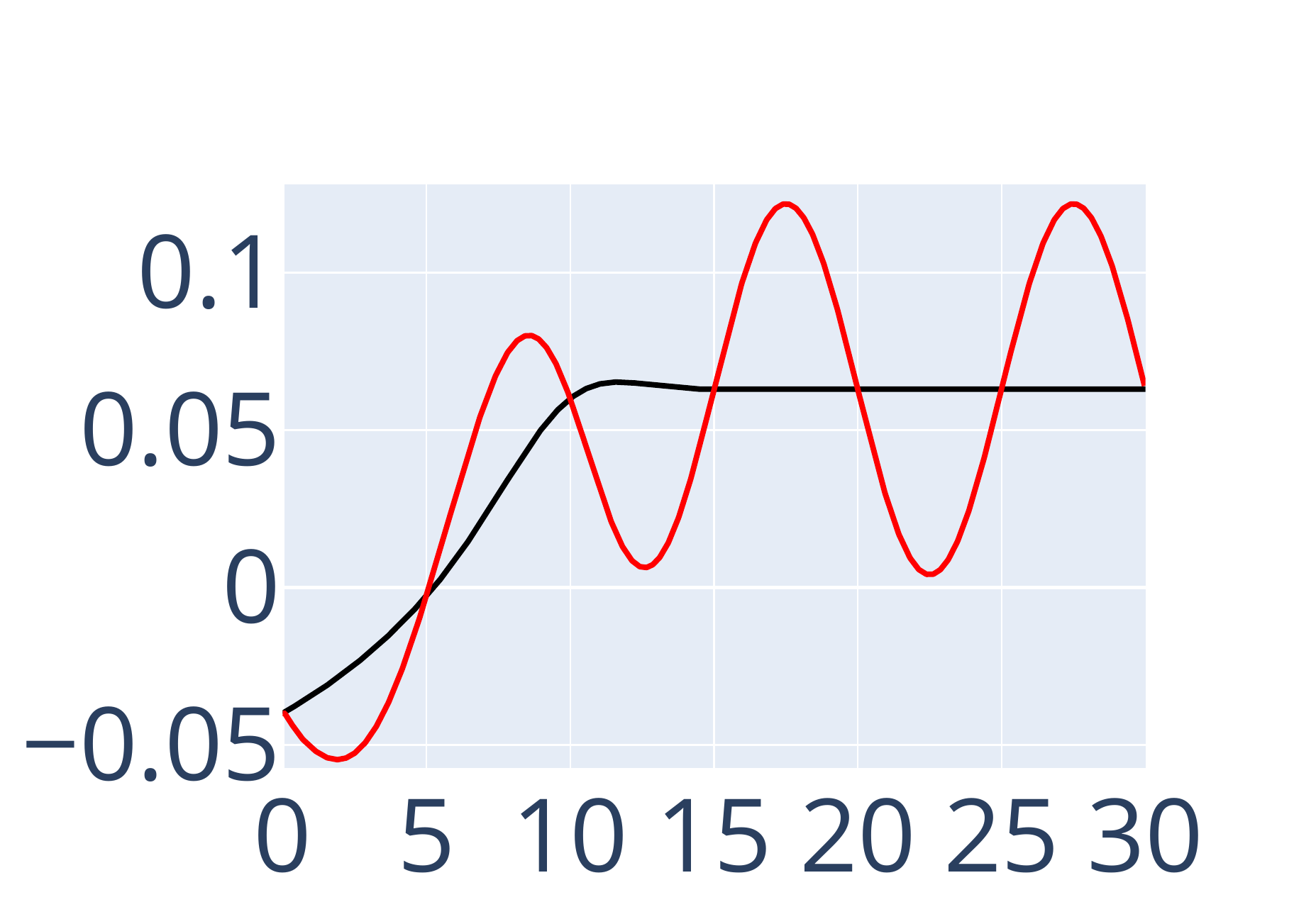}{\mu\>(\rad)}{b)}\hfill
  \graph[3.cm]{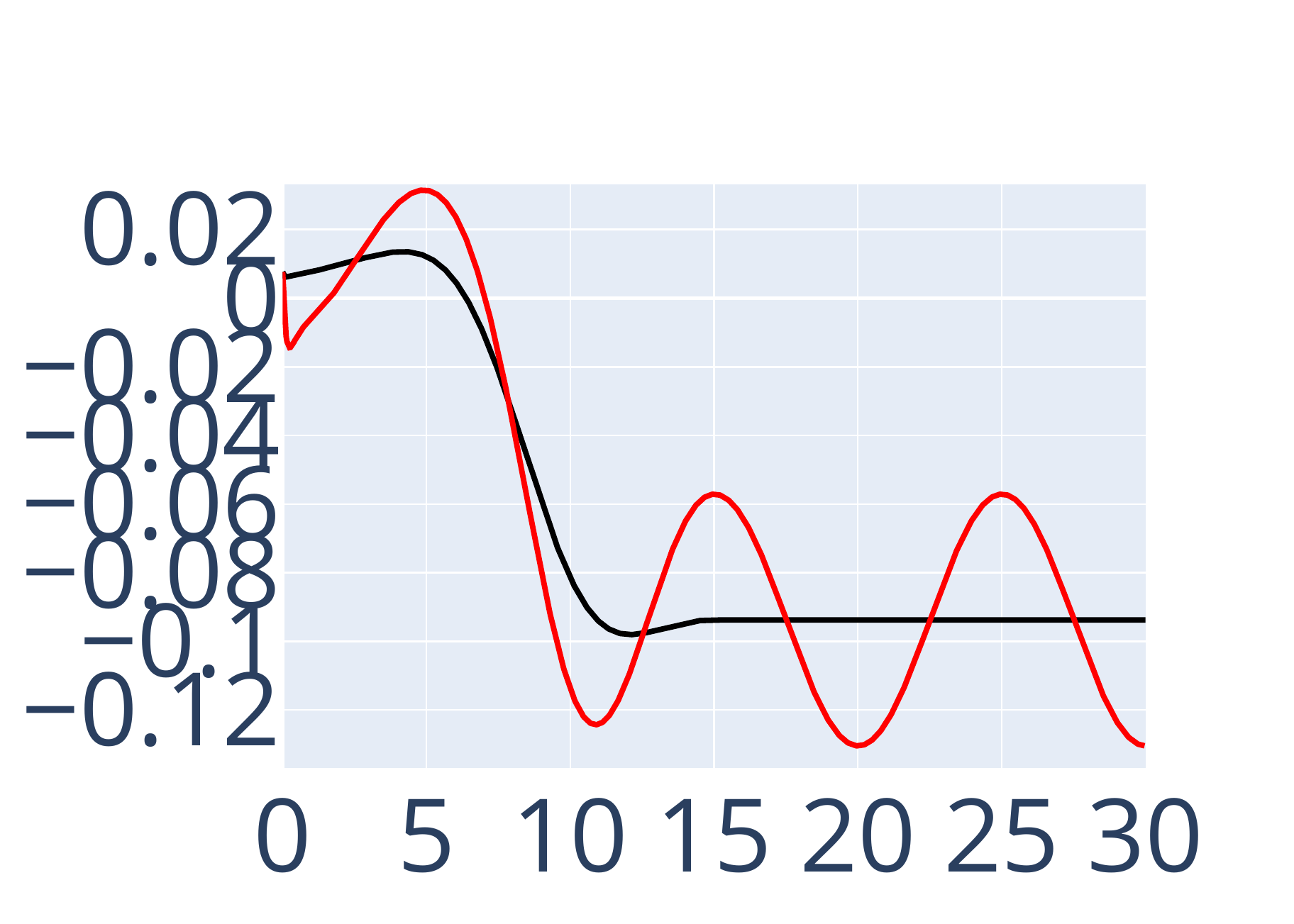}{p\>(\rad/\s)}{c)}\hfill
  \graph[3.cm]{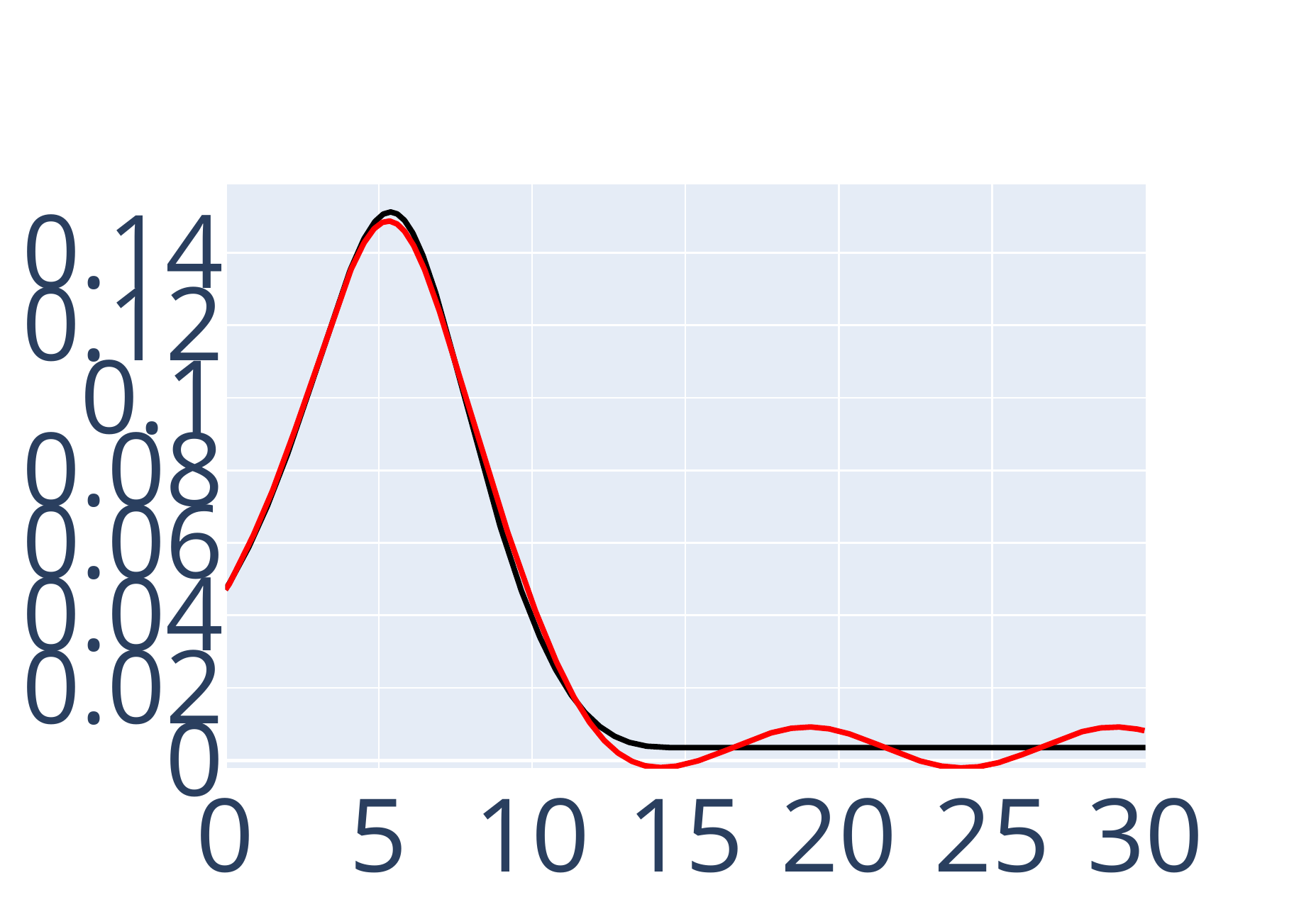}{q\>(\rad/\s)}{d)}\hfill
%  \graph[3.cm]{./TwoReactorsPWSinusoidalWind_r.pdf}{r\>(\rad/\s)}{e)}\hss
}
  \caption{\label{fig::GTM_variable_wind2}
Variable wind. a)--e) Histories of $F,\mu,p,q$}
\end{figure}

\section{Simulations using the full model}\label{sec:full_model}

\subsection{Implementation in Maple}\label{subsec:Maple}

Our implementation computes power series approximations of all state
variables and control at regular time interval. We proceed in the same
way for the feed-back, which is also approximated with power series at
the same time interval for better efficiency during numerical
integration.

The structure of the equations implies the existence of a lazy flat
parametrization, as shown by th.~\ref{th:reg_o_syst_flat}. Moreover,
it means that the requested power series can be computed in a fast way
using Newton's method, when initialized with suitable values. Most of
the time, values for state variables such as $\beta$ or $\alpha$ are
close to $0$. If not, calibration functions can give \textit{e.g.}
the value of $\alpha$, assuming that 
$V$, $\gamma'$, $\chi$, $\beta$ and $\mu$ are constants.
See~\cite[\S~4]{Ollivier2022} for more details.

We denote by $\hat \xi(t)$ the planed function for any state variable
$\xi$, according to the motion planning using the simplified flat system and
the choice of $\hat x$, $\hat y$, $\hat z$ and $\hat \beta$ (or $\hat
\mu$). We also denote by $\updelta\xi$ the difference $\xi-\hat\xi$
between the planed trajectory and the trajectory computed with the
full model.  We did not manage to get a working feed-back without
using integrals $I_{1}$, $I_{2}$, $I_{3}$, $I_{4}$ of
$\cos(\chi)\updelta x+\sin(\chi)\updelta y$,
$-\sin(\chi)\updelta x+\cos(\chi)\updelta y$,
$\updelta z$ and $\updelta \beta$ (or $\updelta \mu$)
respectively. When using the flat outputs $x,y,z,F$, $I_{4}$ is no
longer needed.

\PATCH{For the feed-back, we choose positive real numbers $\lambda_{i,j}$,
with $1\le i\le 4$, and $1\le j\le 5$ for $i=2$ or $i=3$ and $1\le
j\le 3$ for $i=1$ or $i=4$. The value of $\updelta F$, $\updelta
\delta_{l}$, $\updelta \delta_{m}$, $\updelta \delta_{n}$, are computed,
so that $\prod_{i=1}^{3} (d/dt + \lambda_{1,i})I_{1}=0$,
$\prod_{i=1}^{3} (d/dt + \lambda_{4,i})I_{4}=0$, $\prod_{i=1}^{5}
(d/dt + \lambda_{2,i})I_{2}=0$, $\prod_{i=1}^{5} (d/dt +
\lambda_{3,i})I_{3}=0$, using the derivation $d/dt$ of the linearized
simplified system around the planed trajectory. Then, we use the
controls $\hat \delta_{l}+\updelta \delta_{l}$, $\hat \delta_{m}+\updelta
\delta_{m}$, $\hat \delta_{n}+\updelta \delta_{n}$, $\hat F+\updelta F$ in
the numerical integration. If the $\updelta\xi$ are small enough to
behave like the $\dd\xi$ of the linearized system, and the solution of
the full model not too far from the planed solution of the simplified
model, the convergence is granted.}

In practice, the choice of suitable $\lambda_{i,j}$ is difficult and
empirical: two small, the trajectory is lost, two high, increasing
oscillations may appear. We neglect here the dynamics of the actuator,
our goal being to show that the feed-back is able to provide a
solution for the full model, using the trajectory planed with the
simplified one, the linearizing outputs remaining close to their
original values. \PATCH{The feedback used here relies on a linear
  approximation around the planned trajectory. I has the advantage of
  fast computation but will not work any more when perturbations are
  too great. One may refer e.g.\ to~\cite{Buccieri-et-al-2009}
  or~\cite{Zerar2009} for more robust control of flat systems.}
\bigskip

Unless otherwise stated, angles are expressed in radians, lengths in
meters, times in seconds, thrusts in Newtons, masses in kg. Curves in
red correspond to the planed trajectory $\hat\xi$, while curves in
blue correspond to the integration of the full model. For more clarity
$z$ has been replaced by $-z$ to get positive values when drawing
curves.

Computation times are given using Maple 19 with an
Intel processor Core~i5 $2.5$GHz. These are just
indications that can vary from a session to another.

\subsection{Using flat outputs with $\mathbf{\mu}$. Gravity-free flight with the
  F16}\label{subsec:parabolic-F16C}

\begin{figure}[h!]
  \hbox to \hsize{\hss
  \includegraphicsh[width=3.7cm]{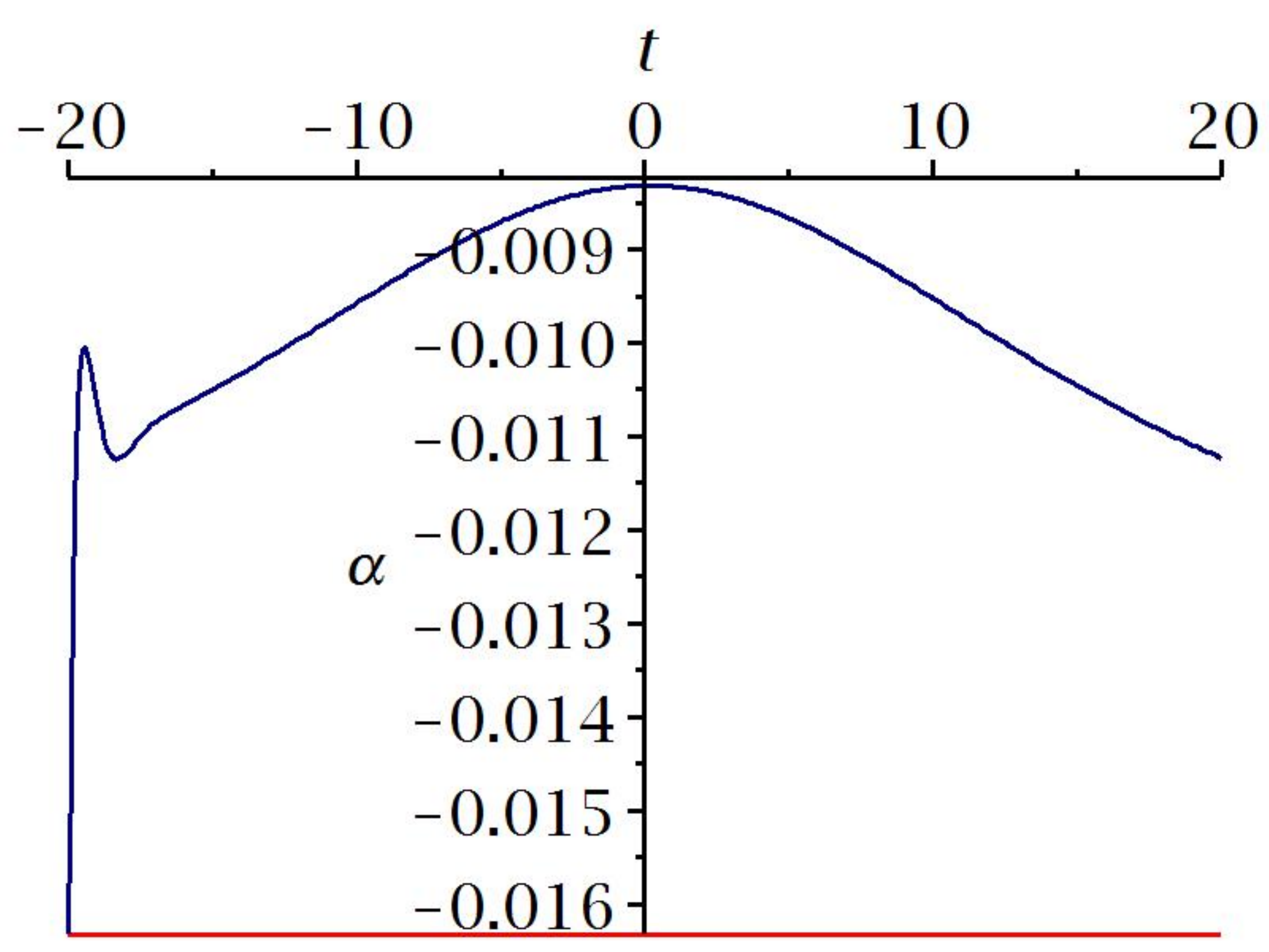}{a)}\hss
  \includegraphicsh[width=3.7cm]{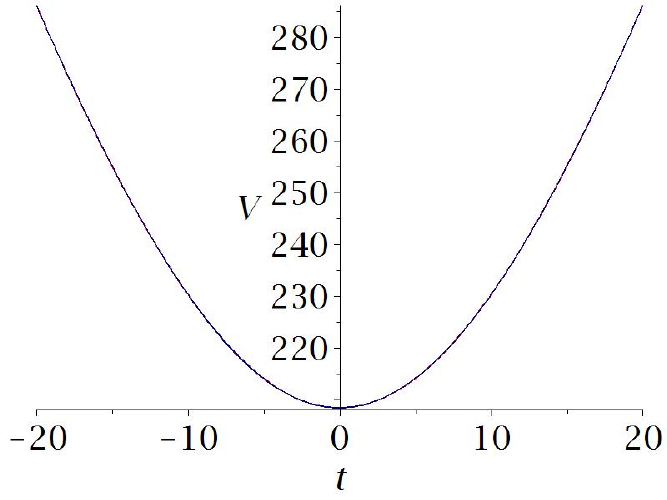}{b)}\hss
  \includegraphicsh[width=3.7cm]{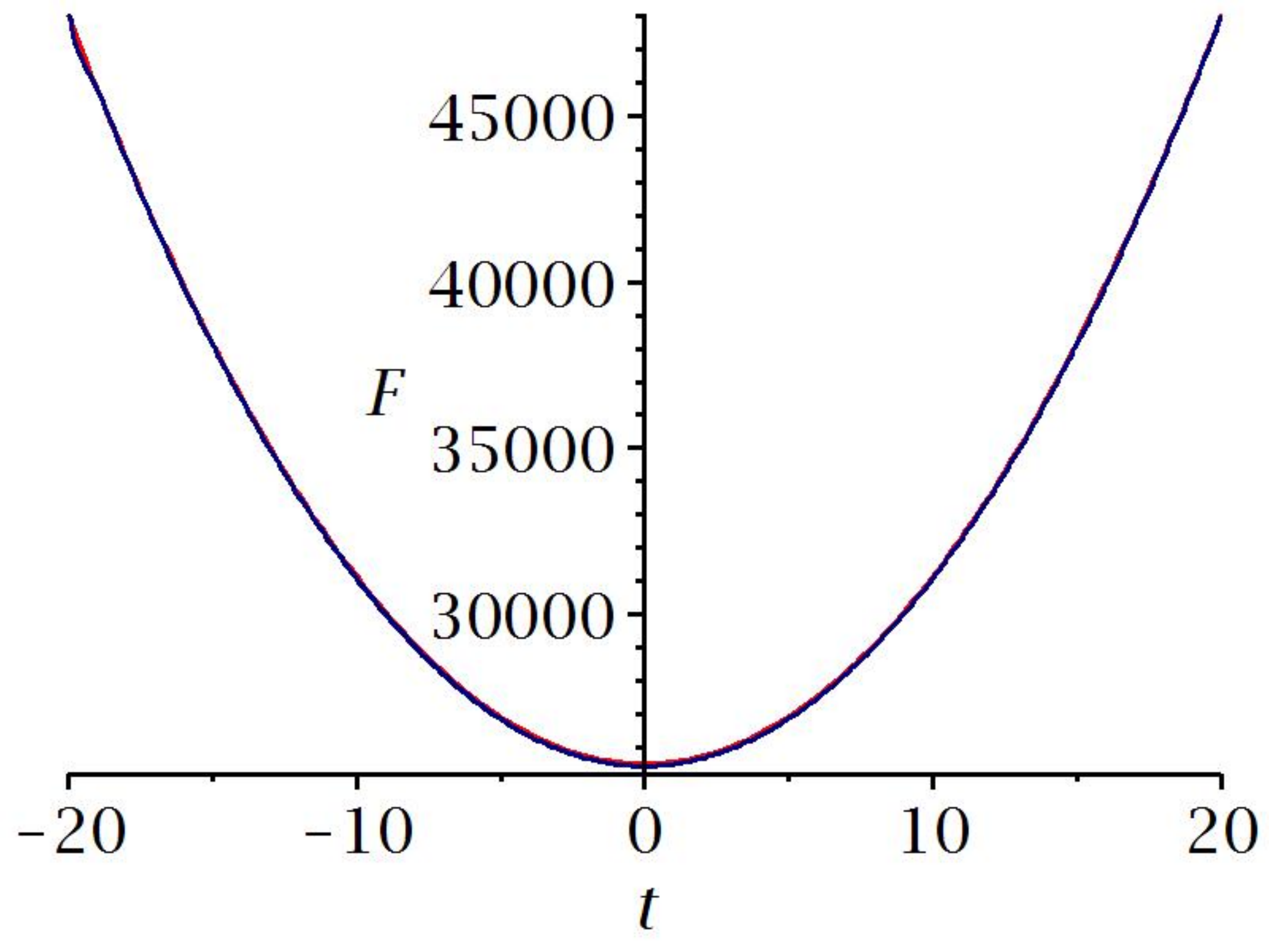}{c)}\hss}
  \hbox to \hsize{\hss
  \includegraphicsh[width=3.7cm]{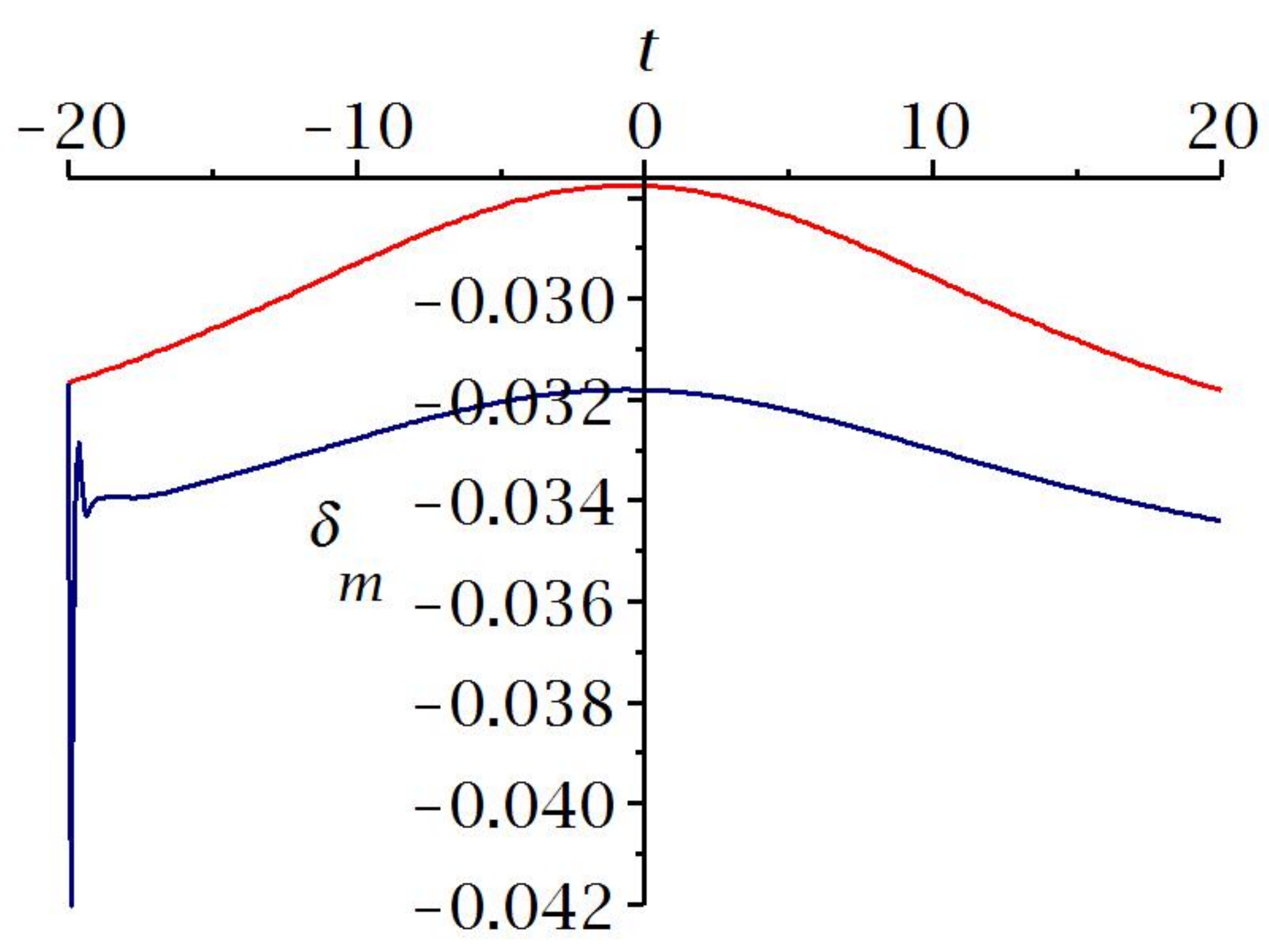}{d)}\hss
  \includegraphicsh[width=3.7cm]{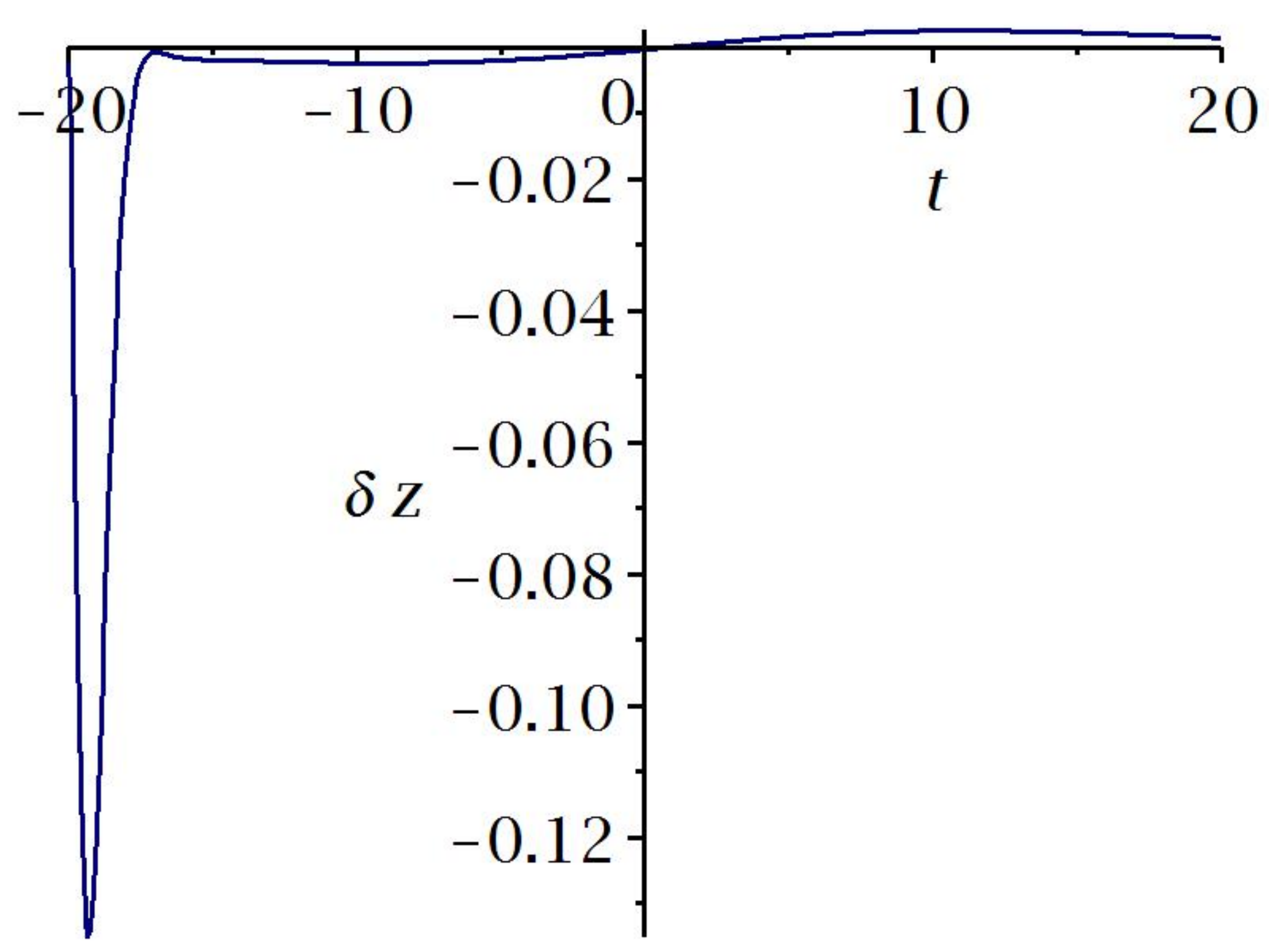}{e)}
\hss}
\caption{\label{fig::F16_parabolic}
  F-16 $0$-g flight. a)~$\alpha\>(\rad)$; b)~$V\>(\m/\s)$;
    c)~$F\>(\Newton)$; d)~$\delta_{m}\>(\rad)$ and e)~$\updelta z\>(\m)$.}
\end{figure}

We experiment here a $0$-g flight with a parabolic trajectory, using
the F-16 model for which $C_{z}$ vanishes for a value of $\alpha$
close to $-0.016$ rad. See fig.~\ref{fig::F16_parabolic}\PATCH{, where
  $\updelta z$ 
  corresponds to the difference between the reference altitude and the
  altitude of the simulated model.}
For this simulation, we used an expression of air density $\rho(z)$
varying with altitude $z$, following Martin~\cite[A.16
  p.~97]{Martin-PHD}.  The total computation time of the simulation is
$371$s.

\begin{equation}
\begin{array}{ll}
  x&=750t\kmh;\quad y=0;\quad
  z=g\frac{t^{2}}{2}-2000;\quad \mu=0;\\
  \lambda_{i,j}&=5.;\\
  \rho(z)&=1.225(1+0.0065z/288.15)^{(9.80665/(287.053\times0.0065)-1)}
\end{array}\label{eq::F16_0G}
\end{equation}

\section{Conclusion}

We have introduced a notion of oudephippical systems, or \=o-systems,
that generalizes many special cases of chained or triangular systems
previously know in the control literature. We further gave algorithms
to test in polynomial time if a given system belongs to this category
and to test regularity conditions that are sufficient to imply that
the system is flat at a some point.

Then, as the flat outputs are state variables, the computation of the
flat parametrization and the design of a feedback to stabilize the
system around the planned trajectory are computationally easy, even
for nontrivial systems, such as the aircraft model used as an
illustration. This model is flat if one neglects some terms and our
simulations show that a suitable feedback is able to compensate model
errors due to this simplification.

The systematic study of all possible flat outputs in our setting made
us discover new flat outputs for the aircraft, some of them of
practical interest. They provide a set of charts with flat
parametrization, that cover all common flight situations, flat
singularities being close to stalling conditions. \PATCH{We are
  moreover able to prove a sufficient condition for flat singularity
  that could be applied to the aircraft model.}

\PATCH{Simulations show that the mathematical framework can handle
  complex models. A systematic study of stabilization issues is left
  for further investigations.}

\section*{Acknowledgments}

The authors thank Jean L{\'e}vine for advices and inspiration.

\section*{Declarations}

\begin{itemize}
\item Funding. F.~Ollivier thanks the French ANR project
  ANR-22-CE48-0016 NODE (Numeric-symbolic resolution of differential
  equations) and the ANR project ANR-22-CE48-0008 OCCAM (Theory and
  practice of differential elimination).
  \item Competing
  interests. No foundings or relations between the authors and
  companies or entities having interests in space or aircraft
  industries or control devices or softwares is likely to call into
  question the objectivity of this work.
\item Code availability. The most recent implementation in Maple for
  aircraft motion planning, including generalized flatness, is
  available at \url{http://www.lix.polytechnique.fr/~ollivier/GFLAT/}.
\item Authors' contributions. Y.J.~Kaminski is responsible for Python
  implementations and F.~Ollivier for Maple
  implementations. Y.J. Kaminski and F. Ollivier have bath contributed
  to the study, conception and realization or to the writing and
  typesetting processes.
\item Data availability. Simulation data for Maple implementations are
  available with the Maple packages.
\end{itemize}

%=========================================
%=========================================
%=========================================
%=========================================

\bibliographystyle{amsplain}
\bibliography{./aircraft-jacobi.bib}

\end{document}